\documentclass[11pt]{amsart}
\usepackage{amsmath,amssymb,amsthm,mathtools}
\usepackage{booktabs,array}
\usepackage{enumitem}
\usepackage{hyperref}
\usepackage[margin=1.15in]{geometry}
\usepackage{tikz}
\usetikzlibrary{arrows.meta,calc,decorations.markings}
\usepackage{pgfplots}
\pgfplotsset{compat=1.18}

\newtheorem{theorem}{Theorem}[section]
\newtheorem{proposition}[theorem]{Proposition}
\newtheorem{lemma}[theorem]{Lemma}
\newtheorem{corollary}[theorem]{Corollary}
\theoremstyle{definition}
\newtheorem{definition}[theorem]{Definition}
\theoremstyle{remark}
\newtheorem{remark}[theorem]{Remark}

\newcommand{\R}{\mathbb R}
\newcommand{\Z}{\mathbb Z}
\newcommand{\Sph}{\mathbb S}
\newcommand{\grad}{\operatorname{grad}}
\newcommand{\Ric}{\operatorname{Ric}}
\newcommand{\tr}{\operatorname{trace}}
\newcommand{\Hess}{\operatorname{Hess}}

\title[Cohomogeneity one biconservative hypersurfaces]{Cohomogeneity one biconservative hypersurfaces in Euclidean spaces and spheres}
\author{S. Montaldo}
\address{Dipartimento di Matematica e Informatica, Universit\`a degli Studi di Cagliari, Via Ospedale 72, 09124 Cagliari, Italy}
\email{montaldo@unica.it}
\author{C. Oniciuc}
\address{Faculty of Mathematics, ``Al.I. Cuza'' University of Iasi, Bd. Carol I no. 11, 700506 Iasi, Romania}
\email{oniciucc@uaic.ro}
\author{A. Ratto}
\address{Dipartimento di Matematica e Informatica, Universit\`a degli Studi di Cagliari, Via Ospedale 72, 09124 Cagliari, Italy}
\email{rattoa@unica.it}
\author{A. Siffert}

\address{Mathematisches Institut, Universit\"at M\"unster, Einsteinstra\ss e 62, 48149 M\"unster, Germany}
\email{asiffert@uni-muenster.de}

\thanks{The authors S.M. and A.R. are members of the Italian National Group G.N.S.A.G.A. of INdAM.  The work was partially supported by the Project {PRoBIKI} funded by Fondazione di Sardegna.
}
\begin{document}
\begin{abstract}

{We study global properties of cohomogeneity one biconservative hypersurfaces in Euclidean spaces and spheres, invariant under the cohomogeneity two actions of Hsiang and Lawson. We show that, for a non-CMC invariant hypersurface, the biconservative equation reduces to a simple relation between the curvature of the profile curve in the orbit space and the volume of the principal orbits. As a consequence, the profiles are the geodesics of a metric on the orbit space conformal to the original one, the conformal factor being a power of the orbit volume. This gives a correspondence between proper biconservative and minimal invariant hypersurfaces when the orbit multiplicities are multiplied by three. In particular, using Hsiang's minimal hyperspheres we construct new proper biconservative embeddings $\Sph^{n-1}\to\Sph^{n}$ for $n=8,11,14,17$.

In Euclidean space, a compactification of the associated planar dynamical system, together with a Lyapunov function, yields for every non-rotational action a one-parameter family, up to homotheties, of complete proper biconservative hypersurfaces with two ends asymptotic in angle to the unique invariant minimal cone. Whether these ends oscillate around the cone is decided, for all the actions at once, by the sign of a single discriminant which depends only on the dimension and on the angle of the orbit space. A weighted identity shows that the distance from the origin along each profile has at most one critical point, which is a strict minimum.

For rotational hypersurfaces in spheres, we identify the reduced metric with one of Otsuki's metrics, for a value of his parameter which is determined by the dimension and is in general not an integer. Otsuki's monotonicity of the period then gives uniqueness for each admissible rotation index, and completes the classification of compact non-CMC rotational biconservative hypersurfaces in spheres of dimension at least four. The same formulation gives direct proofs of the previously known existence and non-embeddedness results.}
\end{abstract}

\maketitle

\section{Introduction}

Biconservative immersions arise by retaining the tangential part of the biharmonic equation. For an orientable hypersurface
\[
\varphi:M^{n-1}\longrightarrow N^n,
\qquad N^n=\R^n\ \text{or}\ \Sph^n,
\]
with unit normal field $\eta$, shape operator $A$, and mean curvature vector
\[
\mathbf H=\frac{f}{n-1}\eta,
\qquad f=\operatorname{trace}A,
\]
the biconservative condition (see Section~\ref{sec:biconservative}) is
\begin{equation}\label{eq:bc-general}
2A(\grad f)+f\,\grad f=0.
\end{equation}
{Throughout the paper we call $f$, the trace of the shape operator, the \emph{mean curvature function} of $\varphi$; it is the unnormalized mean curvature, the mean curvature being $f/(n-1)$.}
Thus, from \eqref{eq:bc-general}, every constant mean curvature hypersurface is biconservative, while the non-CMC solutions provide the genuinely new part of the theory. As in \cite{MOR-BC}, we call such solutions \emph{proper biconservative hypersurfaces}.

{ In the present paper we investigate how symmetry can be used to construct complete and compact non-CMC biconservative hypersurfaces and to determine their global geometry.  We extend two previous lines of investigation. In \cite{MOR-BC} we studied proper biconservative hypersurfaces invariant under the product action $SO(p+1)\times SO(q+1)$, together with the rotational case. That analysis produced a planar dynamical system for the profile curve and a sharp change from oscillatory to nonoscillatory behaviour. 

On the other hand, in \cite{MOR-JGP} we considered the full list of cohomogeneity two isometric actions on Euclidean space which give rise to cohomogeneity one hypersurfaces, and proved that the corresponding biharmonic hypersurfaces are necessarily minimal. The natural question is therefore whether the biconservative equation admits a unified treatment for all those actions and, if so, whether the global phenomena found in the product case persist. The same groups act with cohomogeneity two on the sphere $\Sph^n\subset\R^n\oplus\R$, trivially on the last factor, and we treat both the Euclidean and the spherical cases. 
The actions under consideration are the Hsiang--Lawson cohomogeneity two actions \cite{HL}. A cohomogeneity one $G$-invariant hypersurface is determined by a profile curve in the orbit space. This reduces the problem to a two-dimensional one, but there are two difficulties which must be handled simultaneously. First, the coefficients of the reduced equation become singular at the boundary rays of the orbit space. Secondly, information on trajectories in the phase plane has to be translated back into geometric statements concerning completeness, asymptotics, intersections with the minimal cone, and possible extension across singular orbits.

Our approach separates these two issues. The orbit-volume function gives a uniform reduction of the biconservative equation for every action. We then introduce polar variables and compactify the planar system so that the singular boundary rays become genuine invariant boundary components carrying hyperbolic saddle points.  Then, in the present biconservative problem we determine an explicit Lyapunov function which gives direct global control of the trajectories and replaces the Bendixson--Poincar\'e argument used in the minimal case.
}

\subsection*{Main results}

Throughout the paper $G$ denotes one of the closed subgroups of $O(n)$ listed in Table~\ref{tab:actions}, acting on $\R^n$ with cohomogeneity two, and on $\Sph^n\subset\R^n\oplus\R$ through its action on the first factor. The orbit space $Q=\R^n/G$ is isometric to the flat planar sector
\[
Q=\left\{(r\cos\vartheta,r\sin\vartheta):\ r\ge0,\ 0\le\vartheta\le\frac\pi d\right\},
\]
whose vertex corresponds to the origin of $\R^n$, whose two boundary rays $\vartheta=0$ and $\vartheta=\pi/d$ correspond to the singular orbits, and whose interior parametrizes the principal orbits. We denote by $V(\xi)$ the volume of the principal orbit over an interior point $\xi\in Q$; {with a normalization of the constant factor fixed once and for all after \eqref{eq:volume},} it is a homogeneous polynomial of degree $n-2$. A cohomogeneity one $G$-invariant hypersurface $\Sigma_\gamma$ of $\R^n$ is the preimage of a curve $\gamma$ in $Q$, its \emph{profile curve}, which we always parametrize by arc length; we write $\kappa$ for the signed curvature of $\gamma$ and $\nu$ for its unit normal in $Q$, with the orientation fixed in Section~\ref{sec:equivariant}. In the sphere, the orbit space $Q_{\Sph}=\Sph^n/G$ is the spherical lune over the same sector (Section~\ref{subsec:sphere-setting}), $g_Q$ is the round metric, $\kappa$ is the geodesic curvature, and $V$ is given by the same polynomial, restricted to $Q_{\Sph}\subset\R^3$.

\begin{table}[htbp]
\centering
\small
\renewcommand{\arraystretch}{1.15}
\begin{tabular}{@{}llll@{}}
\toprule
$G$ & $n$ & $d$ & multiplicities \\
\midrule
$SO(n-1)$ & $n\ge3$ & $1$ & $m_0=n-2$ \\
$SO(p+1)\times SO(q+1)$ & $n=p+q+2$, $p,q\ge1$ & $2$ & $m_0=q$, $m_1=p$ \\
$SO(3)$ & $5$ & $3$ & $m_0=m_1=m_2=1$ \\
$SU(3)$ & $8$ & $3$ & $m_i=2$ \\
$Sp(3)$ & $14$ & $3$ & $m_i=4$ \\
$F_4$ & $26$ & $3$ & $m_i=8$ \\
$SO(5)$ & $10$ & $4$ & $m_i=2$ \\
$SO(2)\times SO(m)$ & $n=2m\ge6$ & $4$ & $m_0=m_2=m-2$, $m_1=m_3=1$ \\
$S(U(2)\times U(m))$ & $n=4m\ge8$ & $4$ & $m_0=m_2=2m-3$, $m_1=m_3=2$ \\
$Sp(2)\times Sp(m)$ & $n=8m\ge16$ & $4$ & $m_0=m_2=4m-5$, $m_1=m_3=4$ \\
$U(5)$ & $20$ & $4$ & $m_0=m_2=5$, $m_1=m_3=4$ \\
$U(1)\times \mathrm{Spin}(10)$ & $32$ & $4$ & $m_0=m_2=9$, $m_1=m_3=6$ \\
$G_2$ & $14$ & $6$ & $m_i=2$ \\
$SO(4)$ & $8$ & $6$ & $m_i=1$ \\
\bottomrule
\end{tabular}
\caption{Cohomogeneity two actions on $\R^n$ used in \cite{MOR-JGP}. The representation-theoretic realization of each action is the one given there. The number $m_i$ is the multiplicity of the linear factor $w_{(d,i)}$ in the orbit volume $V$ (see \eqref{eq:wdi} and \eqref{eq:volume}), i.e., the exponent of $|\sin(i\pi/d-\vartheta)|$ in its angular part; it is also the multiplicity of the corresponding principal curvature of the invariant hypersurfaces.}
\label{tab:actions}
\end{table}

Our first result shows that the biconservative condition, which a priori involves all the principal curvatures of $\Sigma_\gamma$ (at most $d+1$ distinct ones) and their multiplicities, can be expressed through the single function $V$, and has a natural interpretation in the geometry of the orbit space.

\begin{theorem}[Biconservative profiles as conformal geodesics]\label{thm:A}
Let $N^n=\R^n$ or $\Sph^n$, and let $\Sigma_\gamma$ be a $G$-invariant hypersurface of $N^n$ whose profile curve $\gamma$ lies in the interior of the orbit space $Q$ (respectively $Q_{\Sph}$), and let $f$ be the trace of its shape operator. On every open subset of $\Sigma_\gamma$ on which $\grad f\neq0$, the hypersurface is biconservative if and only if
\[
3\kappa=\nu(\log V).
\]
Equivalently, $\gamma$, suitably reparametrized, is a geodesic of the conformal metric $V^{2/3}g_Q$ on the interior of the orbit space. Conversely, if the profile satisfies this equation everywhere, then $\Sigma_\gamma$ is biconservative, its mean curvature function is $f=-2\kappa$, and $\Sigma_\gamma$ is proper biconservative as soon as $\kappa$ is not constant. Finally, if $\Sigma_\gamma$ is {smooth,} biconservative and not CMC, then its profile satisfies $3\kappa=\nu(\log V)$ along the whole of $\gamma$.
\end{theorem}

Theorem~\ref{thm:A} should be compared with the classical result of Hsiang and Lawson \cite{HL}, according to which $\Sigma_\gamma$ is minimal if and only if $\gamma$ is a geodesic of $V^{2}g_Q$. Thus, from the point of view of the orbit space, passing from minimal to proper biconservative hypersurfaces amounts to replacing the exponent $2$ by $2/3$. In particular the group enters the problem only through $V$, and this is what makes a unified treatment of all the actions of Table~\ref{tab:actions}, and of both ambient spaces, possible. A first consequence is a correspondence between biconservative and minimal hypersurfaces: if the multiplicities of an action are three times those of another action with the same $d$, the two reduced metrics coincide (Corollary~\ref{cor:tripling}). This happens for $d=1$ when $n\equiv2\pmod 3$, for $d=2$ when $p$ and $q$ are divisible by three, and for $U(1)\times\mathrm{Spin}(10)$, whose profiles are those of minimal $S(U(2)\times U(3))$-invariant hypersurfaces.

\subsection*{Embedded hyperspheres in spheres}

The admissibility conditions on product-action profiles, including smoothness at the singular orbits, are conditions on the profile alone and do not depend on the dimensions of the orbit factors. Applying the correspondence to Hsiang's minimal hyperspheres~\cite{HsiangI,HsiangII} gives the following result, where the admissible profiles are defined in Section~\ref{subsec:transfer}.

\begin{theorem}[Biconservative hyperspheres]\label{thm:E}
Let $p_0,q_0\ge1$. Every admissible profile curve of a non-equatorial $SO(p_0+1)\times SO(q_0+1)$-invariant minimal embedding $\Sph^{p_0+q_0+1}\to\Sph^{p_0+q_0+2}$ is also the profile of an $SO(3p_0+1)\times SO(3q_0+1)$-invariant proper biconservative embedding $\Sph^{3p_0+3q_0+1}\to\Sph^{3p_0+3q_0+2}$. In particular there exist infinitely many proper biconservative embedded hyperspheres 
\begin{itemize}
\item $\Sph^7\to\Sph^8$, invariant under $SO(4)\times SO(4)$;
\item $\Sph^{10}\to\Sph^{11}$, invariant under $SO(4)\times SO(7)$;
\item $\Sph^{13}\to\Sph^{14}$, invariant under $SO(7)\times SO(7)$;
\item $\Sph^{13}\to\Sph^{14}$, invariant under $SO(4)\times SO(10)$;
\item $\Sph^{16}\to\Sph^{17}$, invariant under $SO(4)\times SO(13)$;
\item $\Sph^{16}\to\Sph^{17}$, invariant under $SO(7)\times SO(10)$.
\end{itemize}
\end{theorem}

The construction preserves the profile while changing the dimensions of the orbit factors. Smoothness and embeddedness are therefore transferred from the minimal examples, while the identity $f=-2\kappa$, together with the fact that $\kappa$ is not constant, shows that the resulting hyperspheres are non-CMC. 

\subsection*{The Euclidean space}
{ Let $G$ be one of the {groups} in Table~\ref{tab:actions}
and let $m_i$ be the positive integers which represent the multiplicities listed for the group $G$ in Table~\ref{tab:actions}.
We use polar coordinates $(r,\vartheta)$ on the orbit space $Q$, which is the sector $0\le\vartheta\le\pi/d$. The orbit volume factorizes as
\[
V(r,\vartheta)=r^{n-2}v(\vartheta),
\qquad
v(\vartheta)=\prod_{i=0}^{d-1}\left|\sin\left(\frac{i\pi}{d}-\vartheta\right)\right|^{m_i}.
\]
Thus, along each ray $\vartheta=\mathrm{const}$, $V$ is a power of $r$, while its angular part {$v$} is positive on $(0,\pi/d)$ and vanishes on the two boundary rays, which correspond to the singular orbits. Moreover $\log v$ is strictly concave (see {\eqref{eq:ell-second-new}}), so that $v$ has exactly one critical point $\vartheta_0\in(0,\pi/d)$, a strict maximum. 

The rays from the origin have $\kappa=0$, and along the ray $\vartheta=\mathrm{const}$ the unit normal is $\nu=r^{-1}\partial_\vartheta$, so that $\nu(\log V)=v'(\vartheta)/(r\,v(\vartheta))$. Hence a ray satisfies the minimal equation $\kappa=\nu(\log V)$, or the biconservative equation $3\kappa=\nu(\log V)$, if and only if $\vartheta=\vartheta_0$. The cone $C_G=\Sigma_{\{\vartheta=\vartheta_0\}}$ over this ray is therefore the unique $G$-invariant minimal cone, and the only $G$-invariant cone which is biconservative. It is invariant under homotheties, and it plays the role of an attractor for all the biconservative profiles: the results below show that every profile we construct goes to infinity in the direction of $C_G$, along one or both of its ends.} 

\begin{theorem}[Complete proper biconservative hypersurfaces]\label{thm:B}
Let $G$ be any of the groups in Table~\ref{tab:actions}, and fix $\vartheta_*\in(0,\pi/d)$ and $r_0>0$. Then there exists a unique $G$-invariant biconservative hypersurface $\Sigma$ of $\R^n$ whose profile curve $\gamma$ satisfies $3\kappa=\nu(\log V)$ and touches the circle $r=r_0$ of $Q$ at the point of polar angle $\vartheta_*$.

This hypersurface is complete and proper biconservative, and its profile has the following properties.
\begin{enumerate}[label=\textup{(\roman*)},leftmargin=2em]
\item The profile $\gamma$ is defined for every value $s\in\R$ of the arc length and never leaves the interior of $Q$. Hence $\Sigma$ is a proper immersion of $\R\times G/H$ into $\R^n$, where $G/H$ is the principal orbit.
\item The distance $r(s)$ from the origin is strictly decreasing for $s<0$ and strictly increasing for $s>0$. Thus $\Sigma$ lies outside the open ball of radius $r_0$ centred at the origin and touches its boundary sphere along exactly one principal orbit.
\item Both ends of $\gamma$ go to infinity, and their direction converges to that of the minimal cone:
\begin{equation}\label{eq:profile-asymptotics-new}
\lim_{s\to\pm\infty}r(s)=+\infty,
\qquad
\lim_{s\to\pm\infty}\vartheta(s)=\vartheta_0 .
\end{equation}
More precisely, $r(s)/|s|\to1$ as $s\to\pm\infty$.
\item The mean curvature of $\Sigma$ is not constant: it is different from zero along the orbit of minimal distance from the origin and tends to zero along both ends; more precisely, $f(s)=o(1/|s|)$ as $s\to\pm\infty$.
\end{enumerate}
For $d\ge2$, different values of $\vartheta_*$ correspond, up to homotheties, to different profile curves, so that we obtain a family of complete proper biconservative hypersurfaces depending on one parameter. For $d=1$ all these hypersurfaces are congruent modulo homotheties{, translations along the axis of rotation and the reflection in a hyperplane orthogonal to it}.
\end{theorem}

The key point in the proof of Theorem~\ref{thm:B} is a strict Lyapunov function for the reduced problem. We emphasize that completeness is not read off from the phase portrait: we recover explicitly the original arc-length parameter, in the compactified time of the phase portrait, and prove that it tends to $\pm\infty$ along the two ends; completeness of the induced metric then follows from a general criterion (Lemma~\ref{lem:completeness}).

To describe all solutions at once, we use the invariance of the equation under homotheties to eliminate the distance from the origin. Writing the profile in polar coordinates $(r(s),\vartheta(s))$, let $\delta(s)$ be the oriented angle from the outward radial direction to the unit tangent vector of $\gamma$. Then the pair $(\vartheta,\delta)$ satisfies an autonomous planar system
which, after multiplication by the {factor} $h(\vartheta)=\sin(d\vartheta)/d${, which is positive in the interior of $Q$}, becomes
\begin{equation}\label{eq:compactified1}
\left \{ 
\begin{aligned}
\vartheta'&=3h(\vartheta)\sin\delta,\\
\delta'&=a(\vartheta)\cos\delta-(n+1)h(\vartheta)\sin\delta,
\end{aligned}
\right .
\end{equation}
where $a(\vartheta)=h(\vartheta)\,(\log v(\vartheta))'$.
The coefficients of \eqref{eq:compactified1} now extend smoothly to the closed cylinder
\[
\mathcal C=\left[0,\frac\pi d\right]\times(\R/2\pi\Z),
\]
whose two boundary circles $\vartheta=0$ and $\vartheta=\pi/d$ are invariant and encode the singular orbits. {We call \emph{compactified biconservative vector field} the vector field
\[
X=3h(\vartheta)\sin\delta\,\partial_\vartheta+\bigl(a(\vartheta)\cos\delta-(n+1)h(\vartheta)\sin\delta\bigr)\partial_\delta
\]
on $\mathcal C$.  In the interior of $\mathcal C$ the trajectories of $X$ describe, modulo homotheties, the {profile curves of} biconservative hypersurfaces.}

\begin{theorem}[The compactified phase portrait]\label{thm:C}
The compactified biconservative vector field $X$ on $\mathcal C$ has the following global structure, illustrated in Figure~\ref{fig:compactified-phase}.
\begin{enumerate}[label=\textup{(\roman*)},leftmargin=2em]
\item It has exactly six equilibria: the two interior points
\[
P_0=(\vartheta_0,0),\qquad P_1=(\vartheta_0,\pi),
\]
which both represent the minimal cone $C_G$, traversed away from and towards the origin respectively, and the four boundary points
\[
Q_1=\left(0,\frac\pi2\right),\quad
Q_2=\left(0,\frac{3\pi}{2}\right),\quad
Q_3=\left(\frac\pi d,\frac\pi2\right),\quad
Q_4=\left(\frac\pi d,\frac{3\pi}{2}\right),
\]
which correspond to profiles meeting a singular orbit orthogonally.
\item $P_0$ is a hyperbolic sink and $P_1$ is a hyperbolic source, while the four points $Q_i$ are hyperbolic saddles.
\item Every trajectory through a point of the {segment} $\delta=\pi/2$ in the interior of $\mathcal C$ meets this {segment} only once, tends to $P_1$ in backward time and to $P_0$ in forward time; these are exactly the trajectories of the profiles of Theorem~\ref{thm:B}. More generally, every trajectory meeting the open strip $|\delta|<\pi/2$ tends to $P_0$ in forward time, every trajectory meeting the open strip $\pi/2<\delta<3\pi/2$ tends to $P_1$ in backward time, and the interior of $\mathcal C$ contains no periodic orbit.
\item Each saddle $Q_i$ has exactly one separatrix in the interior of $\mathcal C$, and these four separatrices connect the saddles to the minimal cone:
\[
Q_1\longrightarrow P_0,\qquad Q_4\longrightarrow P_0,\qquad P_1\longrightarrow Q_2,\qquad P_1\longrightarrow Q_3 .
\]
The corresponding profiles leave a singular orbit orthogonally, after a finite arc length and at positive distance from the origin, and their other end is asymptotic in angle to the minimal cone; along each of them the distance from the origin is strictly monotone. For $d\ge2$ the corresponding hypersurfaces{, which lie over the principal orbits,} are embedded and proper biconservative, and their mean curvature has a nonzero limit at the singular orbit{; adding this singular orbit one obtains their closure, which is a complete embedded hypersurface whenever it is smooth across the singular orbit}, while for $d=1$ they are hyperplanes {with one point removed, whose closures are hyperplanes}.
\end{enumerate}
\end{theorem}

Theorems~\ref{thm:B} and~\ref{thm:C} say that all the profiles we construct approach the minimal cone at infinity. The way in which they approach it is decided by the linearization of the vector field $X$ at $P_0$ and $P_1$, which a priori depends on the whole multiplicity pattern of the action. It turns out that, at the critical point $\vartheta_0$, the second derivative of the logarithm of the angular part of $V$ always equals $-d(n-2)$ (Lemma~\ref{lem:ell2-universal-new}), so that the answer depends only on $d$ and on the dimension $n$.

\begin{theorem}[Oscillation around the minimal cone]\label{thm:D}
Let us define the following discriminant
\[
\Delta_d=(n+1)^2-12d(n-2).
\]
\begin{enumerate}[label=\textup{(\roman*)},leftmargin=2em]
\item If $\Delta_d<0$, then $P_0$ is a spiral sink and $P_1$ is a spiral source. Consequently each end of every profile of Theorem~\ref{thm:B}, and the end at infinity of every separatrix profile of Theorem~\ref{thm:C}, crosses the ray of the minimal cone infinitely many times; that is, the hypersurface meets the minimal cone $C_G$ along infinitely many principal orbits, which go off to infinity.
\item If $\Delta_d\ge0$, then $P_0$ and $P_1$ are nodes (improper nodes when $\Delta_d=0$), and each of the above ends crosses the minimal cone only finitely many times: eventually it stays on one side of $C_G$.
\end{enumerate}
In particular, for $d=1$ the second alternative always occurs; for $d=2$, writing $n=p+q+2$, the profiles oscillate if and only if $p+q\le17$; all the actions with $d=3$ and $d=6$ in Table~\ref{tab:actions} are in the oscillating regime; for $d=4$, where all the admissible dimensions are even, the profiles oscillate if and only if $n\le42$.
\end{theorem}

For $d=2$ the threshold of Theorem~\ref{thm:D} is the one found in \cite{MOR-BC}, so that the single discriminant $\Delta_d$ replaces a separate linearization for every representation. For the product action, Theorems~\ref{thm:B} and~\ref{thm:D} also correct part of the statement of \cite[Theorem~4.12]{MOR-BC}, where the hypersurfaces below were claimed to be incomplete cones with a singularity at the origin. More precisely, we have

\begin{corollary}[The case $SO(p+1)\times SO(q+1)$]\label{cor:d2}
Let $p,q\ge1$ and $n=p+q+2$. There exists an infinite family of complete proper $SO(p+1)\times SO(q+1)$-invariant biconservative hypersurfaces in $\R^n$. Their profile curves are defined on the whole real line, stay at positive distance from the origin, and both their ends are asymptotic in angle to the profile of the minimal cone $y=\sqrt{q/p}\,x$. If $p+q\le17$ each end intersects the minimal cone infinitely many times, whereas if $p+q\ge18$ each end eventually does not intersect it. Moreover, along each of these profiles both coordinate functions $x$ and $y$ are bounded below by explicit positive constants.
\end{corollary}

\subsection*{Rotational hypersurfaces in spheres}

The embedded hyperspheres of Theorem~\ref{thm:E} arise from the correspondence with minimal hypersurfaces. A second application concerns rotational hypersurfaces.

For $d=1$, $V^{2/3}g_{\Sph^2}$ is invariant under the rotations of $\Sph^2$ about a pole. This additional symmetry, the spherical analogue of the translations along the axis in the Euclidean case, provides a Clairaut first integral, and the problem becomes integrable. Let $\Phi(E)$ denote the angle by which a profile with Clairaut constant $E$ turns around the pole during one period of its curvature (see \eqref{eq:rot-Phi}).

\begin{theorem}[Closed rotational examples]\label{thm:F}
Let $n\ge3$. The profile curves of the proper biconservative $SO(n-1)$-invariant hypersurfaces of $\Sph^n$ are $\sigma$-elastic curves with $\sigma=(n-2)/(n+1)$; {each of them extends to a maximal profile, defined for every value of the arc length, whose} curvature is periodic and nonconstant, and {which stays} in an annulus around the pole without reaching the singular orbits. The period function $\Phi$ is strictly increasing in $E$, from $\pi$ to $\sqrt2\,\pi$. Consequently, for every pair of coprime natural numbers $l,r$ with $r<2l<\sqrt2\,r$ there is, up to isometries, exactly one compact proper biconservative immersion $\Sph^1\times\Sph^{n-2}\to\Sph^n$ whose profile turns $l$ times around the pole and has $r$ lobes; up to finite coverings, these are all the compact proper biconservative $SO(n-1)$-invariant hypersurfaces of $\Sph^n$, and none of them is embedded. 
\end{theorem}

The existence and non-embeddedness assertions recover \cite{MP-CAG} for $n=3$ and \cite{MOP-JMAA} in higher dimensions. The new conclusion for $n\ge4$ is uniqueness for each admissible pair $(l,r)$, together with exhaustion of the compact proper rotational examples up to finite coverings. It follows by identifying the reduced metric with Otsuki's $O^2_{\mathsf n}$, with real parameter $\mathsf n=(n+1)/3$, and applying his period-monotonicity theorem \cite{Otsuki93}; see Theorem~\ref{thm:Phi-monotone}.

The paper is organized as follows. 
Sections~\ref{sec:biconservative}--\ref{sec:bc-reduction} establish the geometric setting, the reduction to conformal geodesics, and the weighted identity. Sections~\ref{sec:compactified}--\ref{sec:focus-node} treat the Euclidean flow and prove completeness and the oscillation criterion. Section~\ref{sec:d2} provides a detailed analysis of the  product action and gives an independent completeness argument. Section~\ref{sec:spheres} proves the transfer theorem for embedded hyperspheres in $\Sph^n$ and the rotational classification.

\section{Biconservative hypersurfaces}\label{sec:biconservative}

In this section we recall the variational origin of biconservative hypersurfaces and the form of the biconservative equation in space forms; we refer to \cite{FetcuOniciuc} for a recent survey.

\subsection{Harmonic and biharmonic maps}

Let $\varphi:(M,g)\to(N,h)$ be a smooth map between Riemannian manifolds. Harmonic maps are the critical points of the energy functional
\begin{equation}\label{eq:energy}
E(\varphi)=\frac12\int_M|d\varphi|^2\,dv_g ,
\end{equation}
and the corresponding Euler--Lagrange equation is the vanishing of the tension field $\tau(\varphi)=\tr\nabla d\varphi$ \cite{EellsSampson}. Following a suggestion of Eells and Lemaire \cite{EellsLemaire}, a natural generalization is given by the \emph{biharmonic maps}, which are the critical points of the bienergy functional
\begin{equation}\label{eq:bienergy}
E_2(\varphi)=\frac12\int_M|\tau(\varphi)|^2\,dv_g .
\end{equation}
Jiang \cite{Jiang86} showed that the Euler--Lagrange equation of $E_2$ is $\tau_2(\varphi)=0$, where the \emph{bitension field} is
\begin{equation}\label{eq:bitension}
\tau_2(\varphi)=-\Delta\tau(\varphi)-\tr R^N\bigl(d\varphi,\tau(\varphi)\bigr)d\varphi .
\end{equation}
Here $\Delta=-\sum_{i}\bigl(\nabla^\varphi_{e_i}\nabla^\varphi_{e_i}-\nabla^\varphi_{\nabla^M_{e_i}e_i}\bigr)$ is the rough Laplacian on sections of $\varphi^{-1}TN$, computed in a local orthonormal frame $\{e_i\}$ of $M$, and $R^N(X,Y)=\nabla_X\nabla_Y-\nabla_Y\nabla_X-\nabla_{[X,Y]}$. Harmonic maps are clearly biharmonic. An isometric immersion is called a \emph{biharmonic submanifold} if it is a biharmonic map; since for an isometric immersion of an $m$-dimensional manifold $\tau(\varphi)=m\mathbf H$, where $\mathbf H$ is the mean curvature vector field, minimal submanifolds are biharmonic. When the target is the Euclidean space, an isometric immersion is biharmonic if and
only if its position vector satisfies $\Delta^2\varphi=0$, and Chen's conjecture
predicts that every biharmonic submanifold of the Euclidean space is minimal \cite{Chen91,ChenBook}. The conjecture is still open in general, although it has been proved, among other cases, for hypersurfaces with at most two distinct principal curvatures \cite{Dimitric}, for hypersurfaces of $\R^4$ \cite{Defever,HasanisVlachos}, for hypersurfaces with at most three distinct principal curvatures \cite{Fu15}, and for the $G$-invariant hypersurfaces considered in this paper \cite{MOR-JGP}. Very recently, Fu, Hong and Yang posted on arXiv a proof of Chen's conjecture for hypersurfaces in any dimension \cite{FHY} and for surfaces in $\R^4$ with flat normal bundle \cite{FHY1}{; an independent proof for hypersurfaces was also given by Wang \cite{Wang}}.

\subsection{The stress-energy tensor and biconservative hypersurfaces}

According to Hilbert \cite{Hilbert}, the stress-energy tensor associated with a variational problem is a symmetric $2$-covariant tensor $S$ which is conservative, i.e., $\operatorname{div}S=0$, at critical points. For the energy \eqref{eq:energy} it was studied by Baird and Eells \cite{BairdEells} (see also \cite{BairdRatto,Sanini}): the tensor
\[
S=\frac12|d\varphi|^2g-\varphi^*h
\]
satisfies $\operatorname{div}S=-\langle\tau(\varphi),d\varphi\rangle$, so that $\operatorname{div}S=0$ when $\varphi$ is harmonic. For an isometric immersion, $\tau(\varphi)$ is normal and $\operatorname{div}S=0$ holds automatically.

The stress-energy tensor of the bienergy was introduced by Jiang \cite{Jiang87} and further studied in \cite{LMO}. It is given by
\[
S_2(X,Y)=\frac12|\tau(\varphi)|^2\langle X,Y\rangle+\langle d\varphi,\nabla\tau(\varphi)\rangle\langle X,Y\rangle-\langle d\varphi(X),\nabla_Y\tau(\varphi)\rangle-\langle d\varphi(Y),\nabla_X\tau(\varphi)\rangle ,
\]
and it satisfies
\begin{equation}\label{eq:divS2}
\operatorname{div}S_2=-\langle\tau_2(\varphi),d\varphi\rangle .
\end{equation}
Thus the condition $\operatorname{div}S_2=0$ is the vanishing of the tangential part of the bitension field, and in general it is weaker than biharmonicity.

\begin{definition}\label{def:biconservative}
An isometric immersion $\varphi:M\to N$ is \emph{biconservative} if $\operatorname{div}S_2=0$, i.e., if $\tau_2(\varphi)^\top=0$.
\end{definition}

For hypersurfaces, the normal and tangential parts of $\tau_2$ were computed, with contributions from \cite{ChenBook,Oniciuc02,LoubeauMontaldo,Ou,BMO}, as follows.

\begin{theorem}\label{thm:bitension-decomposition}
Let $\varphi:M^{n-1}\to N^n$ be an isometric immersion with unit normal $\eta$, shape operator $A$, and mean curvature vector field $\mathbf H=\frac{f}{n-1}\eta$, $f=\tr A$. Then $\varphi$ is biharmonic if and only if
\begin{align}
&\Delta f+f|A|^2-f\Ric^N(\eta,\eta)=0,\label{eq:bih-normal}\\
&2A(\grad f)+f\grad f-2f\bigl(\Ric^N(\eta)\bigr)^\top=0.\label{eq:bih-tangent}
\end{align}
In particular, $\varphi$ is biconservative if and only if \eqref{eq:bih-tangent} holds.
\end{theorem}

\subsection{Biconservative hypersurfaces in space forms}

When $N^n$ is a space form, in particular $\R^n$ or $\Sph^n$, the Ricci operator is a multiple of the identity, so that $(\Ric^N(\eta))^\top=0$ and \eqref{eq:bih-tangent} reduces to the equation \eqref{eq:bc-general} of the Introduction,
\[
2A(\grad f)+f\grad f=0 .
\]
This equation has two immediate consequences.

First, every hypersurface with constant mean curvature (CMC) is biconservative. The interesting objects are therefore the non-CMC biconservative hypersurfaces, which we call \emph{proper}.

Secondly, on the open set where $\grad f\neq0$, the vector field $\grad f$ is a principal direction with principal curvature $-f/2$:
\begin{equation}\label{eq:gradf-principal}
A(\grad f)=-\frac f2\,\grad f .
\end{equation}
If $k_1=-f/2$ denotes this principal curvature and $k_2,\dots,k_{n-1}$ are the remaining ones, then $f=k_1+\dots+k_{n-1}$ gives the linear relation
\begin{equation}\label{eq:3k1}
3k_1+k_2+\dots+k_{n-1}=0
\end{equation}
on every open set on which $\grad f$ does not vanish. For the invariant hypersurfaces studied in this paper, $\grad f$ is tangent to the profile direction, and \eqref{eq:3k1} becomes the reduced equation \eqref{eq:bc-k}.

Let us briefly recall some known results. In three-dimensional space forms, proper biconservative surfaces are rotational, and they are locally classified \cite{HasanisVlachos,CMOP} {(see \cite{FuLorentz} for the analogous classification in Lorentzian $3$-space forms)}; their global properties (completeness) were studied in \cite{NistorOniciuc}, and closed examples exist in $\Sph^3$ but not in $\R^3$ \cite{MP-CAG}. In higher dimensions, rotational proper biconservative hypersurfaces of $\R^n$ were described in \cite{MOR-BC}, and closed ones of $\Sph^n$ in \cite{MOP-JMAA}. Turgay \cite{Turgay} classified the biconservative hypersurfaces of $\R^n$ with three distinct principal curvatures, under the assumption that the two principal curvatures different from the one in the direction of $\grad f$ have multiplicities greater than one, showing that they are either $SO(p+1)\times SO(q+1)$-invariant or generalized cylinders over rotational hypersurfaces; in this respect, the actions of Table~\ref{tab:actions} with $d\ge3$ provide the natural next class, with $d+1\ge4$ distinct principal curvatures.  The local structure of proper biconservative hypersurfaces of Euclidean spaces was also studied by Bibi, Soret and Ville \cite{BSV}. Rigidity results for compact biconservative hypersurfaces have been studied in \cite{FLO,AndronicKayhan}.

\section{The equivariant setting}\label{sec:equivariant}

\subsection{Euclidean space}\label{subsec:euclid-setting}

In order to make the main results precise, we first recall the orbit-space description used in \cite{MOR-JGP}. The advantage of this framework is that all the actions under consideration can be treated with the same notation: the representation-theoretic information enters only through the angle of the Weyl chamber and the multiplicities of the linear factors in the orbit-volume function.

We use the notation of \cite{MOR-JGP}. Let $G$ be a cohomogeneity two group of isometries acting on $\R^n$. The orbit space can be identified with the cone
\begin{equation}\label{eq:Q}
Q=\left\{(x,y)\in\R^2:\ y\geq 0,\quad x\sin\frac{\pi}{d}-y\cos\frac{\pi}{d}\geq 0\right\},
\end{equation}
where $d\in\{1,2,3,4,6\}$. The orbital distance metric is the Euclidean metric
\[
g_Q=dx^2+dy^2.
\]
For $i=0,\ldots,d-1$, set
\begin{equation}\label{eq:wdi}
w_{(d,i)}(x,y)=x\sin\frac{i\pi}{d}-y\cos\frac{i\pi}{d}.
\end{equation}
If $\xi=(x,y)$ lies in the interior of $Q$, the volume $V(\xi)$ of the principal orbit over $\xi$ is, up to a positive multiplicative constant, determined by
\begin{equation}\label{eq:volume}
V^2(x,y)=\prod_{i=0}^{d-1}\bigl[w_{(d,i)}(x,y)\bigr]^{2m_i},
\end{equation}
where the positive integers $m_i$ are the multiplicities listed in Table~\ref{tab:actions}. {Multiplying $V$ by a positive constant does not change the equations of the profile curves; we fix this constant once and for all by setting
\[
V=\prod_{i=0}^{d-1}\bigl|w_{(d,i)}\bigr|^{m_i}
\]
on the interior of $Q$, so that the actual volume of the orbit differs from $V$ by a positive factor depending only on $G$. Throughout the paper $V$ denotes this normalized function.}

\begin{remark}\label{rem:scope}
Table~\ref{tab:actions} lists, up to orbit equivalence, the cohomogeneity two actions of closed subgroups $G\subset O(n)$, that is, of isometry groups fixing a point; they correspond to the isotropy representations of the symmetric spaces of rank two \cite{HL}. Throughout the paper ``$G$-invariant'' refers to these linear actions. For $d=1$ we list only the multiplicity $m_0=n-2$ of the unique linear factor $w_{(1,0)}=-y$ of $V$. Cohomogeneity two groups of Euclidean isometries without fixed points, such as $SO(k)\times\R^{n-k-1}$ acting on $\R^k\times\R\times\R^{n-k-1}$, produce generalized cylinders and are not considered here.
\end{remark}

Let
\[
\gamma(s)=(x(s),y(s))
\]
be a curve in the interior of $Q$, parametrized by arc length,
\begin{equation}\label{eq:arclength}
\dot x^2+\dot y^2=1.
\end{equation}
We identify $Q$ with the closed Weyl chamber in the plane $\Pi_0\subset\R^n$, which meets every orbit orthogonally; all the points of the interior of $Q$ have the same isotropy group $H$, and $G/H$ is the principal orbit. The corresponding $G$-invariant hypersurface is the immersion
\[
F_\gamma:I\times G/H\to\R^n,
\qquad
F_\gamma(s,gH)=g\cdot\gamma(s),
\]
where $I$ is the interval of definition of $\gamma$; equivalently, $I\times G/H$ is identified with the pullback $\{(s,X)\in I\times\R^n:\ \pi(X)=\gamma(s)\}$ and $F_\gamma(s,X)=X$. We denote it by $\Sigma_\gamma$; its image is $\pi^{-1}(\gamma)$. If $\gamma$ is not injective, $\Sigma_\gamma$ is an immersed, not embedded, hypersurface (see Remark~\ref{rem:self-intersections}). 

We choose the orientation so that the unit normal projects to
\begin{equation}\label{eq:normalQ}
\nu=-\dot y\,\partial_x+\dot x\,\partial_y.
\end{equation}
As in \cite{MOR-JGP}, {the principal curvatures of $\Sigma_\gamma$ are given by the following $d+1$ functions, which need not be distinct at every point, so that $\Sigma_\gamma$ has at most $d+1$ distinct principal curvatures.} The first $d$ are
\begin{equation}\label{eq:ki}
k_i=-\nu\bigl(\log|w_{(d,i)}|\bigr)
=\frac{w_{(d,i)}(\dot y,-\dot x)}{w_{(d,i)}(x,y)},
\qquad i=0,\ldots,d-1,
\end{equation}
with multiplicity $m_i$, while
\begin{equation}\label{eq:kd}
k_d=\ddot y\,\dot x-\ddot x\,\dot y
\end{equation}
has multiplicity one. Thus, with the convention $f=\operatorname{trace}A$,
\begin{equation}\label{eq:f}
f=k_d+\sum_{i=0}^{d-1}m_i k_i.
\end{equation}
Since $V$ is given by \eqref{eq:volume}, formula \eqref{eq:ki} gives the useful identity
\begin{equation}\label{eq:orbitmean}
\sum_{i=0}^{d-1}m_i k_i=-\nu(\log V).
\end{equation}

\subsection{Spheres}\label{subsec:sphere-setting}

Let $G$ be one of the groups of Table~\ref{tab:actions}, and let it act on
\[
 \Sph^n\subset\R^n\oplus\R
\]
by $g\cdot(X,t)=(gX,t)$. Since $G$ acts on $\R^n$ with cohomogeneity two, it acts on $\Sph^n$ with cohomogeneity two as well. Let $\Pi_0\subset\R^n$ be the plane of the Weyl chamber $Q$, with coordinates $(x,y)$; then the three-space $\Pi_0\oplus\R$, with coordinates $(x,y,z)$, meets every orbit orthogonally, and the orbit space is the spherical lune
\begin{equation}\label{eq:QS}
 Q_{\Sph}=\bigl\{(x,y,z)\in\Sph^2:\ (x,y)\in Q\bigr\},
\end{equation}
with the round metric $g_{\Sph^2}$. In the spherical coordinates
\begin{equation}\label{eq:sphere-polar}
 x=\sin r\cos\vartheta,\qquad y=\sin r\sin\vartheta,\qquad z=\cos r,
 \qquad 0\le r\le\pi,\quad 0\le\vartheta\le\frac\pi d,
\end{equation}
the two sides $\vartheta=0$ and $\vartheta=\pi/d$ correspond to singular orbits, and the two poles $r=0,\pi$ are fixed points of $G$. The orbit over $(x,y,z)$ is $(G\cdot(x,y))\times\{z\}$, a copy of the orbit of $(x,y)$ in $\R^n$; hence its volume is still given by \eqref{eq:volume},
\[
 V(x,y,z)=\prod_{i=0}^{d-1}\bigl|w_{(d,i)}(x,y)\bigr|^{m_i}=\sin^{n-2}r\;v(\vartheta),
\]
with
\[
v(\vartheta)=\prod_{i=0}^{d-1}
\left|\sin\left(\frac{i\pi}{d}-\vartheta\right)\right|^{m_i}.
\]
For $d=2$ and $d=1$ we recover, respectively, the action of $SO(p+1)\times SO(q+1)$ on $\R^{p+1}\oplus\R^{q+1}\oplus\R$ and the rotational action of $SO(n-1)$ on $\R\oplus\R^{n-1}\oplus\R$.

Let $\gamma(s)=(x(s),y(s),z(s))$ be a curve in the interior of $Q_{\Sph}$, parametrized by arc length, $T=\gamma'$, and let $\nu=\gamma\times T$ be its unit normal in $\Sph^2$; let $\kappa$ be its geodesic curvature, $\nabla^{\Sph^2}_TT=\kappa\nu$. Let $\Sigma_\gamma$ be the corresponding $G$-invariant hypersurface of $\Sph^n$, and let $\eta$ be the unit normal field which at the points of $\Pi_0\oplus\R$ coincides with $\nu$ and is extended to $\Sigma_\gamma$ by equivariance, $\eta(g\cdot P)=g\cdot\eta(P)$.

\begin{lemma}\label{lem:sphere-curvatures}
With the convention $A(U)=-\nabla^{\Sph^n}_U\eta$, the hypersurface $\Sigma_\gamma$ has the principal curvature $\kappa$, with multiplicity one and principal direction $T$, and the principal curvatures
\[
 k_i=-\nu\bigl(\log|w_{(d,i)}|\bigr)=-\frac{w_{(d,i)}(\nu_1,\nu_2)}{w_{(d,i)}(x,y)},
 \qquad i=0,\dots,d-1,
\]
with multiplicity $m_i$, where $\nu=(\nu_1,\nu_2,\nu_3)$ and $w_{(d,i)}$ is regarded as a linear function on $\R^3$ independent of $z$. In particular \eqref{eq:f} and \eqref{eq:orbitmean} hold verbatim, with $k_d=\kappa$.
\end{lemma}

\begin{proof}
Since $\eta$ is tangent to $\Sph^n$ and orthogonal to $\Sigma_\gamma$, $\nabla^{\Sph^n}_U\eta=D_U\eta$ for tangent $U$, where $D$ is the flat connection of $\R^{n+1}$. Along $\gamma$, which lies in the totally geodesic $\Sph^2=\Sph^n\cap(\Pi_0\oplus\R)$, we have $\eta=\nu$; differentiating $\nu=\gamma\times T$ and using $D_TT=-\gamma+\kappa\nu$ gives $\nu'=\kappa\,\gamma\times\nu=-\kappa T$, so that $A(T)=\kappa T$. At a point $P=(\xi,z)$ of $\gamma$, write $\nu=(\nu_\xi,\nu_3)$ with $\nu_\xi=(\nu_1,\nu_2)\in\Pi_0$. For a vector $U$ tangent to the orbit, the last component of $\eta$ is constant along the orbit, so $D_U\eta=(D_U\widehat\nu_\xi,0)$, where $\widehat\nu_\xi$ is the equivariant extension of the horizontal vector $\nu_\xi$ along the orbit $G\cdot\xi\subset\R^n$. Hence $A(U)$ is the shape operator of the orbit $G\cdot\xi$ of $\R^n$ in the direction $\nu_\xi$, which, as for the Euclidean hypersurfaces of Section~\ref{subsec:euclid-setting} (see \eqref{eq:ki} and \cite{HL,MOR-JGP}), has the eigenvalues $-dw_{(d,i)}(\nu_\xi)/w_{(d,i)}(\xi)$ with multiplicities $m_i$. Since $w_{(d,i)}$ is linear and does not depend on $z$, these are the numbers $k_i$ of the statement.
\end{proof}

\section{Reduction of the biconservative equation}\label{sec:bc-reduction}

By Lemma~\ref{lem:sphere-curvatures}, the reduction is the same in $\R^n$ and $\Sph^n$, with $k_d=\kappa$. The point is that the gradient of the mean curvature is tangent to the profile direction.

For a $G$-invariant hypersurface, $f$ depends only on the arc-length parameter $s$. Hence $\grad f=\dot f\,\partial_s$. Since $\partial_s$ is a principal direction with principal curvature $k_d$, equation \eqref{eq:bc-general} becomes
\begin{equation}\label{eq:bc-factor}
\dot f\,(f+2k_d)=0.
\end{equation}
This is precisely the tangential equation which appears in the biharmonic system of \cite{MOR-JGP}.

On every open set on which $\dot f\neq0$, proper biconservativity is therefore equivalent to
\begin{equation}\label{eq:bc-k}
3k_d+\sum_{i=0}^{d-1}m_i k_i=0.
\end{equation}
Combining \eqref{eq:orbitmean} and \eqref{eq:bc-k} yields the geometric formulation announced in the Introduction. This formulation is more than a convenient rewriting: it shows that the same reduced equation is valid for every Hsiang--Lawson action and explains why the subsequent analysis can be carried out uniformly.

\begin{proposition}\label{prop:geodesic}
Let $U$ be a connected open subset of a $G$-invariant hypersurface of $N^n$ on which $\grad f\neq0$, and let $\gamma(s)$ be its profile curve parametrized by arc length with respect to the orbit-space metric $g_Q$ (flat for $N^n=\R^n$, round for $N^n=\Sph^n$).
Then the immersion restricted to $U$ is biconservative if and only if
\begin{equation}\label{eq:kappaV}
3\kappa=\nu(\log V),
\end{equation}
where $\kappa=k_d$ is the signed (geodesic) curvature of $\gamma$ with respect to $g_Q$ and $\nu$ is the corresponding $g_Q$-unit normal. Equivalently, the trace of $\gamma$ is a geodesic of the conformal metric
\begin{equation}\label{eq:bcmetric}
\widetilde g=V^{2/3}g_Q.
\end{equation}
More precisely, if $\widetilde s$ denotes arc length with respect to $\widetilde g$, then
\begin{equation}\label{eq:tilde-arclength}
 d\widetilde s=V^{1/3}(\gamma(s))\,ds,
 \qquad
 \frac{d\gamma}{d\widetilde s}=V^{-1/3}\dot\gamma,
\end{equation}
and $\gamma(\widetilde s)$ is an affinely parametrized $\widetilde g$-geodesic.
\end{proposition}

\begin{proof}
By \eqref{eq:bc-k} and \eqref{eq:orbitmean}, the biconservative equation is $3\kappa=\nu(\log V)$. To obtain its geodesic interpretation, let $T=\dot\gamma$ and write the Frenet equation as
\begin{equation}\label{eq:Frenet-gQ}
\nabla^Q_TT=\kappa\nu.
\end{equation}
Put
\begin{equation}\label{eq:u-conformal}
u=\frac13\log V,\qquad \widetilde g=e^{2u}g_Q.
\end{equation}
The arc-length parameters satisfy $d\widetilde s=e^u\,ds$, so the $\widetilde g$-unit tangent and normal are
\begin{equation}\label{eq:Tnu-tilde}
\widetilde T=e^{-u}T,\qquad \widetilde\nu=e^{-u}\nu.
\end{equation}
The conformal connection formula is
\begin{equation}\label{eq:conformal-connection}
\widetilde\nabla_XY
=\nabla^Q_XY+X(u)Y+Y(u)X-g_Q(X,Y)\grad_{g_Q}u.
\end{equation}
Hence $\widetilde\nabla_{\widetilde T}\widetilde T=e^{-2u}\bigl(\nabla^Q_TT+T(u)T-\grad_{g_Q}u\bigr)$, and, since $\grad_{g_Q}u=T(u)T+\nu(u)\nu$ and by \eqref{eq:Frenet-gQ}, it follows that
\begin{equation}\label{eq:conformal-acceleration}
\widetilde\nabla_{\widetilde T}\widetilde T
=e^{-2u}\bigl(\kappa-\nu(u)\bigr)\nu
=e^{-u}\bigl(\kappa-\nu(u)\bigr)\widetilde\nu.
\end{equation}
Thus $\gamma$, parametrized by $\widetilde s$, is an affinely parametrized $\widetilde g$-geodesic if and only if $\kappa=\nu(u)=\nu(\log V)/3$, which is \eqref{eq:kappaV}.
\end{proof}

\begin{proof}[{\bf Proof of Theorem~\ref{thm:A}}]
The first two statements are the content of Proposition~\ref{prop:geodesic}. Conversely, assume that the profile satisfies $3\kappa=\nu(\log V)$ everywhere. By \eqref{eq:orbitmean} this is \eqref{eq:bc-k}, so that $f+2k_d=0$ and \eqref{eq:bc-factor} holds on the whole hypersurface; that is, $\Sigma_\gamma$ is biconservative. Moreover, by \eqref{eq:f} and \eqref{eq:bc-k}, $f=k_d-3k_d=-2\kappa$. Hence $f$ is constant if and only if $\kappa$ is constant, and $\Sigma_\gamma$ is proper biconservative as soon as $\kappa$ is not constant. {The} comparison with the minimal case follows from the same computation: $\Sigma_\gamma$ is minimal if and only if $\kappa=-\sum_i m_ik_i=\nu(\log V)$, which by \eqref{eq:conformal-acceleration} with $u=\log V$ is the geodesic equation of $V^2g_Q$.

{It remains to prove the last statement. Let $\Sigma_\gamma$ be smooth, biconservative and not CMC. Then the set $\{\grad f\neq0\}$ is nonempty and \eqref{eq:bc-k} holds there; however, a priori \eqref{eq:bc-factor} also allows CMC pieces, on which $\dot f=0$ but $f+2k_d\neq0$, to be attached to pieces satisfying \eqref{eq:bc-k}. We show that this cannot happen.

By \eqref{eq:f} and \eqref{eq:orbitmean}, along the profile $f=\kappa-\nu(\log V)$, so that the function $\chi=f+2\kappa=3\kappa-\nu(\log V)$ satisfies $\dot f\,\chi=0$ by \eqref{eq:bc-factor}. Let the profile be parametrized by arc length on an open interval $I$ (or periodically, on $\R$), and let $Z=\{s\in I:\ \chi(s)=0\}$, a closed subset of $I$. The open set $A=\{\dot f\neq0\}$ is contained in $Z$, hence in its interior, and it is nonempty because $f$ is not constant and $I$ is connected. Fix $s_1\in A$ and let $C$ be the connected component of the interior of $Z$ containing $s_1$.

On $C$ the profile satisfies $3\kappa=\nu(\log V)$. In the variables $(x,y,\alpha)$ of \eqref{eq:alpha} (and in the analogous variables on the sphere) this is a first-order system with real-analytic right-hand side on the interior of the orbit space, since $V$ is real analytic and positive there. Hence $\gamma|_C$ is real analytic, and so is $f|_C=-2\kappa|_C$. Assume that $b=\sup C$ is an interior point of $I$. The maximal solution $\gamma_P$ of the same system through a point of $C$ coincides with $\gamma$ on $C$ and is defined, and real analytic, on a neighbourhood of $b$, because $\gamma(b)$ lies in the interior of the orbit space. Since $b\in Z$ is not an interior point of $Z$, there are points of $I\setminus Z$ arbitrarily close to $b$. On the open set $I\setminus Z$ we have $\dot f=0$, so that all the derivatives $f^{(j)}$, $j\ge1$, vanish there, and by continuity $f^{(j)}(b)=0$ for every $j\ge1$. On the other hand, $f=-2\kappa_P$ on $C$, where $\kappa_P$ is the curvature of $\gamma_P$; since $f$ is smooth (this is where the $C^\infty$ regularity is used, since all the derivatives of $f$ are involved), its derivatives at $b$ are the limits from the left of those of $-2\kappa_P$, that is, the derivatives of the analytic function $-2\kappa_P$ at $b$. Hence $\kappa_P$ is constant near $b$, and therefore on $C$, so that $f$ is constant on $C$, which contradicts $s_1\in C\cap A$. Thus $\sup C=\sup I$, and in the same way $\inf C=\inf I$. Hence $C=I$, that is, $\chi\equiv0$, i.e., $3\kappa=\nu(\log V)$ along the whole of $\gamma$, and in particular $f=-2\kappa$ on $\Sigma_\gamma$. In the periodic case the same argument applies on $\R$.}
\end{proof}

\begin{corollary}[Tripling]\label{cor:tripling}
Let $G$ and $G'$ be two actions of Table~\ref{tab:actions} with the same $d$, of dimensions $n$ and $n'$ respectively, and assume that their multiplicities satisfy $m_i=3m_i'$ for all $i$. Then, both in the Euclidean and in the spherical case, a curve in the orbit space is the profile of a proper biconservative $G$-invariant hypersurface of $N^n$ if and only if it is the profile of a minimal $G'$-invariant hypersurface of $N^{n'}$ with nonconstant curvature. This applies to $d=1$ with $n-2=3(n'-2)$, to $d=2$ with $(p,q)=(3p',3q')$, and to $G=U(1)\times\mathrm{Spin}(10)$, $G'=S(U(2)\times U(3))$.
\end{corollary}

\begin{proof}
We have $V=(V')^3$ {(with the normalization of $V$ and $V'$ fixed after \eqref{eq:volume}; for the actual orbit volumes this holds up to a positive constant factor, which does not affect the geodesics)}, so $V^{2/3}g_Q=(V')^2g_Q$. By Theorem~\ref{thm:A} and the Hsiang--Lawson {characterization} recalled after Theorem~\ref{thm:A}, the two equations coincide. On the biconservative side $f=-2\kappa$, so the hypersurface is proper exactly when $\kappa$ is not constant. The list of cases follows from Table~\ref{tab:actions}: for $d\ge3$ the only multiplicities divisible by three are $(9,6)$ for $U(1)\times\mathrm{Spin}(10)$, and $(3,2)$ are those of $S(U(2)\times U(m))$ with $m=3$.
\end{proof}

The conformal-geodesic formulation also yields monotone quantities along the profiles. For a smooth function $\phi$ on the interior of the orbit space we write $\dot\phi=\frac{d}{ds}\phi(\gamma(s))=T(\phi)$.

\begin{lemma}[A weighted identity]\label{lem:weighted}
Let $N^n=\R^n$ or $\Sph^n$, and let $\gamma(s)$ be a curve in the interior of the orbit space, parametrized by arc length with respect to $g_Q$, which satisfies $3\kappa=\nu(\log V)$. Let $\phi$ be a smooth function on the interior of the orbit space such that
\[
\Hess_{g_Q}\phi=\lambda\,g_Q,
\qquad
g_Q\bigl(\grad_{g_Q}\phi,\grad_{g_Q}\log V\bigr)=\mu,
\]
for some functions $\lambda$ and $\mu$. Then, along $\gamma$,
\begin{equation}\label{eq:weighted}
\frac{d}{ds}\Bigl(V^{1/3}\dot\phi\Bigr)=V^{1/3}\Bigl(\lambda+\frac{\mu}{3}\Bigr).
\end{equation}
\end{lemma}

\begin{proof}
We have $\frac{d}{ds}V^{1/3}=\frac13V^{1/3}\,T(\log V)$ and, by \eqref{eq:Frenet-gQ},
\[
\frac{d}{ds}\dot\phi=\Hess_{g_Q}\phi(T,T)+g_Q\bigl(\grad_{g_Q}\phi,\nabla^Q_TT\bigr)=\lambda+\kappa\,\nu(\phi)=\lambda+\frac13\nu(\log V)\,\nu(\phi).
\]
Hence
\[
\frac{d}{ds}\Bigl(V^{1/3}\dot\phi\Bigr)=V^{1/3}\Bigl(\lambda+\frac13\bigl[T(\log V)\,T(\phi)+\nu(\log V)\,\nu(\phi)\bigr]\Bigr),
\]
and, since $\{T,\nu\}$ is a $g_Q$-orthonormal frame, the bracket equals $g_Q(\grad_{g_Q}\log V,\grad_{g_Q}\phi)=\mu$.
\end{proof}

Geometrically, $\grad_{g_Q}\phi$ is a conformal vector field of $g_Q$, hence also of $\widetilde g=V^{2/3}g_Q$, and \eqref{eq:weighted} is the derivative of $\widetilde g(\grad_{g_Q}\phi,\widetilde T)=V^{1/3}\dot\phi$ along the $\widetilde g$-geodesic $\gamma$. When $\grad_{g_Q}\phi$ is a Killing field of $\widetilde g$, i.e., $\lambda+\mu/3=0$, \eqref{eq:weighted} is a Clairaut-type first integral. The same computation, with $\kappa=\nu(\log V)$, shows that along the profiles of minimal $G$-invariant hypersurfaces $\frac{d}{ds}(V\dot\phi)=V(\lambda+\mu)$.

\begin{corollary}[Monotonicity of the distance from the origin]\label{cor:monotone}
Let $N^n=\R^n$, and let $\gamma$, defined on an interval $I$, be a profile in the interior of $Q$ which satisfies $3\kappa=\nu(\log V)$. Put $r(s)=|\gamma(s)|$. Then the function
\begin{equation}\label{eq:Lambda}
\Lambda(s)=V^{1/3}\bigl(\gamma(s)\bigr)\,\langle\gamma(s),\dot\gamma(s)\rangle=V^{1/3}\bigl(\gamma(s)\bigr)\,r(s)\,\dot r(s)
\end{equation}
satisfies
\begin{equation}\label{eq:Lambda-dot}
\dot\Lambda=\frac{n+1}{3}\,V^{1/3}>0 .
\end{equation}
Consequently $\dot r$ vanishes at most once on $I$; if $\dot r(s_0)=0$, then $r$ is strictly decreasing on $I\cap(-\infty,s_0]$ and strictly increasing on $I\cap[s_0,+\infty)$. In particular $r$ has no local maximum along $\gamma$, and no profile in the interior of $Q$ is a closed curve.
\end{corollary}

\begin{proof}
Apply Lemma~\ref{lem:weighted} to $\phi=(x^2+y^2)/2$. Then $\Hess_{g_Q}\phi=g_Q$, so that $\lambda=1$, and $\grad_{g_Q}\phi$ is the position vector $\xi$. Since $V$ is positively homogeneous of degree $\sum_im_i=n-2$ by \eqref{eq:volume}, Euler's identity gives $\xi(\log V)=n-2$, that is, $\mu=n-2$. Since $\dot\phi=\langle\gamma,\dot\gamma\rangle=r\dot r$, \eqref{eq:weighted} is \eqref{eq:Lambda-dot}. As $V>0$ and $r>0$ in the interior of $Q$, the zeros of $\dot r$ are the zeros of the strictly increasing function $\Lambda$, and $\Lambda$, hence $\dot r$, changes sign from negative to positive there.
\end{proof}

We now write the equation explicitly in the Euclidean case. Let $\alpha$ be the angle made by the tangent to $\gamma$ with the $x$-axis,
\begin{equation}\label{eq:alpha}
\dot x=\cos\alpha,\qquad \dot y=\sin\alpha.
\end{equation}
Then $k_d=\dot\alpha$, and \eqref{eq:bc-k} becomes
\begin{equation}\label{eq:alpha-general}
3\dot\alpha+
\sum_{i=0}^{d-1}m_i
\frac{\cos\bigl(\alpha-i\pi/d\bigr)}{x\sin(i\pi/d)-y\cos(i\pi/d)}=0.
\end{equation}
Together with \eqref{eq:alpha}, this gives a first-order autonomous system in $(x,y,\alpha)$.

\section{The compactified vector field}\label{sec:compactified}

From now on and up to Section~\ref{sec:d2} we consider the Euclidean case $N^n=\R^n$.

We compactify the Euclidean profile equation so that the singular boundary rays become invariant components of a smooth planar system. This also includes the trajectories corresponding to profiles meeting singular orbits.

\begin{figure}[!htbp]
\centering
\begin{tikzpicture}[line cap=round,
  guide/.style={dashed, line width=.5pt, gray},
  thetaline/.style={densely dotted, line width=.5pt, gray},
  vec/.style={line width=.8pt, -{Latex[length=2mm]}},
  lab/.style={font=\scriptsize, fill=white, inner sep=1pt}]
\draw[line width=1.1pt] (6.6,0) -- (0,0) -- (0,6.6);
\node[font=\scriptsize, anchor=north] at (4.7,0) {$\vartheta=0$};
\node[font=\scriptsize, anchor=east] at (0,4.7) {$\vartheta=\frac\pi d$};
\node[anchor=north east] at (6.6,-0.05) {$x$};
\node[anchor=north east] at (-0.05,6.6) {$y$};
\node[font=\small] at (5.6,1.6) {$Q$};
\node[font=\scriptsize, anchor=north east] at (0,0) {$O$};
\draw[thetaline] (0,0) -- (4.55,4.55) node[lab, anchor=south west, rotate=45] {$\vartheta=\vartheta_0$};
\draw[guide] (2,0) arc[start angle=0, end angle=90, radius=2];
\node[font=\scriptsize,anchor=north] at (2,0) {$r_0$};
\draw[line width=1.1pt] plot coordinates {(4.632,4.270) (4.590,4.212) (4.548,4.154) (4.506,4.096) (4.464,4.038) (4.422,3.980) (4.381,3.921) (4.340,3.862) (4.299,3.804) (4.259,3.745) (4.218,3.685) (4.178,3.626) (4.138,3.567) (4.099,3.507) (4.059,3.447) (4.020,3.387) (3.981,3.327) (3.943,3.266) (3.905,3.206) (3.867,3.145) (3.829,3.084) (3.792,3.023) (3.755,2.962) (3.718,2.900) (3.682,2.838) (3.646,2.777) (3.610,2.714) (3.575,2.652) (3.540,2.590) (3.505,2.527) (3.471,2.464) (3.437,2.401) (3.403,2.338) (3.370,2.274) (3.337,2.210) (3.305,2.147) (3.273,2.083) (3.241,2.018) (3.210,1.954) (3.179,1.889) (3.148,1.825) (3.118,1.760) (3.088,1.695) (3.058,1.629) (3.029,1.564) (3.000,1.499) (2.971,1.433) (2.942,1.368) (2.913,1.302) (2.885,1.236) (2.856,1.171) (2.827,1.105) (2.798,1.040) (2.769,0.974) (2.739,0.909) (2.708,0.845) (2.675,0.781) (2.642,0.717) (2.606,0.655) (2.568,0.595) (2.526,0.537) (2.479,0.482) (2.427,0.433) (2.369,0.392) (2.303,0.363) (2.233,0.350) (2.162,0.356) (2.095,0.381) (2.036,0.421) (1.986,0.472) (1.945,0.531) (1.913,0.595) (1.888,0.662) (1.869,0.731) (1.855,0.801) (1.847,0.872) (1.843,0.944) (1.842,1.015) (1.846,1.087) (1.852,1.158) (1.862,1.229) (1.874,1.300) (1.889,1.370) (1.906,1.439) (1.926,1.508) (1.947,1.577) (1.971,1.644) (1.996,1.711) (2.022,1.778) (2.051,1.844) (2.080,1.909) (2.111,1.974) (2.143,2.038) (2.177,2.101) (2.211,2.164) (2.247,2.226) (2.283,2.288) (2.320,2.349) (2.358,2.410) (2.397,2.470) (2.437,2.529) (2.477,2.589) (2.518,2.647) (2.560,2.706) (2.602,2.764) (2.644,2.821) (2.687,2.878) (2.731,2.935) (2.775,2.992) (2.820,3.048) (2.864,3.104) (2.910,3.159) (2.955,3.215) (3.001,3.270) (3.047,3.324) (3.094,3.379) (3.140,3.433) (3.187,3.487) (3.235,3.541) (3.282,3.595) (3.330,3.648) (3.378,3.701) (3.426,3.754) (3.474,3.807) (3.523,3.860) (3.571,3.913) (3.620,3.965) (3.669,4.017) (3.718,4.070) (3.767,4.122) (3.817,4.174) (3.866,4.225) (3.916,4.277) (3.965,4.329) (4.015,4.380) (4.065,4.432) (4.115,4.483) (4.165,4.534) (4.215,4.585) (4.266,4.636)};
\node[font=\small, anchor=west] at (4.3,2.9) {$\gamma$};
\fill (1.892,0.648) circle (1.6pt);
\node[font=\scriptsize, anchor=north] at (1.59,0.8) {$\gamma(0)$};
\coordinate (X) at (1.858,1.202);
\draw[line width=.5pt] (0,0) -- (X);
\node[lab, anchor=south east] at (0.85,0.6) {$r$};
\draw[line width=.5pt] (0.9,0) arc[start angle=0, end angle=32.9, radius=0.9];
\node[font=\scriptsize] at (16.5:1.15) {$\vartheta$};
\draw[guide] (X) -- ++(32.9:1.5);
\draw[guide] (X) -- ++(1.4,0);
\draw[vec] (X) -- ++(82.03:1.35) node[font=\small, anchor=south] {$T$};
\draw[vec] (X) -- ++(172.03:0.95) node[font=\small, anchor=east] {$\nu$};
\draw[line width=.5pt] ($(X)+(1.05,0)$) arc[start angle=0, end angle=82.03, radius=1.05];
\node[font=\scriptsize] at ($(X)+(27:1.27)$) {$\alpha$};
\draw[line width=.5pt] ($(X)+(32.9:0.6)$) arc[start angle=32.9, end angle=82.03, radius=0.6];
\node[font=\scriptsize] at ($(X)+(58:0.42)$) {$\delta$};
\fill (X) circle (1.4pt);
\end{tikzpicture}
\caption{The orbit space $Q$ of $\R^n$ and the variables used in Sections~\ref{sec:compactified}--\ref{sec:focus-node} (here $d=2$; the profile shown is a profile of Theorem~\ref{thm:B} for $p=q=4$). The thick rays $\vartheta=0$ and $\vartheta=\pi/d$ correspond to the singular orbits, and the dotted ray $\vartheta=\vartheta_0$ to the minimal cone $C_G$. At a point of the profile $\gamma$ with polar coordinates $(r,\vartheta)$, $\alpha$ is the angle of the unit tangent $T$ with the direction $\vartheta=0$, $\nu$ is the unit normal, and $\delta=\alpha-\vartheta$ is the angle from the outward radial direction (dashed) to $T$. Thus $\delta=0$ (resp.\ $\delta=\pi$) means that $\gamma$ moves radially away from (resp.\ towards) the origin, $\delta=\pi/2$ that $\gamma$ is tangent to the circle $r=\mathrm{const}$, as at the point $\gamma(0)$ where the profile touches the circle $r=r_0$ (dashed).}
\label{fig:orbit-space-angles}
\end{figure}
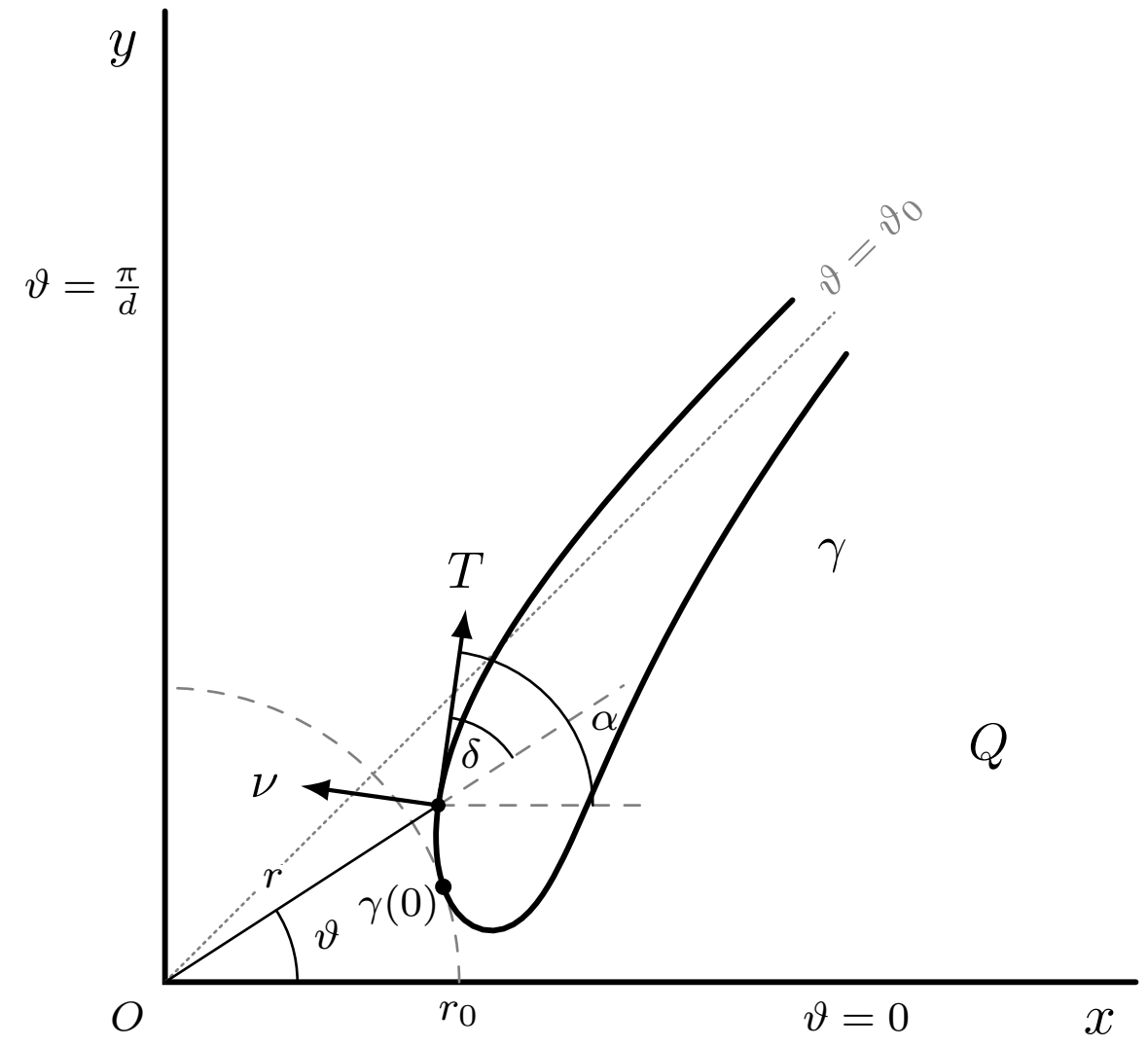

We now pass to polar variables in the orbit space. Write
\begin{equation}\label{eq:polar-new}
x=r\cos\vartheta,\qquad y=r\sin\vartheta,
\qquad 0<\vartheta<\frac{\pi}{d},
\end{equation}
and set
\begin{equation}\label{eq:delta-new}
\delta=\alpha-\vartheta.
\end{equation}
Thus $\delta$ is the oriented angle from the radial direction to the tangent of the profile curve. These variables are illustrated in Figure~\ref{fig:orbit-space-angles}.

 Define
\begin{equation}\label{eq:ell-new}
\ell(\vartheta)=\log v(\vartheta),
\qquad
v(\vartheta)=\prod_{i=0}^{d-1}
\left|\sin\left(\frac{i\pi}{d}-\vartheta\right)\right|^{m_i}.
\end{equation}
We now write the biconservative equation \eqref{eq:kappaV} in these variables. Along the profile, $\dot r=\cos\delta$ and $\dot\vartheta=\sin\delta/r$, the unit tangent and normal are $T=\cos\delta\,e_r+\sin\delta\,e_\vartheta$ and $\nu=-\sin\delta\,e_r+\cos\delta\,e_\vartheta$, where $e_r=\partial_r$ and $e_\vartheta=r^{-1}\partial_\vartheta$, and the curvature is $\kappa=\dot\alpha=\dot\delta+\dot\vartheta$. Since $V(r,\vartheta)=r^{n-2}v(\vartheta)$,
\[
\nu(\log V)=-(n-2)\frac{\sin\delta}{r}+\ell'(\vartheta)\frac{\cos\delta}{r} .
\]
Hence the biconservative equation $3\kappa=\nu(\log V)$, that is, 
\[
3\dot\delta+3\sin\delta/r=\bigl(\ell'(\vartheta)\cos\delta-(n-2)\sin\delta\bigr)/r,
\]
 is equivalent to the system
\begin{equation}\label{eq:polar-s-new}
\dot r=\cos\delta,
\qquad
\dot\vartheta=\frac{\sin\delta}{r},
\qquad
\dot\delta=\frac{1}{3r}
\left(\ell'(\vartheta)\cos\delta-(n+1)\sin\delta\right),
\end{equation}
where dot denotes differentiation with respect to arc length $s$.

Introducing a parameter $t$ by
\begin{equation}\label{eq:tparam-new}
\frac{dt}{ds}=\frac{1}{3r},
\end{equation}
we obtain the autonomous system
\begin{equation}\label{eq:raw-planar}
\left \{
\begin{aligned}
\vartheta'&=3\sin\delta,\\
\delta'&=\ell'(\vartheta)\cos\delta-(n+1)\sin\delta,
\end{aligned}
\right .
\end{equation}
where prime denotes differentiation with respect to $t$.

The system \eqref{eq:raw-planar} is perfectly suited to the interior of the Weyl chamber, but $\ell'$ diverges at its boundary. Following the philosophy of the compact phase portraits used for invariant minimal hypersurfaces in \cite{ABPRS}, it is useful to remove this artificial singularity by a further positive time change.

Set
\begin{equation}\label{eq:h-a}
h(\vartheta)=\frac{1}{d}\sin(d\vartheta),
\qquad
a(\vartheta)=h(\vartheta)\ell'(\vartheta).
\end{equation}
Since $h>0$ on $(0,\pi/d)$, multiplying the vector field by $h$ does not change its trajectories in the interior. We therefore consider
\begin{equation}\label{eq:compactified}
\left \{ 
\begin{aligned}
\vartheta'&=3h(\vartheta)\sin\delta,\\
\delta'&=a(\vartheta)\cos\delta-(n+1)h(\vartheta)\sin\delta.
\end{aligned}
\right .
\end{equation}
From now on prime will refer to the compactified time parameter.

The first issue is to check that the apparent singularity has indeed disappeared. The following elementary lemma is what allows us to regard the boundary rays as ordinary invariant boundary components of a smooth planar flow.

\begin{lemma}\label{lem:a-smooth}
The function $a$ extends smoothly to the closed interval $[0,\pi/d]$. Moreover,
\begin{equation}\label{eq:a-boundary}
a(0)=\mu_-,\qquad a\left(\frac\pi d\right)=-\mu_+,
\end{equation}
where
\[
\mu_-=m_0,
\qquad
\mu_+=\begin{cases}
m_1,& d\ge2,\\
m_0,& d=1.
\end{cases}
\]
In particular, $a(0)>0$ and $a(\pi/d)<0$.
\end{lemma}

\begin{proof}
From \eqref{eq:ell-new},
\[
\ell'(\vartheta)=\sum_{i=0}^{d-1}m_i
\cot\left(\vartheta-\frac{i\pi}{d}\right).
\]
For $d\ge2$, on the interval $[0,\pi/d]$ only the terms $i=0$ and $i=1$ can be singular, and multiplication by $\sin(d\vartheta)$ cancels the corresponding simple poles. More precisely,
\[
\lim_{\vartheta\to 0^+}
\frac{1}{d}\sin(d\vartheta)\cot\vartheta=1,
\qquad
\lim_{\vartheta\to (\pi/d)^{-}}
\frac{1}{d}\sin(d\vartheta)
\cot\left(\vartheta-\frac\pi d\right)=-1.
\]
All the remaining terms are smooth up to the boundary. The case $d=1$ follows directly from
\[
a(\vartheta)=m_0\sin\vartheta\cot\vartheta=m_0\cos\vartheta.
\]
\end{proof}

Having regularized the vector field at the boundary, we next identify the distinguished direction in the interior. This is the direction of the minimal cone, and its uniqueness follows from a simple but important concavity property of $\ell$. Indeed,
\begin{equation}\label{eq:ell-second-new}
\ell''(\vartheta)=
-\sum_{i=0}^{d-1}m_i
\csc^2\left(\vartheta-\frac{i\pi}{d}\right)<0.
\end{equation}
Furthermore,
\[
\ell'(\vartheta)\longrightarrow +\infty
\quad\text{as }\vartheta\to 0^+,
\qquad
\ell'(\vartheta)\longrightarrow -\infty
\quad\text{as }\vartheta\to (\pi/ d)^-.
\]
Thus there is a unique $\vartheta_0\in(0,\pi/d)$ such that
\begin{equation}\label{eq:theta0-new}
\ell'(\vartheta_0)=0.
\end{equation}
The radial ray $\vartheta=\vartheta_0$ is the profile of the unique minimal cone associated with the action.

Since the variable $\delta$ is $2\pi$-periodic, it is enough to work in the closed cylinder
\[
\mathcal C=\left[0,\frac\pi d\right]\times
\bigl(\mathbb R/2\pi\mathbb Z\bigr).
\]
The compactified vector field extends smoothly to $\mathcal C$.

We can now identify all equilibria of the compactified problem. The two interior points represent the same minimal cone with opposite orientations; the four boundary points play the role of the boundary saddles familiar from the phase portraits of invariant minimal hypersurfaces. Their hyperbolicity is essential for the existence and uniqueness of the separatrices considered in the next section.

\begin{proposition}\label{prop:singular-points}
Modulo the $2\pi$-periodicity in $\delta$, the singular points of \eqref{eq:compactified} are the two interior points
\begin{equation}\label{eq:P01-new}
P_0=(\vartheta_0,0),
\qquad
P_1=(\vartheta_0,\pi),
\end{equation}
and the four boundary points
\begin{equation}\label{eq:Q-points}
Q_1=\left(0,\frac\pi2\right),\quad
Q_2=\left(0,\frac{3\pi}{2}\right),\quad
Q_3=\left(\frac\pi d,\frac\pi2\right),\quad
Q_4=\left(\frac\pi d,\frac{3\pi}{2}\right).
\end{equation}
All four boundary points are hyperbolic saddles.
\end{proposition}

\begin{proof}
In the interior, the first equation in \eqref{eq:compactified} gives $\sin\delta=0$, and the second then gives $a(\vartheta)=0$. Since $h>0$, this is equivalent to $\ell'(\vartheta)=0$, hence to $\vartheta=\vartheta_0$. On the boundary $h=0$, so a singular point must satisfy $a(\vartheta)\cos\delta=0$. By Lemma~\ref{lem:a-smooth}, $a$ does not vanish at either endpoint, hence $\cos\delta=0$.

For example, at $Q_1$ the linearization is
\[
DX(Q_1)=
\begin{pmatrix}
3&0\\
-(n+1)&-\mu_-
\end{pmatrix},
\]
whose eigenvalues are $3$ and $-\mu_-$. At $Q_2$ they are $-3$ and $\mu_-$, while at $Q_3$ and $Q_4$ they are respectively $-3,\mu_+$ and $3,-\mu_+$. Thus all $Q_i$ are hyperbolic saddles.
\end{proof}

\begin{remark}\label{rem:boundary-meaning}
The vertical boundary components $\vartheta=0$ and $\vartheta=\pi/d$ are invariant for the compactified flow. The saddle points $Q_i$ encode profile curves which meet a singular boundary ray orthogonally: at $Q_1,Q_2$ one has $\alpha=\delta\equiv\pi/2\pmod\pi$, while at $Q_3,Q_4$ one has $\alpha=\pi/d+\pi/2\pmod\pi$.
\end{remark}

\section{Global phase portrait and separatrices}\label{sec:global-phase}

The preceding section gives the local structure of the compactified flow, but the main existence theorem requires a global statement: trajectories issued from the middle line $\delta=\pi/2$ must be shown to converge to the minimal-cone equilibria rather than to a boundary saddle, a periodic orbit, or a more complicated invariant set. In the minimal problem of \cite{ABPRS} this type of question is handled by a detailed region analysis together with Bendixson and Poincar\'e--Bendixson arguments. In the present biconservative setting there is a simpler global mechanism, namely a strict Lyapunov function.

Consider the two open strips
\begin{align*}
\mathcal R_1&=\left\{0<\vartheta<\frac\pi d,
-\frac\pi2<\delta<\frac\pi2\right\},\\
\mathcal R_2&=\left\{0<\vartheta<\frac\pi d,
\frac\pi2<\delta<\frac{3\pi}2\right\}.
\end{align*}
The first is positively invariant, while the second is negatively invariant. Indeed, along $\delta=\pi/2$ the second component of the raw system \eqref{eq:raw-planar} is $-(n+1)$, whereas along $\delta=-\pi/2$ it is $n+1$.

\begin{lemma}\label{lem:Lyap-new}
On $\mathcal R_1$ the function
\begin{equation}\label{eq:Lyap-new}
\mathcal L(\vartheta,\delta)
=\ell(\vartheta)+3\log(\cos\delta)
\end{equation}
is strictly increasing along every nonstationary trajectory. More precisely, for the raw system \eqref{eq:raw-planar},
\begin{equation}\label{eq:Lyapder-new}
\frac{d\mathcal L}{dt}
=3(n+1)\frac{\sin^2\delta}{\cos\delta}\ge0,
\end{equation}
with equality only when $\delta=0$. On $\mathcal R_2$, the function
\[
\ell(\vartheta)+3\log|\cos\delta|
\]
is non-increasing in forward time, hence non-decreasing in backward time.
\end{lemma}

\begin{proof}
Using \eqref{eq:raw-planar},
\[
\frac{d}{dt}\ell(\vartheta)=3\ell'(\vartheta)\sin\delta,
\]
while
\[
\frac{d}{dt}\bigl(3\log|\cos\delta|\bigr)
=-3\tan\delta\,\delta'
=-3\ell'(\vartheta)\sin\delta
+3(n+1)\frac{\sin^2\delta}{\cos\delta}.
\]
Adding gives the result. Multiplication of the vector field by the positive factor $h(\vartheta)$ does not change the sign of the derivative in the interior.
\end{proof}

The Lyapunov function immediately turns the local direction of the vector field along $\delta=\pi/2$ into a global connection between the two minimal-cone equilibria. This is the dynamical core of the complete-family theorem stated in the Introduction.

Recall that, for a trajectory $\Gamma=\{z(\tau):\tau\in\R\}$ of a smooth vector field on the compact cylinder $\mathcal C$, the \emph{$\omega$-limit set} $\Omega(\Gamma)$ is the set of points $p\in\mathcal C$ such that $z(\tau_k)\to p$ for some sequence $\tau_k\to+\infty$. The \emph{$\alpha$-limit set} $\mathcal A(\Gamma)$ is defined in the same way with $\tau_k\to-\infty$; we write $\mathcal A$ to avoid confusion with the angle $\alpha$. Since $\mathcal C$ is compact, both sets are nonempty, compact, connected and invariant under the flow.

\begin{theorem}\label{thm:global-trajectory}
Fix
\[
\vartheta_*\in\left(0,\frac\pi d\right)
\]
and let $\Gamma_{\vartheta_*}$ be the trajectory of \eqref{eq:compactified} through
\begin{equation}\label{eq:initial-middle}
(\vartheta_*,\pi/2).
\end{equation}
Then
\begin{equation}\label{eq:alpha-omega}
\mathcal A(\Gamma_{\vartheta_*})=\{P_1\},
\qquad
\Omega(\Gamma_{\vartheta_*})=\{P_0\}.
\end{equation}
Equivalently,
\[
(\vartheta,\delta)\longrightarrow(\vartheta_0,\pi)
\quad\text{in backward time},
\qquad
(\vartheta,\delta)\longrightarrow(\vartheta_0,0)
\quad\text{in forward time}.
\]
Moreover, the trajectory meets the line $\delta=\pi/2$ exactly once.
\end{theorem}

\begin{figure}[!htbp]
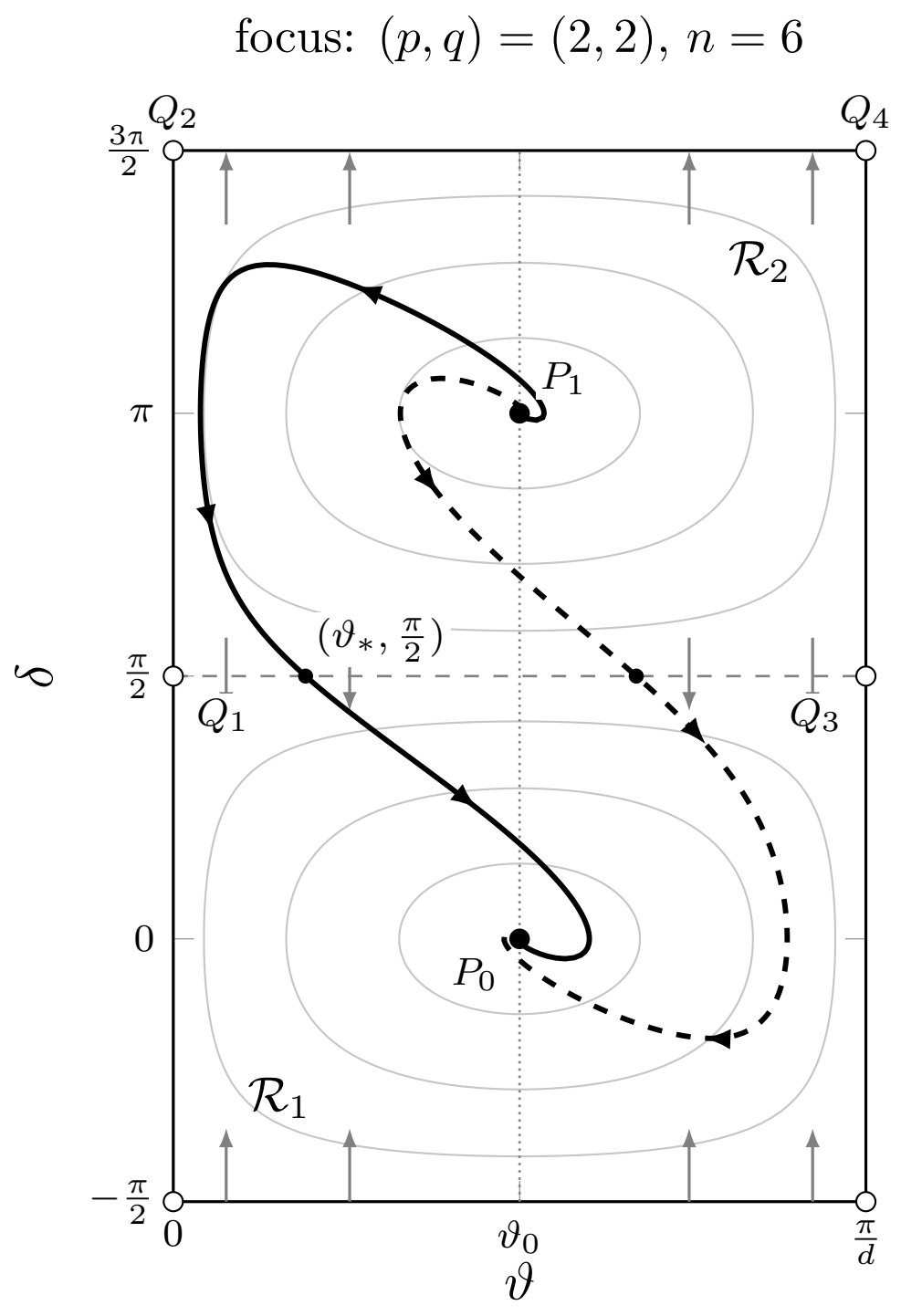

\centering

\caption{The mechanism of Theorem~\ref{thm:global-trajectory} for $d=2$, $p=q=2$ ($n=6$, $\Delta_2=-47$, focus case), on the fundamental domain $-\pi/2\le\delta\le3\pi/2$, whose horizontal sides are identified. Small grey arrows: the vector field enters the strip $\mathcal R_1$ and leaves $\mathcal R_2$ across $\delta=\pm\pi/2$. Thin grey curves: level sets of the Lyapunov functions of Lemma~\ref{lem:Lyap-new}. Thick curves: the trajectories $\Gamma_{\vartheta_*}$ for $\vartheta_*=0.30$ (solid) and $\vartheta_*=1.05$ (dashed), from the source $P_1$ to the sink $P_0$. White circles: the boundary saddles $Q_i$.}
\label{fig:global-trajectory}
\end{figure}

\begin{proof}
Let
\[
z(\tau)=(\vartheta(\tau),\delta(\tau))
\]
denote the solution of the compactified system \eqref{eq:compactified} such that
\[
z(0)=(\vartheta_*,\pi/2).
\]
Here $\tau$ denotes the compactified time parameter. The mechanism of the proof is illustrated in Figure~\ref{fig:global-trajectory}. Since the vector field in \eqref{eq:compactified} is smooth on the compact cylinder $\mathcal C$, the solution $z(\tau)$ is defined for all $\tau\in\mathbb R$.

At the initial point one has $h(\vartheta_*)>0$, and therefore
\[
\vartheta'(0)=3h(\vartheta_*)>0,
\qquad
\delta'(0)=-(n+1)h(\vartheta_*)<0.
\]
Consequently, for $\tau>0$ sufficiently small the trajectory enters $\mathcal R_1$, whereas for $\tau<0$ sufficiently close to zero it lies in $\mathcal R_2$.

We first study the positive semi-orbit. Notice that $\mathcal R_1$ is positively invariant. Indeed, along its upper horizontal boundary $\delta=\pi/2$ one has
\[
\delta'=-(n+1)h(\vartheta)<0,
\]
while along its lower horizontal boundary $\delta=-\pi/2$ one has
\[
\delta'=(n+1)h(\vartheta)>0
\]
for $0<\vartheta<\pi/d$. Thus the vector field points strictly into $\mathcal R_1$ along both horizontal sides. The two vertical sides are invariant for the compactified flow because $h=0$ there. By uniqueness for the smooth system \eqref{eq:compactified}, a trajectory starting in the interior cannot meet one of these invariant vertical sides in finite time. It follows that
\[
z(\tau)\in\mathcal R_1\qquad\text{for every }\tau>0.
\]

Fix now $\tau_1>0$. Since $z(\tau_1)\in\mathcal R_1$, the number
\[
c_1=\mathcal L(z(\tau_1))
\]
is finite. On $\mathcal R_1$, differentiation with respect to compactified time gives, by Lemma~\ref{lem:Lyap-new},
\begin{equation}\label{eq:Lyap-compact-time}
\frac{d\mathcal L}{d\tau}
=h(\vartheta)\,
3(n+1)\frac{\sin^2\delta}{\cos\delta}\ge0.
\end{equation}
Hence
\[
\mathcal L(z(\tau))\ge c_1,
\qquad \tau\ge\tau_1.
\]
On the other hand,
\[
\mathcal L(\vartheta,\delta)
=\ell(\vartheta)+3\log(\cos\delta)
\le \ell(\vartheta_0),
\]
because $\ell$ has its unique maximum at $\vartheta_0$ and $\log(\cos\delta)\le0$ in $\mathcal R_1$. Moreover,
\[
\mathcal L(\vartheta,\delta)\longrightarrow-\infty
\]
when $\vartheta\to0$, $\vartheta\to\pi/d$, or $\delta\to\pm\pi/2$. Therefore the superlevel set
\[
K_{c_1}=
\left\{(\vartheta,\delta)\in\mathcal R_1:
\mathcal L(\vartheta,\delta)\ge c_1\right\}
\]
is a compact subset of $\mathcal R_1$. Thus the whole positive tail
\[
\{z(\tau):\tau\ge\tau_1\}
\]
is contained in the compact set $K_{c_1}$. In particular, its $\omega$-limit set is nonempty, compact, and invariant.

Since $\mathcal L(z(\tau))$ is monotone increasing and bounded above, it converges to a finite limit, say $L_+$. Every point of the $\omega$-limit set has $\mathcal L$-value $L_+$. By invariance of the $\omega$-limit set, $\mathcal L$ must remain constant along every orbit contained in it. Therefore the right-hand side of \eqref{eq:Lyap-compact-time} must vanish on the $\omega$-limit set. Because $K_{c_1}$ is contained in the interior of $\mathcal R_1$, there $h(\vartheta)>0$ and $\cos\delta>0$; hence
\[
\sin\delta=0.
\]
Since $-\pi/2<\delta<\pi/2$, this gives $\delta=0$. Consequently
\[
\Omega(\Gamma_{\vartheta_*})\subset
\{(\vartheta,0):0<\vartheta<\pi/d\}.
\]
However, on the line $\delta=0$ the compactified system reads
\[
\vartheta'=0,
\qquad
\delta'=a(\vartheta)=h(\vartheta)\ell'(\vartheta).
\]
Thus a point of this line can belong to an invariant subset only if $a(\vartheta)=0$. Since $h>0$ in the interior and $\ell'$ has the unique zero $\vartheta_0$, the largest invariant subset of $\{\delta=0\}$ is the singleton $\{P_0\}$. It follows that
\[
\Omega(\Gamma_{\vartheta_*})=\{P_0\},
\]
and therefore
\[
(\vartheta(\tau),\delta(\tau))
\longrightarrow(\vartheta_0,0)
\qquad\text{as }\tau\to+\infty.
\]

We next consider the negative semi-orbit. The strip $\mathcal R_2$ is negatively invariant. Indeed, along its lower boundary $\delta=\pi/2$ the compactified vector field has $\delta'=-(n+1)h(\vartheta)<0$, and along its upper boundary $\delta=3\pi/2$ it has $\delta'=(n+1)h(\vartheta)>0$. Thus, after reversing the direction of time, the vector field points strictly into $\mathcal R_2$ along both horizontal sides; as before, the vertical sides are invariant. Hence
\[
z(\tau)\in\mathcal R_2\qquad\text{for every }\tau<0.
\]

To determine the backward limit, put on $\mathcal R_2$
\[
\mathcal L_2(\vartheta,\delta)
=\ell(\vartheta)+3\log|\cos\delta|.
\]
Along the compactified flow,
\[
\frac{d\mathcal L_2}{d\tau}
=h(\vartheta)\,
3(n+1)\frac{\sin^2\delta}{\cos\delta}\le0,
\]
because $\cos\delta<0$ on $\mathcal R_2$. Hence $\mathcal L_2$ is non-decreasing when the trajectory is followed backward in time. The same superlevel-set argument shows that, after any fixed negative time, the negative semi-orbit is contained in a compact subset of $\mathcal R_2$. Its $\alpha$-limit set is therefore nonempty, compact, and invariant, and on it the preceding derivative must vanish. This gives $\sin\delta=0$, and now the only possibility in $\mathcal R_2$ is $\delta=\pi$. On that line,
\[
\vartheta'=0,
\qquad
\delta'=-a(\vartheta),
\]
so the only invariant point is $P_1=(\vartheta_0,\pi)$. Therefore
\[
\mathcal A(\Gamma_{\vartheta_*})=\{P_1\},
\qquad
(\vartheta(\tau),\delta(\tau))
\longrightarrow(\vartheta_0,\pi)
\quad\text{as }\tau\to-\infty.
\]

Finally, choose the continuous lift of $\delta$ for which $\delta(0)=\pi/2$. We have proved that
\[
-\frac\pi2<\delta(\tau)<\frac\pi2
\quad\text{for }\tau>0,
\]
whereas
\[
\frac\pi2<\delta(\tau)<\frac{3\pi}2
\quad\text{for }\tau<0.
\]
Hence the trajectory can meet the line $\delta=\pi/2$ only at $\tau=0$. This proves the last assertion and completes the proof.
\end{proof}

The trajectories of Theorem~\ref{thm:global-trajectory} stay in the interior of $\mathcal C$, i.e., the corresponding profiles never meet a singular orbit. The only interior trajectories which converge to a point of the boundary circles are those converging to one of the saddles $Q_i$, that is, the branches of their one-dimensional stable and unstable manifolds contained in the interior of $\mathcal C$; the remaining branches lie in the invariant boundary circles, or outside $\mathcal C$, and have no geometric meaning. 
In the minimal problem studied in \cite{ABPRS}, analogous separatrices organize the global classification of profiles which meet singular orbits. Here the Lyapunov function makes their global limits particularly transparent.

\begin{proposition}\label{prop:separatrices}
Each boundary saddle has a unique one-dimensional separatrix contained in the interior of the cylinder. Up to the periodicity in $\delta$, these separatrices have the following connections:
\begin{equation}\label{eq:separatrix-connections}
Q_1\longrightarrow P_0,
\qquad
Q_4\longrightarrow P_0,
\qquad
P_1\longrightarrow Q_2,
\qquad
P_1\longrightarrow Q_3.
\end{equation}
The arrows indicate the direction of the compactified flow.
\end{proposition}

\begin{proof}
The stable manifold theorem gives the local one-dimensional stable and unstable manifolds of every $Q_i$. At $Q_1$ the eigenvalue $-\mu_-$ has the vertical eigenvector $(0,1)$, tangent to the invariant boundary $\vartheta=0$, while the eigenvalue $3$ has the eigenvector $(3+\mu_-,-(n+1))$. Hence the unstable manifold of $Q_1$ has exactly one branch in the interior, along which $\vartheta$ increases and $\delta$ decreases from $\pi/2$; this branch therefore enters $\mathcal R_1$. The same computation at the other saddles shows that the interior branch is: the stable branch of $Q_2$, tangent to $(3+\mu_-,-(n+1))$ at $(0,3\pi/2)$, which lies in $\mathcal R_2$; the stable branch of $Q_3$, tangent to $(-(3+\mu_+),n+1)$ at $(\pi/d,\pi/2)$, which lies in $\mathcal R_2$; and the unstable branch of $Q_4$, tangent to $(-(3+\mu_+),n+1)$ at $(\pi/d,-\pi/2)$, which enters $\mathcal R_1$.

Consider the interior unstable branch of $Q_1$, and fix a point $z_1$ on it. By positive invariance it stays in $\mathcal R_1$ in forward time, and the superlevel-set argument in the proof of Theorem~\ref{thm:global-trajectory}, applied with $c_1=\mathcal L(z_1)$, shows that its forward limit is $P_0$. The same argument applies to the interior branch of $Q_4$. For $Q_2$ and $Q_3$ one argues in backward time in $\mathcal R_2$ with the function $\mathcal L_2$, obtaining the backward limit $P_1$.
\end{proof}

The proofs of Theorem~\ref{thm:global-trajectory} and of Proposition~\ref{prop:separatrices} only use the fact that the trajectory has a point in $\mathcal R_1$ (respectively in $\mathcal R_2$). We record the general statement; note that the absence of periodic orbits is a direct consequence of the Lyapunov functions, without any Poincar\'e--Bendixson argument.

\begin{corollary}\label{cor:general-limits}
Let $\Gamma$ be a trajectory of \eqref{eq:compactified} through a point of the interior of $\mathcal C$.
\begin{enumerate}[label=\textup{(\roman*)},leftmargin=2em]
\item If $\Gamma$ has a point in $\mathcal R_1$, then $\Omega(\Gamma)=\{P_0\}$.
\item If $\Gamma$ has a point in $\mathcal R_2$, then $\mathcal A(\Gamma)=\{P_1\}$.
\item The interior of $\mathcal C$ contains no periodic orbit of \eqref{eq:compactified}.
\end{enumerate}
\end{corollary}

\begin{proof}
(i) If $z(\tau_1)\in\mathcal R_1$, then $z(\tau)\in\mathcal R_1$ for every $\tau\ge\tau_1$ by positive invariance, and the argument in the proof of Theorem~\ref{thm:global-trajectory}, with $c_1=\mathcal L(z(\tau_1))$, applies verbatim. (ii) is analogous, in backward time, with the function $\mathcal L_2$.

(iii) Let $\Gamma$ be a periodic orbit in the interior of $\mathcal C$. The complement of $\mathcal R_1\cup\mathcal R_2$ in the interior of $\mathcal C$ consists of the two {segments} $\delta=\pi/2$ and $\delta=3\pi/2$, on which $\delta'=\mp(n+1)h(\vartheta)\neq0$; hence $\Gamma$ meets $\mathcal R_1$ or $\mathcal R_2$. If $\Gamma$ meets $\mathcal R_1$, then, being periodic and $\mathcal R_1$ being positively invariant, $\Gamma\subset\mathcal R_1$. Since $\mathcal L$ is non-decreasing along $\Gamma$ and takes the same value after one period, \eqref{eq:Lyap-compact-time} gives $\sin\delta\equiv0$, i.e., $\delta\equiv0$ along $\Gamma$; then $\vartheta'=0$ and $\delta'=a(\vartheta)=0$, so that $\Gamma=\{P_0\}$, which is not a periodic orbit. If $\Gamma$ meets $\mathcal R_2$, the same argument with $\mathcal L_2$ gives $\Gamma=\{P_1\}$.
\end{proof}

\section{Geometry of the profile curves}\label{sec:profiles}

The phase portrait is only an auxiliary object: the geometric content lies in the corresponding curves in the orbit space and in the hypersurfaces reconstructed from them. In this section we recover the radial variable and arc length from the planar trajectories to prove completeness, describe the ends, and analyse contact with singular orbits.

We begin with the elementary identities relating the phase variables to the profile curve. The first identities are immediate from \eqref{eq:polar-s-new}:
\begin{equation}\label{eq:profile-basic}
\dot r=\cos\delta,
\qquad
\dot\vartheta=\frac{\sin\delta}{r}.
\end{equation}
Moreover, since $\alpha=\vartheta+\delta$,
\begin{equation}\label{eq:kappa-delta}
\kappa=\dot\alpha
=\frac{1}{3r}
\left(\ell'(\vartheta)\cos\delta-(n-2)\sin\delta\right).
\end{equation}

The first application is the main global existence result. Theorem~\ref{thm:global-trajectory} controls the angular variables; the only remaining point is to recover the radial variable and the original arc-length parameter carefully. 

{In the proof of Theorem~\ref{thm:B} we shall use} the following general completeness criterion, which does not require properness.

\begin{lemma}[Completeness criterion]\label{lem:completeness}
Let $N^n=\R^n$ or $\Sph^n$, and let $\gamma:\R\to Q$ \textup{(}respectively $Q_{\Sph}$\textup{)} be a curve in the interior of the orbit space, parametrized by arc length. Then the metric induced by $F_\gamma$ on $\R\times G/H$ is complete.
\end{lemma}

\begin{proof}
At $(s,gH)$, the vector $dF_\gamma(\partial_s)=g\cdot\dot\gamma(s)$ is a unit vector tangent to $g\cdot\Pi_0$ \textup{(}respectively $g\cdot(\Pi_0\oplus\R)$\textup{)}, which meets the orbit through $g\cdot\gamma(s)$ orthogonally, while $dF_\gamma$ maps the tangent space of $\{s\}\times G/H$ onto the tangent space of that orbit. Hence the induced metric is
\[
F_\gamma^*\langle\,,\rangle=ds^2+g_s,
\]
where $g_s$ is the metric induced on $G/H$ by the orbit map $gH\mapsto g\cdot\gamma(s)$. Consequently, the length of any piecewise $C^1$ curve $c(u)=(s(u),p(u))$, $u\in[0,1]$, is at least $\int_0^1|s'(u)|\,du\ge|s(1)-s(0)|$, so that the distance between $(s_1,p_1)$ and $(s_2,p_2)$ is at least $|s_1-s_2|$. Therefore every bounded subset of $\R\times G/H$ is contained in a set $[a,b]\times G/H$, which is compact because $G$ is compact. By the Hopf--Rinow theorem the induced metric is complete.
\end{proof}

\begin{figure}[!htbp]
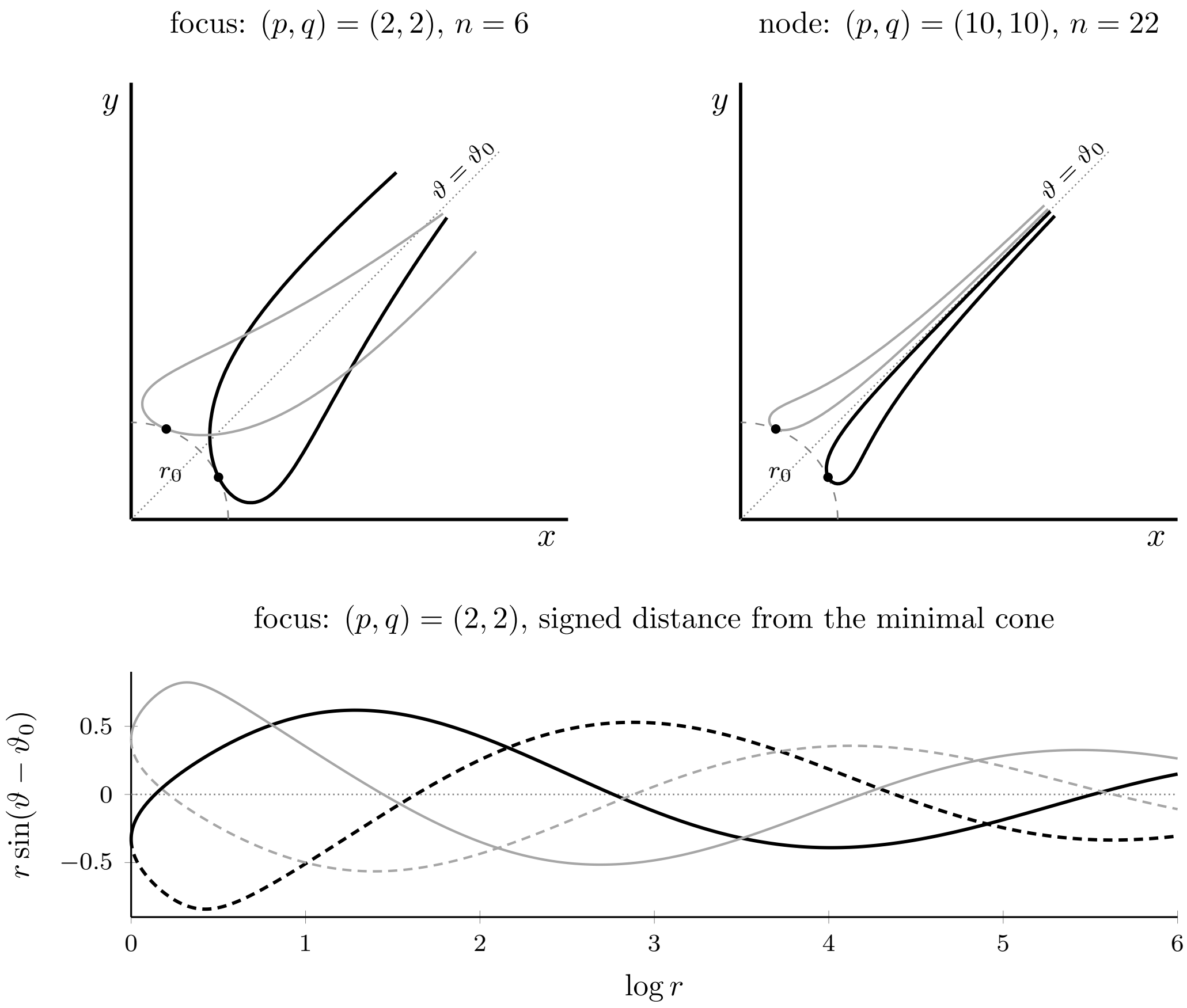

\centering

\caption{Profile curves of Theorem~\ref{thm:B} for the action $SO(p+1)\times SO(q+1)$ ($d=2$, so that $Q$ is the quarter plane and $\vartheta_0=\pi/4$ when $p=q$), with $r_0=1$ and $\vartheta_*=0.45$ (solid) and $\vartheta_*=1.2$ (grey). Each profile touches the circle $r=r_0$ (dashed) at the point of polar angle $\vartheta_*$ (dots), stays outside it and in the interior of $Q$, and both its ends go to infinity in the direction of the minimal cone $\vartheta=\vartheta_0$ (dotted). In the focus case (left) the ends oscillate around the cone direction and cross it infinitely many times, in the node case (right) they approach it without oscillating (Theorem~\ref{thm:D}); since $\vartheta-\vartheta_0$ decays exponentially in $\log r$, only the first crossings are visible at this scale. The coordinate axes are the singular boundary rays of $Q$. Bottom: in the focus case, the signed distance $r\sin(\vartheta-\vartheta_0)$ of the profile from the minimal-cone ray as a function of $\log r$, along the end $s>0$ (solid lines) and along the end $s<0$ (dashed lines). Its zeros are the crossings with the minimal cone.
}
\label{fig:thmB}
\end{figure}

\begin{proof}[{\bf Proof of Theorem~\ref{thm:B}}]

{\emph{Step 1: the profile.} By Proposition~\ref{prop:geodesic}, in the polar variables of Section~\ref{sec:compactified} the equation $3\kappa=\nu(\log V)$ is equivalent to the system \eqref{eq:polar-s-new}, whose right-hand side is smooth on $\{r>0,\ 0<\vartheta<\pi/d\}$. A curve touches the circle $r=r_0$ at the point of polar angle $\vartheta_*$ if and only if, at that point, $r=r_0$, $\vartheta=\vartheta_*$ and the tangent is orthogonal to the radial direction, i.e., $\delta=\pi/2$ or $\delta=3\pi/2$; the two values correspond to the two orientations of the same curve. Fixing the orientation, the profile is therefore the maximal solution of \eqref{eq:polar-s-new}, defined on an open interval $I\ni0$, with the initial data
\begin{equation}\label{eq:profile-init-new}
r(0)=r_0,
\qquad
\vartheta(0)=\vartheta_*,
\qquad
\delta(0)=\frac\pi2 ,
\end{equation}
and it exists and is unique by the standard theory of ordinary differential equations.

\emph{Step 2: radial behaviour.} By Corollary~\ref{cor:monotone}, the function $\Lambda=V^{1/3}r\dot r=V^{1/3}r\cos\delta$ is strictly increasing on $I$, and $\Lambda(0)=0$ because $\delta(0)=\pi/2$. Since $V>0$ and $r>0$, it follows that $\dot r<0$ for $s<0$ and $\dot r>0$ for $s>0$. Thus $r$ is strictly decreasing on $I\cap(-\infty,0]$ and strictly increasing on $I\cap[0,+\infty)$; in particular
\[
r(s)\ge r_0\quad\text{on } I,\qquad\text{with equality only at } s=0 .
\]

\emph{Step 3: the profile is defined for every $s$.} Let $\tau$ be the time of the compactified system \eqref{eq:compactified}, with $\tau=0$ at $(\vartheta_*,\pi/2)$. Since \eqref{eq:compactified} is obtained from \eqref{eq:raw-planar} by multiplying the vector field by $h(\vartheta)$, we have $dt/d\tau=h(\vartheta)$, and by \eqref{eq:tparam-new} and \eqref{eq:profile-basic}
\[
\frac{ds}{d\tau}=3h(\vartheta)\,r,
\qquad
\frac{dr}{d\tau}=3h(\vartheta)\,r\cos\delta .
\]
Now let $(\vartheta(\tau),\delta(\tau))$ be the trajectory $\Gamma_{\vartheta_*}$ of Theorem~\ref{thm:global-trajectory}; since $\mathcal C$ is compact and its boundary is invariant, it is defined for every $\tau\in\R$. The two formulas
\begin{equation}\label{eq:r-tau}
r(\tau)=r_0\exp\Bigl(3\int_0^\tau h\bigl(\vartheta(\sigma)\bigr)\cos\delta(\sigma)\,d\sigma\Bigr),
\qquad
s(\tau)=3\int_0^\tau h\bigl(\vartheta(\sigma)\bigr)\,r(\sigma)\,d\sigma
\end{equation}
then define $r$ and $s$ for every $\tau\in\R$, and the curve with polar coordinates $(r(\tau),\vartheta(\tau))$, parametrized by $s$, solves \eqref{eq:polar-s-new} with the initial data \eqref{eq:profile-init-new}; by uniqueness it is the profile, and $I=s(\R)$. By Theorem~\ref{thm:global-trajectory} the closure of $\Gamma_{\vartheta_*}$ in $\mathcal C$ is $\Gamma_{\vartheta_*}\cup\{P_0,P_1\}$, a compact subset of the interior of $\mathcal C$. Hence there are $0<\vartheta_-\le\vartheta_+<\pi/d$ such that $\vartheta(\tau)\in[\vartheta_-,\vartheta_+]$ for every $\tau$, and therefore $h(\vartheta(\tau))\ge h_*:=\min_{[\vartheta_-,\vartheta_+]}h>0$. Together with Step~2 this gives
\[
\frac{ds}{d\tau}\ge3h_*r_0>0,
\]
so that $s$ is an increasing bijection of $\R$ onto $\R$. Thus $I=\R$, the profile is defined for every value of the arc length, and it stays in the sector $\vartheta_-\le\vartheta\le\vartheta_+$, in particular in the interior of $Q$. Moreover $s\to\pm\infty$ if and only if $\tau\to\pm\infty$.

\emph{Step 4: asymptotics.} By Theorem~\ref{thm:global-trajectory} and Step~3,
\[
(\vartheta,\delta)\to(\vartheta_0,0)\ \text{ as } s\to+\infty,
\qquad
(\vartheta,\delta)\to(\vartheta_0,\pi)\ \text{ as } s\to-\infty .
\]
Hence $\vartheta(s)\to\vartheta_0$, and $\dot r=\cos\delta\to\pm1$ as $s\to\pm\infty$. By l'H\^opital's rule, $r(s)/s\to\pm1$ as $s\to\pm\infty$, that is, $r(s)/|s|\to1$; in particular $r(s)\to+\infty$. This proves (iii). Typical profiles are shown, for $d=2$, in Figure~\ref{fig:thmB}; the way in which the two ends approach the direction $\vartheta_0$ of the minimal cone, and in particular their crossings with it, is discussed in Section~\ref{sec:focus-node}.

\emph{Step 5: completeness and properness.} The hypersurface $\Sigma=\Sigma_\gamma$ is the immersion $F_\gamma:\R\times G/H\to\R^n$, $F_\gamma(s,gH)=g\cdot\gamma(s)$; its image may have self-intersections (Remark~\ref{rem:self-intersections}). By Lemma~\ref{lem:completeness} its induced metric is complete. Since $|F_\gamma(s,gH)|=r(s)$, $r(s)\to+\infty$ as $s\to\pm\infty$ and $G/H$ is compact, the preimage under $F_\gamma$ of every compact subset of $\R^n$ is closed and contained in a set $[a,b]\times G/H$, hence compact; i.e., $F_\gamma$ is proper. This proves (i). Moreover, by Step~2, $|F_\gamma(s,gH)|=r(s)\ge r_0$, with equality if and only if $s=0$; hence $\Sigma$ lies outside the open ball of radius $r_0$ and touches its boundary sphere exactly along the principal orbit $G\cdot\gamma(0)$. This proves (ii).

\emph{Step 6: mean curvature.} By Theorem~\ref{thm:A}, $\Sigma$ is biconservative and $f=-2\kappa$, where, by \eqref{eq:kappa-delta},
\[
f=-\frac{2}{3r}\bigl(\ell'(\vartheta)\cos\delta-(n-2)\sin\delta\bigr).
\]
At $s=0$ we have $\cos\delta=0$ and $\sin\delta=1$, so that $f(0)=2(n-2)/(3r_0)\neq0$. As $s\to\pm\infty$, the bracket tends to $\pm\ell'(\vartheta_0)=0$, since $\vartheta_0$ is the critical point of $\ell$, while $r(s)/|s|\to1$ by Step~4; hence $|s|\,f(s)\to0$, i.e., $f(s)=o(1/|s|)$. In particular $f$ is not constant, and $\Sigma$ is proper biconservative by Theorem~\ref{thm:A}. This proves (iv).

\emph{Step 7: the parameters $r_0$ and $\vartheta_*$.} If $\gamma$ satisfies $3\kappa=\nu(\log V)$, then so does $\lambda\gamma$ for every $\lambda>0$: the curvature is multiplied by $1/\lambda$ and, since $V$ is positively homogeneous by \eqref{eq:volume}, $\log V(\lambda\xi)=\log V(\xi)+\mathrm{const}$, so that $\nu(\log V)$ is also multiplied by $1/\lambda$. Since $\lambda\gamma$ touches the circle $r=\lambda r_0$ at the same polar angle, by uniqueness the hypersurface obtained from $(r_0,\vartheta_*)$ is the image of the one obtained from $(1,\vartheta_*)$ under the homothety $\xi\mapsto r_0\xi$. A homothety does not change $\vartheta$ and $\delta$, so homothetic profiles give the same trajectory of \eqref{eq:compactified}. By the last assertion of Theorem~\ref{thm:global-trajectory}, $\Gamma_{\vartheta_*}$ meets the segment $\delta=\pi/2$ only at $(\vartheta_*,\pi/2)$; hence different values of $\vartheta_*$ give different trajectories, and the corresponding hypersurfaces are not homothetic.}

{\emph{Step 8: the rotational case $d=1$.} For $d=1$ the group $SO(n-1)$ acts on $\R^n=\R\times\R^{n-1}$ and commutes with the translations along the axis $\R\times\{0\}$; in the orbit space these are the horizontal translations $(x,y)\mapsto(x+c,y)$, which preserve $V=y^{n-2}$ and hence equation \eqref{eq:alpha-general}. Moreover, in this case \eqref{eq:alpha-general} reads $3\dot\alpha=(n-2)\cos\alpha/y$, and a direct computation shows that
\begin{equation}\label{eq:d1-first-integral}
y^{(n-2)/3}\cos\alpha
\end{equation}
is a first integral. It is the case $\phi=x$ of Lemma~\ref{lem:weighted}: here $\lambda=0$ and, since $V=y^{n-2}$ does not depend on $x$, $\mu=0$, so that $V^{1/3}\dot x=y^{(n-2)/3}\cos\alpha$ is constant. Consequently, for $d=1$ all the profiles constructed above are congruent to each other modulo homotheties, translations along the axis and the reflection $x\mapsto-x$: the parameter $\vartheta_*$ is not essential in this case, and we recover the complete ``catenary type'' rotational hypersurfaces of \cite[Theorem~6.5]{MOR-BC}. The first integral also shows that $x$ is strictly monotone along each profile and that $dx/dy=O(y^{-(n-2)/3})$ along both ends. Hence, for $n>5$, both ends are asymptotic to hyperplanes $x=\mathrm{const}$ parallel to the minimal ``cone'' $x=0$, but in general distinct from it. This is the reason why, in the statement, we only claim that the ends are asymptotic \emph{in angle} to the minimal cone. For $d\ge2$ the translations along the axis are no longer symmetries and $\vartheta_*$ is a genuine parameter.}
\end{proof}

\begin{remark}\label{rem:self-intersections}
The profiles of Theorem~\ref{thm:B} need not be injective, so the hypersurfaces of Theorem~\ref{thm:B} are, in general, only immersed. For instance, let $d=2$ and $p=q$, so that $\ell$ is symmetric with respect to $\vartheta_0=\pi/4$, and take $\vartheta_*=\vartheta_0$. By uniqueness, the reflection of $Q$ in the minimal-cone ray maps the profile to itself with the opposite orientation: $r(-s)=r(s)$ and $\vartheta(-s)=2\vartheta_0-\vartheta(s)$. Hence $\gamma(-s)=\gamma(s)$ whenever $\vartheta(s)=\vartheta_0$. In the focus case (for instance $p=q=1$, where $\Delta_2=-23$) this happens for infinitely many values of $s>0$ by Theorem~\ref{thm:D}, and the profile has infinitely many self-intersections along the minimal cone.
\end{remark}

\begin{remark}\label{rem:radial-graph}
By Theorem~\ref{thm:B}(ii) the distance $r$ is strictly decreasing on $(-\infty,0]$ and strictly increasing on $[0,+\infty)$ along the profiles of Theorem~\ref{thm:B}. Hence each half of the profile is a radial graph: there are continuous functions $\Theta_\pm:[r_0,+\infty)\to(0,\pi/d)$, smooth on $(r_0,+\infty)$, with $\Theta_+(r_0)=\Theta_-(r_0)=\vartheta_*$, such that
\[
\gamma\bigl([0,+\infty)\bigr)=\bigl\{(r,\Theta_+(r)):\ r\ge r_0\bigr\},
\qquad
\gamma\bigl((-\infty,0]\bigr)=\bigl\{(r,\Theta_-(r)):\ r\ge r_0\bigr\}.
\]
Consequently $\gamma$ is injective if and only if $\Theta_+(r)\neq\Theta_-(r)$ for every $r>r_0$. In this case $F_\gamma$ is injective, because distinct points of the interior of $Q$ correspond to disjoint orbits and the isotropy group of every interior point is $H$, and, being a proper injective immersion, it is an embedding. In the situation of Remark~\ref{rem:self-intersections}, $\Theta_-=2\vartheta_0-\Theta_+$, and the self-intersections occur exactly at the radii at which $\Theta_+(r)=\vartheta_0$. This reformulation may be useful for the problem of deciding for which multiplicities the hypersurfaces of Theorem~\ref{thm:B} are embedded.
\end{remark}

The preceding family consists of two-ended profiles staying in the regular part of the orbit space. The separatrices of Proposition~\ref{prop:separatrices} describe the complementary phenomenon in which one end approaches a singular stratum. Their phase-plane endpoint immediately encodes the angle of contact with the boundary.

\begin{proposition}\label{prop:boundary-profile}
A profile corresponding to an interior separatrix of one of the boundary saddles $Q_i$ reaches the corresponding boundary ray of the orbit space after a finite arc length, at a point at positive distance from the origin, and meets it orthogonally. At its other end it is defined for an infinite arc length, $r\to+\infty$, and it is asymptotic in angle to the minimal cone. Along the whole profile the distance $r$ from the origin is strictly monotone.
\end{proposition}

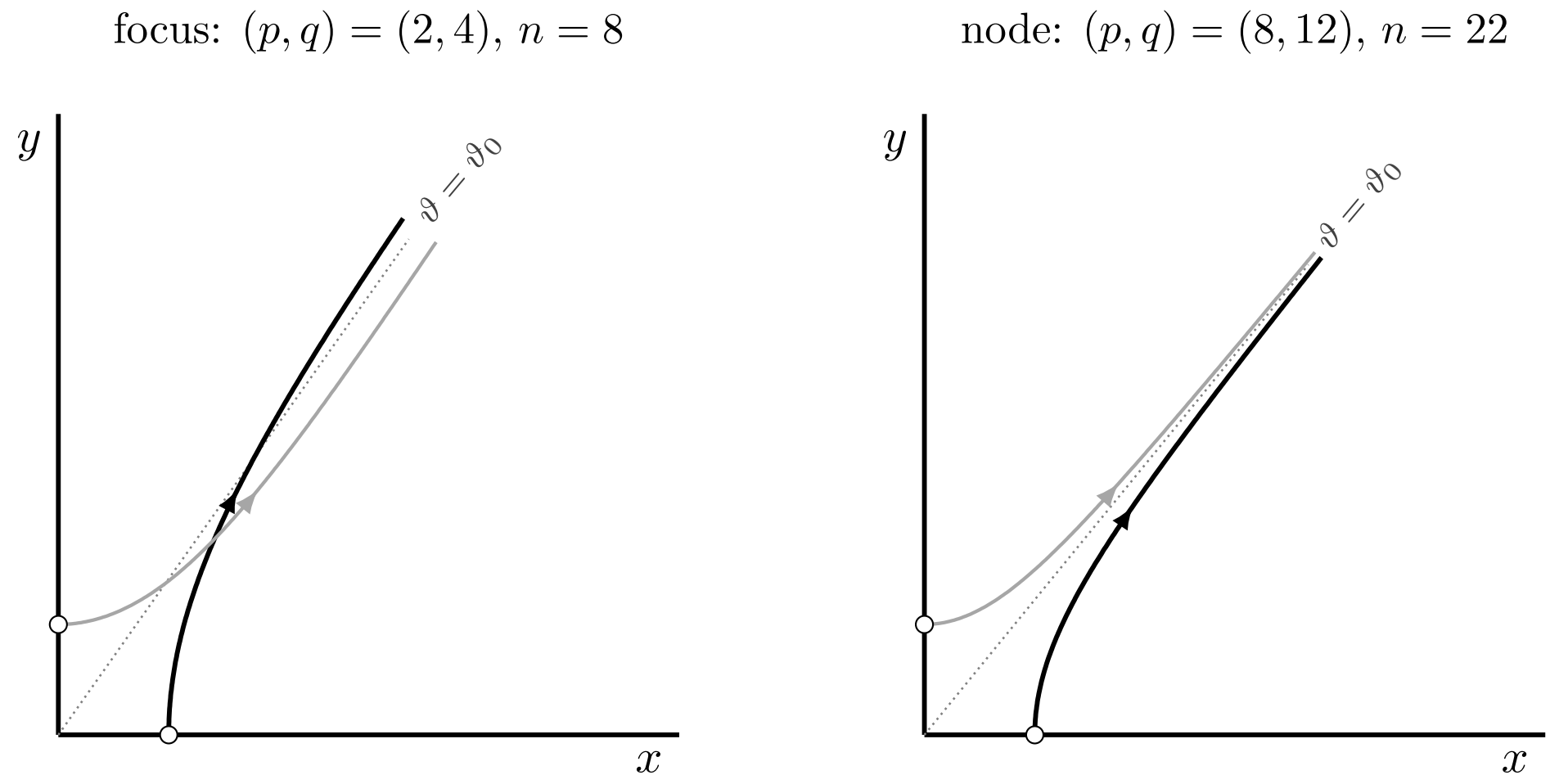
\begin{figure}[!htbp]
\centering
\begin{tikzpicture}[
  guide/.style={dashed, line width=.5pt, gray},
  thetaline/.style={densely dotted, line width=.5pt, gray},
  saddle/.style={circle, draw, fill=white, inner sep=1.4pt},
  lab/.style={font=\scriptsize, fill=white, inner sep=1pt},
  sepc/.style={line width=1.1pt, black, postaction={decorate}, decoration={markings, mark=at position 0.45 with {\arrow{Latex[length=2mm]}}}}]
\begin{axis}[at={(0cm,0)}, width=6.6cm, height=6.6cm, scale only axis=false, axis equal image,
  xmin=0, xmax=4.5, ymin=0, ymax=4.5, axis lines=left, axis line style={line width=1.1pt, -},
  xtick=\empty, ytick=\empty, xlabel={$x$}, ylabel={$y$}, x label style={at={(axis description cs:1,0)}, anchor=north east}, y label style={at={(axis description cs:0,1)}, anchor=north east, rotate=-90}, 
  tick label style={font=\scriptsize}, title={focus: $(p,q)=(2,4)$, $n=8$}, title style={font=\small,yshift=4pt}, clip=false]
  \addplot[thetaline] coordinates {(0,0) (2.540,3.593)};
  \node[font=\scriptsize, gray!50!black, anchor=south west, rotate=50] at (axis cs:2.740,3.493) {$\vartheta=\vartheta_0$};
  \addplot[sepc] coordinates {(0.800,0.000) (0.800,0.019) (0.800,0.038) (0.801,0.057) (0.801,0.076) (0.802,0.095) (0.802,0.114) (0.803,0.133) (0.804,0.152) (0.805,0.171) (0.806,0.190) (0.808,0.209) (0.809,0.228) (0.811,0.247) (0.813,0.266) (0.814,0.285) (0.816,0.304) (0.819,0.323) (0.821,0.342) (0.823,0.361) (0.826,0.379) (0.828,0.398) (0.831,0.417) (0.834,0.436) (0.837,0.455) (0.840,0.474) (0.843,0.492) (0.846,0.511) (0.849,0.530) (0.853,0.548) (0.857,0.567) (0.860,0.586) (0.864,0.604) (0.868,0.623) (0.872,0.642) (0.876,0.660) (0.881,0.679) (0.885,0.697) (0.889,0.716) (0.894,0.734) (0.899,0.753) (0.903,0.771) (0.908,0.790) (0.913,0.808) (0.918,0.826) (0.923,0.845) (0.929,0.863) (0.934,0.881) (0.939,0.899) (0.945,0.918) (0.951,0.936) (0.956,0.954) (0.962,0.972) (0.968,0.990) (0.974,1.008) (0.980,1.026) (0.986,1.044) (0.992,1.062) (0.998,1.080) (1.005,1.098) (1.011,1.116) (1.017,1.134) (1.024,1.152) (1.031,1.170) (1.037,1.188) (1.044,1.205) (1.051,1.223) (1.058,1.241) (1.065,1.259) (1.072,1.276) (1.079,1.294) (1.086,1.312) (1.093,1.329) (1.100,1.347) (1.108,1.364) (1.115,1.382) (1.123,1.399) (1.130,1.417) (1.138,1.434) (1.145,1.452) (1.153,1.469) (1.161,1.487) (1.168,1.504) (1.176,1.521) (1.184,1.539) (1.192,1.556) (1.200,1.573) (1.208,1.590) (1.216,1.608) (1.224,1.625) (1.232,1.642) (1.241,1.659) (1.249,1.676) (1.257,1.694) (1.266,1.711) (1.274,1.728) (1.282,1.745) (1.291,1.762) (1.299,1.779) (1.308,1.796) (1.317,1.813) (1.325,1.830) (1.334,1.847) (1.343,1.864) (1.351,1.881) (1.360,1.897) (1.369,1.914) (1.378,1.931) (1.387,1.948) (1.396,1.965) (1.404,1.982) (1.413,1.998) (1.422,2.015) (1.432,2.032) (1.441,2.049) (1.450,2.065) (1.459,2.082) (1.468,2.099) (1.477,2.115) (1.486,2.132) (1.496,2.149) (1.505,2.165) (1.514,2.182) (1.524,2.198) (1.533,2.215) (1.542,2.231) (1.552,2.248) (1.561,2.265) (1.571,2.281) (1.580,2.298) (1.590,2.314) (1.599,2.331) (1.609,2.347) (1.618,2.363) (1.628,2.380) (1.637,2.396) (1.647,2.413) (1.657,2.429) (1.666,2.445) (1.676,2.462) (1.686,2.478) (1.696,2.495) (1.705,2.511) (1.715,2.527) (1.725,2.543) (1.735,2.560) (1.745,2.576) (1.754,2.592) (1.764,2.609) (1.774,2.625) (1.784,2.641) (1.794,2.657) (1.804,2.674) (1.814,2.690) (1.824,2.706) (1.834,2.722) (1.844,2.738) (1.854,2.755) (1.864,2.771) (1.874,2.787) (1.884,2.803) (1.894,2.819) (1.904,2.835) (1.915,2.851) (1.925,2.867) (1.935,2.884) (1.945,2.900) (1.955,2.916) (1.965,2.932) (1.975,2.948) (1.986,2.964) (1.996,2.980) (2.006,2.996) (2.016,3.012) (2.027,3.028) (2.037,3.044) (2.047,3.060) (2.058,3.076) (2.068,3.092) (2.078,3.108) (2.088,3.124) (2.099,3.140) (2.109,3.156) (2.120,3.172) (2.130,3.188) (2.140,3.204) (2.151,3.220) (2.161,3.236) (2.172,3.252) (2.182,3.268) (2.192,3.284) (2.203,3.299) (2.213,3.315) (2.224,3.331) (2.234,3.347) (2.245,3.363) (2.255,3.379) (2.266,3.395) (2.276,3.411) (2.287,3.426) (2.297,3.442) (2.308,3.458) (2.318,3.474) (2.329,3.490) (2.339,3.506) (2.350,3.522) (2.361,3.537) (2.371,3.553) (2.382,3.569) (2.392,3.585) (2.403,3.601) (2.413,3.616) (2.424,3.632) (2.435,3.648) (2.445,3.664) (2.456,3.680) (2.467,3.695) (2.477,3.711) (2.488,3.727) (2.499,3.743)};
  \addplot[sepc, line width=.8pt, gray!70] coordinates {(0.000,0.800) (0.018,0.800) (0.037,0.801) (0.055,0.802) (0.073,0.803) (0.092,0.804) (0.110,0.806) (0.128,0.808) (0.146,0.811) (0.164,0.813) (0.182,0.817) (0.200,0.820) (0.218,0.824) (0.236,0.828) (0.254,0.832) (0.272,0.837) (0.289,0.842) (0.307,0.847) (0.324,0.853) (0.342,0.858) (0.359,0.864) (0.376,0.871) (0.393,0.877) (0.410,0.884) (0.427,0.891) (0.444,0.898) (0.461,0.906) (0.477,0.914) (0.494,0.922) (0.510,0.930) (0.527,0.938) (0.543,0.947) (0.559,0.956) (0.575,0.965) (0.591,0.974) (0.607,0.983) (0.622,0.993) (0.638,1.002) (0.653,1.012) (0.669,1.022) (0.684,1.032) (0.699,1.043) (0.714,1.053) (0.729,1.064) (0.744,1.074) (0.759,1.085) (0.774,1.096) (0.788,1.107) (0.803,1.118) (0.817,1.130) (0.831,1.141) (0.846,1.153) (0.860,1.164) (0.874,1.176) (0.888,1.188) (0.902,1.200) (0.916,1.212) (0.929,1.224) (0.943,1.236) (0.957,1.248) (0.970,1.261) (0.984,1.273) (0.997,1.286) (1.010,1.298) (1.024,1.311) (1.037,1.324) (1.050,1.336) (1.063,1.349) (1.076,1.362) (1.089,1.375) (1.102,1.388) (1.115,1.401) (1.127,1.414) (1.140,1.428) (1.153,1.441) (1.165,1.454) (1.178,1.468) (1.191,1.481) (1.203,1.494) (1.215,1.508) (1.228,1.521) (1.240,1.535) (1.252,1.549) (1.265,1.562) (1.277,1.576) (1.289,1.590) (1.301,1.604) (1.313,1.617) (1.325,1.631) (1.337,1.645) (1.349,1.659) (1.361,1.673) (1.373,1.687) (1.384,1.701) (1.396,1.715) (1.408,1.729) (1.420,1.743) (1.431,1.757) (1.443,1.771) (1.455,1.786) (1.466,1.800) (1.478,1.814) (1.489,1.828) (1.501,1.842) (1.512,1.857) (1.524,1.871) (1.535,1.885) (1.547,1.900) (1.558,1.914) (1.569,1.929) (1.581,1.943) (1.592,1.957) (1.603,1.972) (1.615,1.986) (1.626,2.001) (1.637,2.015) (1.648,2.030) (1.659,2.044) (1.670,2.059) (1.682,2.074) (1.693,2.088) (1.704,2.103) (1.715,2.117) (1.726,2.132) (1.737,2.147) (1.748,2.161) (1.759,2.176) (1.770,2.191) (1.781,2.205) (1.792,2.220) (1.803,2.235) (1.814,2.249) (1.824,2.264) (1.835,2.279) (1.846,2.294) (1.857,2.308) (1.868,2.323) (1.879,2.338) (1.889,2.353) (1.900,2.368) (1.911,2.383) (1.922,2.397) (1.933,2.412) (1.943,2.427) (1.954,2.442) (1.965,2.457) (1.975,2.472) (1.986,2.487) (1.997,2.501) (2.008,2.516) (2.018,2.531) (2.029,2.546) (2.040,2.561) (2.050,2.576) (2.061,2.591) (2.071,2.606) (2.082,2.621) (2.093,2.636) (2.103,2.651) (2.114,2.666) (2.124,2.681) (2.135,2.696) (2.145,2.711) (2.156,2.726) (2.167,2.741) (2.177,2.756) (2.188,2.771) (2.198,2.786) (2.209,2.801) (2.219,2.816) (2.230,2.831) (2.240,2.846) (2.251,2.861) (2.261,2.876) (2.272,2.891) (2.282,2.906) (2.292,2.921) (2.303,2.936) (2.313,2.951) (2.324,2.966) (2.334,2.981) (2.345,2.996) (2.355,3.011) (2.365,3.027) (2.376,3.042) (2.386,3.057) (2.397,3.072) (2.407,3.087) (2.418,3.102) (2.428,3.117) (2.438,3.132) (2.449,3.147) (2.459,3.162) (2.469,3.177) (2.480,3.193) (2.490,3.208) (2.501,3.223) (2.511,3.238) (2.521,3.253) (2.532,3.268) (2.542,3.283) (2.552,3.299) (2.563,3.314) (2.573,3.329) (2.583,3.344) (2.594,3.359) (2.604,3.374) (2.614,3.389) (2.625,3.405) (2.635,3.420) (2.645,3.435) (2.655,3.450) (2.666,3.465) (2.676,3.480) (2.686,3.495) (2.697,3.511) (2.707,3.526) (2.717,3.541) (2.728,3.556) (2.738,3.571)};
  \node[saddle] at (axis cs:0.8,0) {};
  \node[saddle] at (axis cs:0,0.8) {};
\end{axis}
\begin{axis}[at={(7.0cm,0)}, width=6.6cm, height=6.6cm, scale only axis=false, axis equal image,
  xmin=0, xmax=4.5, ymin=0, ymax=4.5, axis lines=left, axis line style={line width=1.1pt, -},
  xtick=\empty, ytick=\empty, xlabel={$x$}, ylabel={$y$}, x label style={at={(axis description cs:1,0)}, anchor=north east}, y label style={at={(axis description cs:0,1)}, anchor=north east, rotate=-90},
  tick label style={font=\scriptsize}, title={node: $(p,q)=(8,12)$, $n=22$}, title style={font=\small,yshift=4pt}, clip=false]
  \addplot[thetaline] coordinates {(0,0) (2.783,3.408)};
  \node[font=\scriptsize, gray!50!black, anchor=south west, rotate=50] at (axis cs:2.983,3.308) {$\vartheta=\vartheta_0$};
  \addplot[sepc] coordinates {(0.800,0.000) (0.800,0.019) (0.800,0.037) (0.801,0.056) (0.802,0.075) (0.803,0.093) (0.804,0.112) (0.806,0.131) (0.807,0.149) (0.809,0.168) (0.812,0.186) (0.814,0.205) (0.817,0.223) (0.819,0.242) (0.822,0.260) (0.826,0.279) (0.829,0.297) (0.833,0.315) (0.837,0.334) (0.841,0.352) (0.845,0.370) (0.849,0.388) (0.854,0.406) (0.859,0.425) (0.864,0.443) (0.869,0.461) (0.874,0.478) (0.880,0.496) (0.885,0.514) (0.891,0.532) (0.897,0.550) (0.903,0.567) (0.910,0.585) (0.916,0.602) (0.923,0.620) (0.929,0.637) (0.936,0.655) (0.943,0.672) (0.950,0.689) (0.957,0.707) (0.965,0.724) (0.972,0.741) (0.980,0.758) (0.987,0.775) (0.995,0.792) (1.003,0.809) (1.011,0.826) (1.019,0.843) (1.027,0.860) (1.035,0.876) (1.044,0.893) (1.052,0.910) (1.061,0.927) (1.069,0.943) (1.078,0.960) (1.087,0.976) (1.095,0.993) (1.104,1.009) (1.113,1.026) (1.122,1.042) (1.131,1.058) (1.141,1.075) (1.150,1.091) (1.159,1.107) (1.168,1.123) (1.178,1.139) (1.187,1.155) (1.197,1.172) (1.206,1.188) (1.216,1.204) (1.226,1.220) (1.235,1.236) (1.245,1.252) (1.255,1.267) (1.265,1.283) (1.275,1.299) (1.285,1.315) (1.294,1.331) (1.304,1.347) (1.315,1.362) (1.325,1.378) (1.335,1.394) (1.345,1.410) (1.355,1.425) (1.365,1.441) (1.376,1.456) (1.386,1.472) (1.396,1.488) (1.407,1.503) (1.417,1.519) (1.427,1.534) (1.438,1.550) (1.448,1.565) (1.459,1.581) (1.469,1.596) (1.480,1.612) (1.490,1.627) (1.501,1.642) (1.512,1.658) (1.522,1.673) (1.533,1.689) (1.544,1.704) (1.554,1.719) (1.565,1.734) (1.576,1.750) (1.587,1.765) (1.597,1.780) (1.608,1.796) (1.619,1.811) (1.630,1.826) (1.641,1.841) (1.651,1.856) (1.662,1.872) (1.673,1.887) (1.684,1.902) (1.695,1.917) (1.706,1.932) (1.717,1.947) (1.728,1.963) (1.739,1.978) (1.750,1.993) (1.761,2.008) (1.772,2.023) (1.783,2.038) (1.794,2.053) (1.805,2.068) (1.816,2.083) (1.827,2.098) (1.839,2.113) (1.850,2.128) (1.861,2.143) (1.872,2.158) (1.883,2.173) (1.894,2.188) (1.905,2.203) (1.917,2.218) (1.928,2.233) (1.939,2.248) (1.950,2.263) (1.961,2.278) (1.973,2.293) (1.984,2.308) (1.995,2.323) (2.006,2.338) (2.018,2.353) (2.029,2.368) (2.040,2.382) (2.052,2.397) (2.063,2.412) (2.074,2.427) (2.085,2.442) (2.097,2.457) (2.108,2.472) (2.119,2.487) (2.131,2.502) (2.142,2.516) (2.153,2.531) (2.165,2.546) (2.176,2.561) (2.188,2.576) (2.199,2.591) (2.210,2.605) (2.222,2.620) (2.233,2.635) (2.245,2.650) (2.256,2.665) (2.267,2.679) (2.279,2.694) (2.290,2.709) (2.302,2.724) (2.313,2.739) (2.324,2.753) (2.336,2.768) (2.347,2.783) (2.359,2.798) (2.370,2.812) (2.382,2.827) (2.393,2.842) (2.405,2.857) (2.416,2.871) (2.428,2.886) (2.439,2.901) (2.451,2.916) (2.462,2.930) (2.474,2.945) (2.485,2.960) (2.497,2.975) (2.508,2.989) (2.520,3.004) (2.531,3.019) (2.543,3.034) (2.554,3.048) (2.566,3.063) (2.577,3.078) (2.589,3.092) (2.600,3.107) (2.612,3.122) (2.623,3.136) (2.635,3.151) (2.646,3.166) (2.658,3.181) (2.670,3.195) (2.681,3.210) (2.693,3.225) (2.704,3.239) (2.716,3.254) (2.727,3.269) (2.739,3.283) (2.751,3.298) (2.762,3.313) (2.774,3.327) (2.785,3.342) (2.797,3.357) (2.808,3.371) (2.820,3.386) (2.832,3.401) (2.843,3.415) (2.855,3.430) (2.866,3.445) (2.878,3.459)};
  \addplot[sepc, line width=.8pt, gray!70] coordinates {(0.000,0.800) (0.018,0.800) (0.036,0.801) (0.055,0.802) (0.073,0.804) (0.091,0.806) (0.109,0.808) (0.127,0.811) (0.145,0.814) (0.162,0.818) (0.180,0.822) (0.198,0.827) (0.215,0.831) (0.233,0.837) (0.250,0.842) (0.267,0.848) (0.284,0.855) (0.301,0.861) (0.318,0.868) (0.335,0.875) (0.351,0.883) (0.368,0.891) (0.384,0.899) (0.400,0.907) (0.416,0.916) (0.432,0.924) (0.448,0.933) (0.464,0.943) (0.479,0.952) (0.495,0.962) (0.510,0.971) (0.525,0.981) (0.540,0.991) (0.555,1.002) (0.570,1.012) (0.585,1.023) (0.600,1.033) (0.615,1.044) (0.629,1.055) (0.643,1.066) (0.658,1.078) (0.672,1.089) (0.686,1.100) (0.700,1.112) (0.714,1.123) (0.728,1.135) (0.742,1.147) (0.756,1.159) (0.770,1.171) (0.783,1.183) (0.797,1.195) (0.811,1.207) (0.824,1.219) (0.838,1.231) (0.851,1.244) (0.864,1.256) (0.878,1.268) (0.891,1.281) (0.904,1.293) (0.917,1.306) (0.930,1.319) (0.943,1.331) (0.956,1.344) (0.969,1.357) (0.982,1.370) (0.995,1.382) (1.008,1.395) (1.021,1.408) (1.034,1.421) (1.046,1.434) (1.059,1.447) (1.072,1.460) (1.084,1.473) (1.097,1.486) (1.110,1.499) (1.122,1.512) (1.135,1.526) (1.147,1.539) (1.160,1.552) (1.172,1.565) (1.185,1.579) (1.197,1.592) (1.210,1.605) (1.222,1.618) (1.234,1.632) (1.247,1.645) (1.259,1.659) (1.271,1.672) (1.284,1.685) (1.296,1.699) (1.308,1.712) (1.320,1.726) (1.333,1.739) (1.345,1.753) (1.357,1.766) (1.369,1.780) (1.381,1.793) (1.394,1.807) (1.406,1.820) (1.418,1.834) (1.430,1.847) (1.442,1.861) (1.454,1.875) (1.466,1.888) (1.478,1.902) (1.490,1.915) (1.502,1.929) (1.515,1.943) (1.527,1.956) (1.539,1.970) (1.551,1.984) (1.563,1.997) (1.575,2.011) (1.587,2.025) (1.599,2.038) (1.610,2.052) (1.622,2.066) (1.634,2.079) (1.646,2.093) (1.658,2.107) (1.670,2.121) (1.682,2.134) (1.694,2.148) (1.706,2.162) (1.718,2.176) (1.730,2.189) (1.742,2.203) (1.753,2.217) (1.765,2.231) (1.777,2.245) (1.789,2.258) (1.801,2.272) (1.813,2.286) (1.825,2.300) (1.836,2.314) (1.848,2.327) (1.860,2.341) (1.872,2.355) (1.884,2.369) (1.896,2.383) (1.907,2.397) (1.919,2.411) (1.931,2.424) (1.943,2.438) (1.954,2.452) (1.966,2.466) (1.978,2.480) (1.990,2.494) (2.002,2.508) (2.013,2.521) (2.025,2.535) (2.037,2.549) (2.049,2.563) (2.060,2.577) (2.072,2.591) (2.084,2.605) (2.096,2.619) (2.107,2.633) (2.119,2.646) (2.131,2.660) (2.143,2.674) (2.154,2.688) (2.166,2.702) (2.178,2.716) (2.189,2.730) (2.201,2.744) (2.213,2.758) (2.225,2.772) (2.236,2.786) (2.248,2.800) (2.260,2.814) (2.271,2.827) (2.283,2.841) (2.295,2.855) (2.306,2.869) (2.318,2.883) (2.330,2.897) (2.341,2.911) (2.353,2.925) (2.365,2.939) (2.377,2.953) (2.388,2.967) (2.400,2.981) (2.412,2.995) (2.423,3.009) (2.435,3.023) (2.447,3.037) (2.458,3.051) (2.470,3.065) (2.481,3.079) (2.493,3.093) (2.505,3.107) (2.516,3.121) (2.528,3.135) (2.540,3.149) (2.551,3.163) (2.563,3.176) (2.575,3.190) (2.586,3.204) (2.598,3.218) (2.610,3.232) (2.621,3.246) (2.633,3.260) (2.645,3.274) (2.656,3.288) (2.668,3.302) (2.679,3.316) (2.691,3.330) (2.703,3.344) (2.714,3.358) (2.726,3.372) (2.738,3.386) (2.749,3.400) (2.761,3.414) (2.772,3.428) (2.784,3.442) (2.796,3.456) (2.807,3.470) (2.819,3.484) (2.830,3.498)};
  \node[saddle] at (axis cs:0.8,0) {};
  \node[saddle] at (axis cs:0,0.8) {};
\end{axis}
\end{tikzpicture}
\caption{The separatrix profiles of Proposition~\ref{prop:boundary-profile} for the two actions of Figure~\ref{fig:compactified-phase}: the profiles of the separatrices $Q_1\to P_0$ (solid) and $Q_4\to P_0$ (grey), normalized by a homothety so that they start at distance $0.8$ from the origin (white circles). Each meets a singular boundary ray orthogonally and goes to infinity in the direction of the minimal cone (dotted). By Theorem~\ref{thm:D} the ends cross the cone infinitely many times in the focus case (left; only the first crossing is visible at this scale) and finitely many times in the node case (right). The separatrices $P_1\to Q_2$ and $P_1\to Q_3$ give the same curves with the opposite orientation.}
\label{fig:separatrix-profiles}
\end{figure}

Two such profiles are shown in Figure~\ref{fig:separatrix-profiles}.

\begin{proof}
We consider the separatrix leaving $Q_1$; the other three cases are identical up to the symmetries $\delta\mapsto\delta+\pi$ (reversal of the orientation of the profile) and $\vartheta\mapsto\pi/d-\vartheta$ in the local analysis. Along the profile, \eqref{eq:polar-s-new} gives
\[
\frac{dr}{d\vartheta}=r\cot\delta,
\qquad
\frac{ds}{d\vartheta}=\frac{r}{\sin\delta},
\]
wherever $\sin\delta\neq0$. Near $Q_1$ the separatrix is a smooth curve tangent to the eigenvector $(3+\mu_-,-(n+1))$, so it can be written as $\delta=\pi/2+O(\vartheta)$ for $0<\vartheta\le\vartheta_1$, with $\vartheta_1>0$ small. Hence $\cot\delta=O(\vartheta)$ and $\sin\delta\ge1/2$ on this portion. Integrating the first identity from $\vartheta_1$ down to $0$ shows that $\log r$ has a finite limit as $\vartheta\to0^+$, i.e., $r\to r_b\in(0,+\infty)$. The second identity then shows that the arc length needed to reach $\vartheta=0$ is finite. Thus the profile reaches the boundary ray $\vartheta=0$ at the point at distance $r_b>0$ from the origin, and there $\alpha=\vartheta+\delta=\pi/2$, so its tangent is orthogonal to the ray $\vartheta=0$. At $Q_2$ one finds in the same way $\alpha\equiv\pi/2\pmod\pi$, while at $Q_3,Q_4$ one has $\vartheta=\pi/d$ and $\alpha\equiv\pi/d+\pi/2\pmod\pi$, which is orthogonal to the second boundary ray.

At the other end the separatrix converges to $P_0$ (or to $P_1$ in backward time) by \eqref{eq:separatrix-connections}, and the argument used in the proof of Theorem~\ref{thm:B} shows that $r\to+\infty$, that the arc length is infinite, and that $\vartheta\to\vartheta_0$.

Finally, by Corollary~\ref{cor:monotone} the function $\Lambda=V^{1/3}r\cos\delta$ is strictly increasing along the profile. Along the separatrix leaving $Q_1$, as the profile approaches the boundary ray $\vartheta=0$, $V=r^{n-2}v(\vartheta)\to0$ because $m_0\ge1$, while $r\to r_b$ and $|\cos\delta|\le1$; hence $\Lambda\to0$ at that end. Therefore $\Lambda>0$, i.e., $\cos\delta>0$, along the whole profile, and $r$ is strictly increasing. The same argument applies to the separatrix leaving $Q_4$, while along the separatrices $P_1\to Q_2$ and $P_1\to Q_3$, parametrized in the direction of the flow, $\Lambda$ increases to $0$, so that $\cos\delta<0$ and $r$ is strictly decreasing.
\end{proof}

\begin{remark}[The hypersurfaces generated by the separatrices]\label{rem:separatrix-hypersurfaces}
Let $\gamma$ be the profile of the separatrix leaving $Q_1$, parametrized by arc length on $(0,+\infty)$, with $\gamma(s)\to\xi_b=(r_b,0)$ as $s\to0^+$. Since the separatrix is tangent at $Q_1$ to the eigenvector $(3+\mu_-,-(n+1))$, with $\mu_-=m_0$, along it we have $\cos\delta=\frac{n+1}{m_0+3}\,\vartheta+O(\vartheta^2)$, while $\ell'(\vartheta)=m_0/\vartheta+O(1)$ by the expression of $\ell'$ in the proof of Lemma~\ref{lem:a-smooth}. Hence, by \eqref{eq:kappa-delta}, the curvature of the profile and the mean curvature $f=-2\kappa$ have finite limits at the singular orbit:
\begin{equation}\label{eq:f-singular-orbit}
\lim_{s\to0^+}\kappa(s)=-\frac{n-2-m_0}{(m_0+3)\,r_b},
\qquad
\lim_{s\to0^+}f(s)=\frac{2(n-2-m_0)}{(m_0+3)\,r_b}.
\end{equation}
On the other hand, $f\to0$ along the other end, because the bracket in \eqref{eq:kappa-delta} tends to $\ell'(\vartheta_0)=0$ at $P_0$ while $r\to+\infty$. The same holds for the separatrix leaving $Q_4$, with $m_0$ replaced by $\mu_+$ (the sign of $\kappa$ depends on the orientation), while the separatrices $P_1\to Q_2$ and $P_1\to Q_3$ give the same two curves with the opposite orientation.

\emph{The case $d=1$.} Here $n-2-m_0=0$, and in fact the separatrices are straight lines. Indeed, along a profile which reaches the axis $y=0$, the first integral {\eqref{eq:d1-first-integral}} must vanish, so that $\cos\alpha\equiv0$ and the profile is the half-line $x=\pm r_b$. The corresponding hypersurfaces {$\Sigma_\gamma$} are the hyperplanes $x=\pm r_b$ orthogonal to the axis of rotation{, with their point on the axis removed; this point is the fixed point of $G$ over $\xi_b$, and adding it one obtains the whole hyperplanes,} which are totally geodesic. Thus, for $d=1$, the separatrices do not produce proper biconservative hypersurfaces.

\emph{The case $d\ge2$.} Now $n-2-m_0=\sum_{i\ge1}m_i>0$, so that $f$ is not constant by \eqref{eq:f-singular-orbit}, and $\Sigma_\gamma$ is proper biconservative by Theorem~\ref{thm:A}. Since $r$ is strictly monotone, $\gamma$ is a radial graph over $(r_b,+\infty)$; in particular it is simple and its closure meets $\partial Q$ only at $\xi_b$. Hence $\Sigma_\gamma$, {which lies over the principal orbits,} is embedded, as in Remark~\ref{rem:radial-graph}, and its closure in $\R^n$ is $\Sigma_\gamma\cup G\cdot\xi_b$, where $G\cdot\xi_b$ is the singular orbit of dimension $n-2-m_0$ along which the $m_0$-dimensional spheres of the principal orbits collapse. If $\Sigma_\gamma$ extends smoothly across $G\cdot\xi_b$ (see Remark~\ref{rem:smooth-extension}), then $\overline\Sigma_\gamma=\Sigma_\gamma\cup G\cdot\xi_b$ is a properly embedded, hence complete, proper biconservative hypersurface, diffeomorphic to the total space of a vector bundle of rank $m_0+1$ over $G\cdot\xi_b$. Its mean curvature along $G\cdot\xi_b$ is given by \eqref{eq:f-singular-orbit}, and it has a single end, which is asymptotic in angle to the minimal cone $C_G$ and crosses it infinitely many times if and only if $\Delta_d<0$ (Theorem~\ref{thm:D}). For instance, for $SO(p+1)\times SO(q+1)$ the separatrix leaving $Q_1$ ends on the ray $y=0$, where the factor $\Sph^q$ collapses; then $f=2p/((q+3)r_b)$ along the singular orbit $\Sph^p(r_b)\times\{0\}$, and $\overline\Sigma_\gamma$ is diffeomorphic to $\Sph^p\times\R^{q+1}$. Since each saddle has a unique separatrix, modulo homotheties every action gives only two hypersurfaces of this kind, one for each type of singular orbit.
\end{remark}

We can now collect the dynamical results of Sections~\ref{sec:compactified} and~\ref{sec:global-phase}.

\begin{figure}[!htbp]
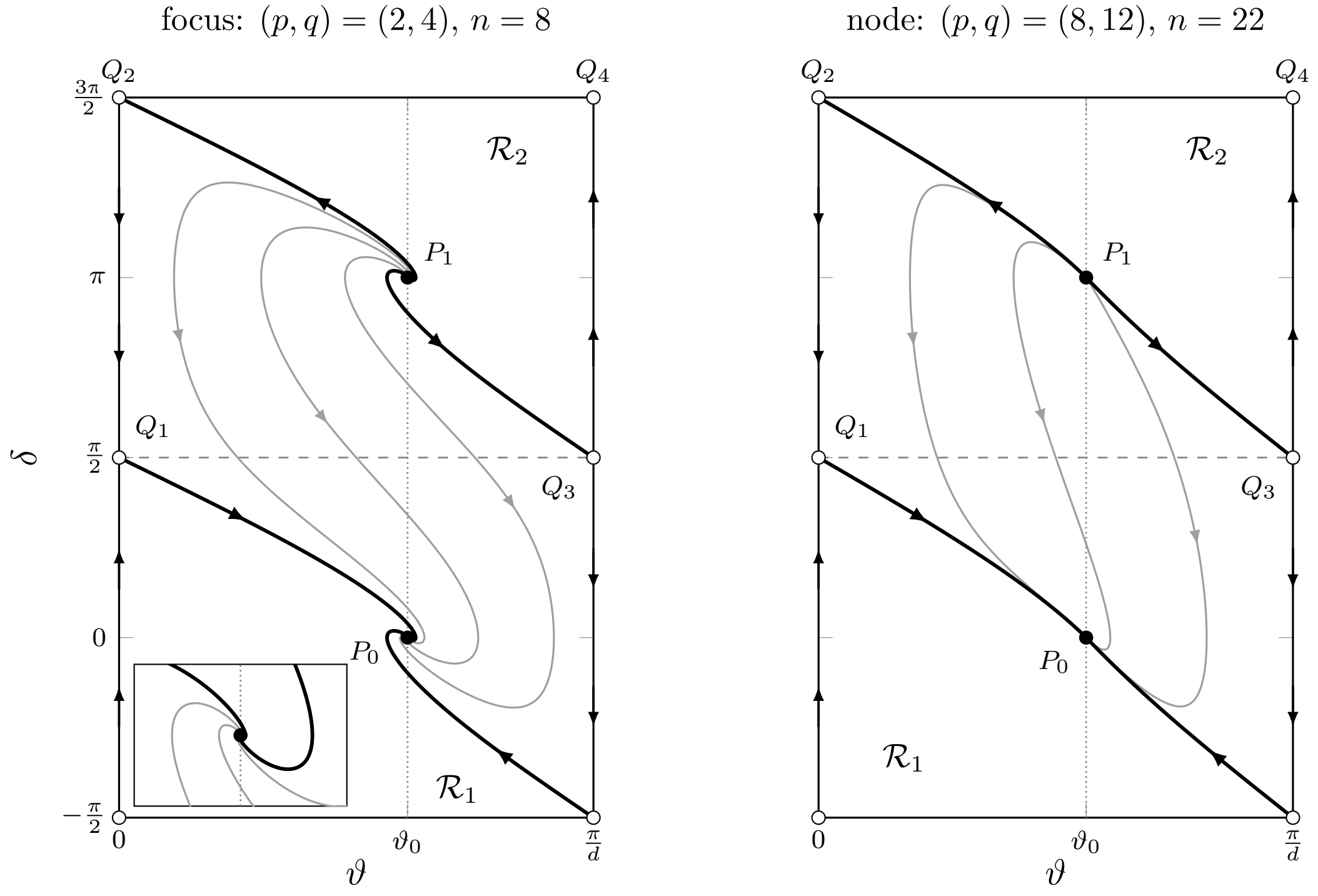

\centering

\caption{The compactified phase portrait of Theorem~\ref{thm:C} for $SO(p+1)\times SO(q+1)$ ($d=2$), in a focus case (left, $\Delta_2=-63$) and in a node case (right, $\Delta_2=49$), on the fundamental domain $-\pi/2\le\delta\le3\pi/2$, whose horizontal sides are identified. Thick black curves: the separatrices $Q_1\to P_0$, $Q_4\to P_0$, $P_1\to Q_2$, $P_1\to Q_3$. Thin grey curves: the trajectories through $(\vartheta_*,\pi/2)$, $\vartheta_*=\tfrac14\tfrac\pi d,\tfrac12\tfrac\pi d,\tfrac34\tfrac\pi d$, which give profiles of Theorem~\ref{thm:B}. Arrows on the vertical sides: the flow on the invariant boundary. Inset: a neighbourhood of $P_0$ in the focus case 
magnified about ten times. The trajectories spiral into $P_0$ (and, in backward time, into $P_1$), but each half-turn reduces the size of the spiral by a factor of about $1/35$, so that at the scale of the main picture only the first approach is visible.}
\label{fig:compactified-phase}
\end{figure}

\begin{proof}[{\bf Proof of Theorem~\ref{thm:C}}]
The phase portrait described in the statement is shown, for $d=2$, in Figure~\ref{fig:compactified-phase}. By Lemma~\ref{lem:a-smooth} the compactified vector field \eqref{eq:compactified} extends smoothly to $\mathcal C$, and its boundary circles are invariant because $h$ vanishes there. Its equilibria are listed in Proposition~\ref{prop:singular-points}, which also shows that the $Q_i$ are hyperbolic saddles. The points $P_0$ and $P_1$ lie on the ray $\vartheta=\vartheta_0$ and have $\delta=0$ and $\delta=\pi$, i.e. the tangent to the profile is the outward or the inward radial direction: they represent the minimal cone traversed away from or towards the origin. At the boundary points $\cos\delta=0$, i.e. the profile is orthogonal to a boundary ray (Remark~\ref{rem:boundary-meaning}). This proves~(i).

At $P_0$ the linearization of \eqref{eq:compactified} is $h(\vartheta_0)J_0$, with $J_0$ given in \eqref{eq:Jac-new}; since the vector field vanishes at $P_0$, the factor $h$ does not contribute further terms. The matrix $J_0$ has trace $-(n+1)<0$ and determinant $-3\ell''(\vartheta_0)>0$ by \eqref{eq:ell-second-new}, so both its eigenvalues have negative real part and $P_0$ is a hyperbolic sink. At $P_1$ the linearization is $-h(\vartheta_0)J_0$, because $\cos\delta$ and $\sin\delta$ change sign; hence $P_1$ is a hyperbolic source. This proves~(ii).

Statement~(iii) is Theorem~\ref{thm:global-trajectory}, and the identification with the profiles of Theorem~\ref{thm:B} follows from the initial conditions \eqref{eq:profile-init-new}; the additional assertions in (iii) are Corollary~\ref{cor:general-limits}. Statement~(iv) is Proposition~\ref{prop:separatrices}, together with Proposition~\ref{prop:boundary-profile} and Remark~\ref{rem:separatrix-hypersurfaces} for its geometric interpretation.
\end{proof}

\begin{remark}\label{rem:smooth-extension}
Orthogonality to a singular boundary ray is the natural orbit-space condition for a smooth extension across the corresponding singular orbit. For the $d=2$ action this is the familiar condition in the $O(m)\times O(n)$ setting. For the general Hsiang--Lawson actions, a complete classification of which separatrix profiles reconstruct to globally smooth hypersurfaces across the singular orbit requires the usual representation-dependent smoothness conditions and will be treated separately.
\end{remark}

\section{Focus versus node: oscillation and crossings of the minimal cone}\label{sec:focus-node}

The global analysis proves that both ends of the complete profiles converge to the minimal-cone direction, but it does not say how this convergence takes place. This is a local question and it is determined by the linearization at $P_0$ or $P_1$.

We now linearize at the minimal cone. For the raw system \eqref{eq:raw-planar}, the Jacobian at $P_0$ is
\begin{equation}\label{eq:Jac-new}
J_0=
\begin{pmatrix}
0&3\\
\ell''(\vartheta_0)&-(n+1)
\end{pmatrix}.
\end{equation}
The compactified system has the same linearization up to multiplication by the positive number $h(\vartheta_0)$, so the dynamical type is unchanged. The characteristic polynomial is
\begin{equation}\label{eq:char-new}
\lambda^2+(n+1)\lambda-3\ell''(\vartheta_0)=0.
\end{equation}

At this stage the answer still appears to depend on the full multiplicity pattern of the action. The next identity is therefore the key simplification behind the universal threshold.

\begin{lemma}\label{lem:ell2-universal-new}
For every action in Table~\ref{tab:actions},
\begin{equation}\label{eq:ell2-universal-new}
\ell''(\vartheta_0)=-d(n-2).
\end{equation}
\end{lemma}

\begin{proof}
For $d=1$, $\ell=(n-2)\log\sin\vartheta$ and $\vartheta_0=\pi/2$. For $d=2$, writing $m_0+m_1=n-2$ and using the critical-point relation gives $\ell''(\vartheta_0)=-2(n-2)$. For $d=3$ and $d=6$ all multiplicities are equal, and the product formula for the sine gives, up to an additive constant,
\[
\ell(\vartheta)=\frac{n-2}{d}\log|\sin(d\vartheta)|,
\]
which yields the formula. For $d=4$, the multiplicities alternate, say $m_0=m_2=a$ and $m_1=m_3=b$, and up to a constant
\[
\ell(\vartheta)=a\log\sin(2\vartheta)+b\log\cos(2\vartheta).
\]
Using the critical-point relation $a\cot(2\vartheta_0)=b\tan(2\vartheta_0)$ gives
\[
\ell''(\vartheta_0)=-8(a+b)=-4(n-2).
\]
\end{proof}

Consequently the discriminant is
\begin{equation}\label{eq:disc-new}
\Delta_d=(n+1)^2-12d(n-2).
\end{equation}

We can now prove Theorem~\ref{thm:D}. The proof turns a local dynamical classification into a geometric statement because crossing the vertical line $\vartheta=\vartheta_0$ is exactly crossing the minimal-cone ray in the orbit space.

\begin{proof}[{\bf Proof of Theorem~\ref{thm:D}}]
By \eqref{eq:char-new} and Lemma~\ref{lem:ell2-universal-new}, the eigenvalues of $J_0$ are the roots of $\lambda^2+(n+1)\lambda+3d(n-2)=0$, whose discriminant is $\Delta_d$; at $P_1$ the eigenvalues are their opposites. Thus $P_0$ is a spiral sink and $P_1$ a spiral source if $\Delta_d<0$, while both are nodes if $\Delta_d\ge0$. {Since the vector field is real analytic, we can use the classical description of the trajectories near a hyperbolic equilibrium of a perturbed linear system in the plane \cite[Chapter~VIII]{Hartman}, which, unlike the topological equivalence given by the Hartman--Grobman theorem, also controls the directions of approach: near a focus the polar angle around the equilibrium tends to $\pm\infty$ along every nonstationary trajectory, while every trajectory converging to a node, proper or improper, is eventually tangent to an eigendirection of the linearization.} In the focus case every nonstationary trajectory sufficiently close to $P_0$ or $P_1$ winds around the equilibrium, hence crosses the vertical line $\vartheta=\vartheta_0$ infinitely often. This is illustrated, for $d=2$ and $p=q=2$, in the bottom panel of Figure~\ref{fig:thmB}. In the orbit space this line corresponds exactly to the minimal-cone ray. In the node case, including the improper node $\Delta_d=0$, the eigenvalues are real and every eigenvector $(3,\lambda)$ of $J_0$ has nonzero $\vartheta$-component. Hence every trajectory converging to $P_0$ is eventually tangent to a direction transversal to the line $\vartheta=\vartheta_0$, so that $\vartheta-\vartheta_0$ eventually has constant sign. Therefore only finitely many crossings can occur near each end; the same argument applies at $P_1$ in backward time. {Moreover, all the crossings are transversal: on the line $\vartheta=\vartheta_0$ we have $\vartheta'=3h(\vartheta_0)\sin\delta$, which vanishes only at $P_0$ and $P_1$. Hence the crossings are isolated in the compactified time, only finitely many of them occur in every compact time interval, and in the node case each end of a trajectory converging to $P_0$ or $P_1$ crosses the line $\vartheta=\vartheta_0$ only finitely many times.} These local statements apply to the ends of the profiles of Theorem~\ref{thm:B}, which converge to $P_0$ and $P_1$, and to the ends at infinity of the separatrices of Theorem~\ref{thm:C}. Since $r\to+\infty$ along these ends, the corresponding crossings with the minimal cone, i.e. the principal orbits along which $\Sigma$ meets $C_G$, go off to infinity. The explicit thresholds for the various values of $d$ are listed below.
\end{proof}

For convenience we spell out the consequences for the five values of $d$:
\begin{itemize}[leftmargin=2em]
\item $d=1$: $\Delta_1=(n-5)^2$; the rotational case is nonoscillatory, with a repeated eigenvalue (improper node) at $n=5$.
\item $d=2$: the profile is oscillatory for $4\le n\le19$ and nonoscillatory for $n\ge20$.
\item $d=3$: all the cases in Table~\ref{tab:actions}, namely $n=5,8,14,26$, are oscillatory.
\item $d=4$: here every admissible dimension is even, and $\Delta_4=n^2-46n+97$ is negative for $3\le n\le43$ and positive for $n\ge44$. Hence the profile is oscillatory for the admissible dimensions $n\le42$ and nonoscillatory for the admissible dimensions $n\ge44$.
\item $d=6$: both cases in Table~\ref{tab:actions}, namely $n=8,14$, are oscillatory.
\end{itemize}

{In the oscillatory regime $\Delta_d<0$, Theorem~\ref{thm:D} shows that every end of a complete profile crosses the minimal cone infinitely many times, but it says nothing about where these crossings occur. Since $P_0$ and $P_1$ are hyperbolic foci, the linearization determines the oscillation quantitatively: the next result shows that, in the radial variable, the crossings are asymptotically in geometric progression, with a ratio depending only on $n$ and $d$, and that their amplitude decays like a fixed power of $r$.

\begin{proposition}[Spacing of cone crossings]\label{thm:cone-spacing}
Suppose that
\[
\Delta_d=(n+1)^2-12d(n-2)<0.
\]
Let $\gamma$ be a profile from Theorem~\ref{thm:B} or a separatrix from
Theorem~\ref{thm:C}, and consider an end of $\gamma$ tending to infinity.
There are constants $A,B\in\mathbb R$, with $A^2+B^2>0$, such that
\begin{equation}\label{eq:cone-spacing-expansion}
\vartheta(r)-\vartheta_0
=r^{-(n+1)/6}\left(
A\cos\left(\frac{\sqrt{-\Delta_d}}6\log r\right)
+B\sin\left(\frac{\sqrt{-\Delta_d}}6\log r\right)
+o(1)\right).
\end{equation}
The remainder is $o(1)$ also after differentiation with respect to
$\log r$. If $r_k$ are the radii of consecutive intersections of this
end with the minimal-cone ray $\vartheta=\vartheta_0$, in increasing
order, then
\begin{equation}\label{eq:cone-spacing-ratio}
\lim_{k\to\infty}\frac{r_{k+1}}{r_k}
=\exp\left(\frac{6\pi}{\sqrt{-\Delta_d}}\right).
\end{equation}
Moreover, the number $N(R)$ of these intersections between a fixed
sufficiently large radius and $R$ satisfies
\begin{equation}\label{eq:cone-spacing-count}
N(R)=\frac{\sqrt{-\Delta_d}}{6\pi}\log R+O(1).
\end{equation}
\end{proposition}

\begin{proof}
Consider first an end converging to $P_0=(\vartheta_0,0)$.
We use the time variable of \eqref{eq:raw-planar}. Set
\[
z=\begin{pmatrix}\vartheta-\vartheta_0\\ \delta\end{pmatrix}.
\]
Equations \eqref{eq:Jac-new}--\eqref{eq:disc-new} give
\[
z'=J_0z+F(z),\qquad F(z)=O(|z|^2),
\]
where the eigenvalues of $J_0$ are $-\alpha\pm i\beta$ with
\[
\alpha=\frac{n+1}{2},\qquad
\beta=\frac{\sqrt{-\Delta_d}}2.
\]
Choose a real invertible matrix $S$ such that
\[
S^{-1}J_0S=-\alpha I+\beta J,
\qquad J=\begin{pmatrix}0&-1\\1&0\end{pmatrix},
\]
and write $u=S^{-1}z$. Then
\[
u'=(-\alpha I+\beta J)u+G(u),\qquad G(u)=O(|u|^2).
\]
Since $u(t)\to0$, exponential stability gives, for every sufficiently
small $\varepsilon>0$,
\[
|u(t)|\leq C_\varepsilon e^{-(\alpha-\varepsilon)t}
\]
for large $t$. In particular $\int^\infty|u(t)|\,dt<\infty$.
The solution is nonstationary, so $u(t)\neq0$ and
\[
\frac{d}{dt}\bigl(\log|u(t)|+\alpha t\bigr)
=\frac{\langle u(t),G(u(t))\rangle}{|u(t)|^2}
=O(|u(t)|).
\]
It follows that $e^{\alpha t}|u(t)|$ has a positive finite limit.
Furthermore,
\[
w(t)=e^{\alpha t}e^{-\beta Jt}u(t)
\]
satisfies
\[
|w'(t)|\leq Ce^{\alpha t}|u(t)|^2
\leq C_\varepsilon e^{-(\alpha-2\varepsilon)t}.
\]
Taking $2\varepsilon<\alpha$, we obtain $w(t)\to w_\infty\neq0$
and $w'(t)\to0$. Thus
\[
z(t)=e^{-\alpha t}S e^{\beta Jt}
\bigl(w_\infty+o(1)\bigr).
\]
The first component of the leading term is a nonzero sinusoid.
Indeed, if it vanished identically, the first equation of the
linearized system, $z_1'=3z_2$, would imply $w_\infty=0$.
Consequently
\begin{equation}\label{eq:cone-spacing-time}
\vartheta(t)-\vartheta_0
=e^{-\alpha t}\bigl(C\cos(\beta t)+D\sin(\beta t)+o(1)\bigr),
\qquad C^2+D^2>0,
\end{equation}
with the corresponding differentiated expansion.

{By \eqref{eq:polar-s-new} and \eqref{eq:tparam-new}, the radial equation is}
\[
(\log r)'=3\cos\delta=3+O(e^{-2\alpha t}).
\]
Hence
\[
\log r=3t+c+O(e^{-2\alpha t})
\]
for some constant $c$. Substitution in
\eqref{eq:cone-spacing-time} gives
\eqref{eq:cone-spacing-expansion}, including its differentiated
form.

Set $x=\log r$ and $\omega=\sqrt{-\Delta_d}/6$. Then
\[
H(x):=e^{(n+1)x/6}\bigl(\vartheta(e^x)-\vartheta_0\bigr)
=K\cos(\omega x+\phi)+\eta(x),
\]
where $K>0$ and $\eta(x),\eta'(x)\to0$. For sufficiently large
$x$, each small fixed neighbourhood of a zero of the leading cosine
contains exactly one zero of $H$, and there are no other zeros:
existence follows from the change of sign, while uniqueness follows
from the fixed sign of $H'$ in that neighbourhood. Therefore the
successive zeros satisfy
\[
x_k=\frac{(k+\tfrac12)\pi-\phi}{\omega}+o(1),
\]
with a suitable indexing. Thus
\[
x_{k+1}-x_k\to\frac{\pi}{\omega},\qquad
N(R)=\frac{\omega}{\pi}\log R+O(1),
\]
which proves \eqref{eq:cone-spacing-ratio} and
\eqref{eq:cone-spacing-count}.

For an end converging to $P_1=(\vartheta_0,\pi)$ in backward
time, reverse time and replace $\delta$ by $\delta-\pi$.
The transformed equations have the same linearization and radial
growth as at $P_0$. The same argument applies.
\end{proof}

{\begin{remark}\label{rem:spacing-example}
For the profiles of Figure~\ref{fig:thmB}, i.e., $d=2$ and $p=q=2$, we have $n=6$ and $\Delta_2=49-96=-47$. By Proposition~\ref{thm:cone-spacing}, $\vartheta-\vartheta_0$ decays like $r^{-7/6}$, so that the signed distance $r\sin(\vartheta-\vartheta_0)$ from the minimal-cone ray, plotted in the bottom panel, decays like $r^{-1/6}$. Consecutive crossings with the minimal cone are asymptotically spaced in $\log r$ by $6\pi/\sqrt{47}\approx2.75$, that is, $r_{k+1}/r_k\to e^{6\pi/\sqrt{47}}\approx15.6$, and $N(R)\approx0.36\log R$. This agrees with the bottom panel of Figure~\ref{fig:thmB}, where the distance between consecutive zeros is already about $2.75$.
\end{remark}}

}

\section{A detailed analysis of the case $d=2$}\label{sec:d2}

For the product action, the compactified system agrees with the one studied in \cite{MOR-BC}. We identify this agreement, correct the incompleteness assertion there, and give an independent quantitative proof of completeness. The authors S.M., C.O. and A.R. wish to thank Jule Kalbertod for pointing out this error.

We now compare our general theory with the family studied in \cite{MOR-BC}. In that paper the ambient space was written as
\[
\R^n=\R^{p+1}\times\R^{q+1},
\qquad n=p+q+2,
\]
and the equivariant immersion was
\begin{equation}\label{eq:d2imm-new}
\varphi_{p,q}:\Sph^p\times\Sph^q\times I\longrightarrow
\R^{p+1}\times\R^{q+1},
\qquad
(w,z,s)\longmapsto\bigl(x(s)w,y(s)z\bigr).
\end{equation}
With $\dot x=\cos\alpha$ and $\dot y=\sin\alpha$, equation \eqref{eq:alpha-general} becomes
\begin{equation}\label{eq:d2eq-new}
3\dot\alpha+p\frac{\sin\alpha}{x}
-q\frac{\cos\alpha}{y}=0,
\end{equation}
which is exactly the proper biconservative equation used in \cite{MOR-BC}.

The normalization chosen in \eqref{eq:h-a} is convenient also because it makes the compactified system agree \emph{exactly}, rather than only up to a constant time-rescaling, with the vector field used in the $d=2$ analysis of \cite{MOR-BC}. Indeed, when $d=2$,
\[
h(\vartheta)=\frac12\sin(2\vartheta)=\sin\vartheta\cos\vartheta,
\]
and therefore the first equation in \eqref{eq:compactified} becomes
\begin{equation}\label{eq:d2theta-old}
\vartheta'
=3\sin\vartheta\cos\vartheta\,
\sin(\alpha-\vartheta),
\end{equation}
which is the first equation of \cite[(4.7)]{MOR-BC}. Moreover,
\[
\ell'(\vartheta)=q\cot\vartheta-p\tan\vartheta,
\]
so that
\[
a(\vartheta)
=q\cos^2\vartheta-p\sin^2\vartheta.
\]
Since $\alpha=\vartheta+\delta$, adding the two equations of \eqref{eq:compactified} gives
\begin{align*}
\alpha'
&=a(\vartheta)\cos\delta
-(n-2)h(\vartheta)\sin\delta\\
&=q\cos\alpha\cos\vartheta
-p\sin\alpha\sin\vartheta,
\end{align*}
where $n=p+q+2$. Hence the whole compactified system reduces exactly to the system (4.7) of \cite{MOR-BC}.

The minimal cone is
\begin{equation}\label{eq:d2minimal-new}
y=\sqrt{\frac qp}\,x,
\qquad
\vartheta_0=\arctan\sqrt{\frac qp}.
\end{equation}
Since here $n=p+q+2$, the focus condition $\Delta_2<0$ of Theorem~\ref{thm:D}, i.e. $n^2-22n+49<0$, is equivalent to
\[
p+q\le17,
\]
while the node condition is equivalent to $p+q\ge18$. Thus the threshold from the general linearization is precisely the one found in the detailed $d=2$ phase portrait.

\begin{remark}\label{rem:tripling-euclid}
By Corollary~\ref{cor:tripling}, for $(p,q)=(3p_0,3q_0)$ the profiles of Theorem~\ref{thm:B} are the profiles of $SO(p_0+1)\times SO(q_0+1)$-invariant minimal hypersurfaces of $\R^{p_0+q_0+2}$, studied in \cite{ABPRS}. The thresholds agree: $p+q\le17$ means $p_0+q_0\le5$, i.e., $(p_0+1)+(q_0+1)\le7$, which is the classical range in which the profiles of $O(m)\times O(n)$-invariant minimal hypersurfaces oscillate around the minimal cone.
\end{remark}

\begin{proof}[Proof of Corollary~\ref{cor:d2}]
For $d=2$ the group of Table~\ref{tab:actions} is $SO(p+1)\times SO(q+1)$, the minimal cone is \eqref{eq:d2minimal-new}, and the family is given by Theorem~\ref{thm:B}. The statement on the intersections with the minimal cone is Theorem~\ref{thm:D}, by the computation of the threshold above. The lower bounds for $x$ and $y$ are given by Proposition~\ref{prop:elementary} below. The comparison with \cite{MOR-BC} is explained in the following remark.
\end{proof}

\begin{remark}\label{rem:erratum-new}
In \cite[Theorem~1.1 and Theorem~4.12]{MOR-BC} we stated that the proper $SO(p+1)\times SO(q+1)$-invariant biconservative hypersurfaces produced by the phase-plane analysis are cones which are not complete, because their profile curves were claimed to reach the origin of the orbit space at a finite value $s_0$ of the arc length. This conclusion is incorrect.

The error is in the final part of the proof of \cite[Theorem~4.12]{MOR-BC}. There we used the functions
\[
I=y^{q/3}\cos\alpha,
\qquad
J=x^{p/3}\sin\alpha,
\]
which are increasing along every solution of \eqref{eq:d2eq-new}, and we argued that a profile defined for all $s$ smaller than a given $s^*$ would contain a point at which $I<0$ and $J<0$, and then, since $I$ and $J$ must both become positive, that the profile must pass through the origin. This last implication is not justified: in the interior of $Q$ the signs of $I$ and $J$ are the signs of $\cos\alpha$ and $\sin\alpha$, so that $I$ and $J$ can change sign through a rotation of the tangent vector, without $x$ or $y$ vanishing. In fact, for the solutions considered there the opposite holds: if $\alpha(0)=\vartheta_*+\pi/2$ with $0<\vartheta_*<\pi/2$, then $I(0)<0<J(0)$, and the monotonicity of $I$ and $J$ prevents the profile from reaching the origin (Proposition~\ref{prop:elementary}).

The corrected picture is the following, and it is the content of Corollary~\ref{cor:d2}. Choose $0<\vartheta_*<\pi/2$ and $r_0>0$, and impose
\[
x(0)=r_0\cos\vartheta_*,\qquad
 y(0)=r_0\sin\vartheta_*,\qquad
\alpha(0)=\vartheta_*+\frac\pi2.
\]
Then the corresponding solution of \eqref{eq:d2eq-new} is defined for all $s\in\R$; the radial function is bounded away from zero, has a unique minimum at $s=0$, and tends to $+\infty$ as $s\to\pm\infty$. Both ends are asymptotic in angle to the minimal cone \eqref{eq:d2minimal-new}. Hence one obtains an infinite family of complete proper $SO(p+1)\times SO(q+1)$-invariant biconservative hypersurfaces in $\R^{p+q+2}$. Moreover, if $p+q\le17$, the profile curve intersects the minimal-cone ray infinitely many times at each end, whereas if $p+q\ge18$ it is eventually nonoscillatory.

Accordingly, \cite[Theorem~1.1 and Theorem~4.12]{MOR-BC} must be replaced by the statement above: the words ``(cones)'' and ``None of these hypersurfaces is complete'' have to be deleted, and the profile curves are defined on the whole real line instead of on half-lines ending at the origin. Remarks~1.2 and~4.13 of \cite{MOR-BC}, which concern the alleged conical singularity at the origin, have to be deleted. The remaining results of \cite{MOR-BC}, in particular the oscillation threshold, the nonexistence of proper biharmonic $SO(p+1)\times SO(q+1)$-invariant hypersurfaces and the rotational case \cite[Theorem~6.5]{MOR-BC}, are not affected. The argument of the present paper shows that the same mechanism is not peculiar to $d=2$.
\end{remark}

The monotonicity of $I$ and $J$ actually gives an elementary and quantitative proof that these profiles are defined on the whole real line and that the corresponding hypersurfaces are complete, independent of the phase portrait, which uses only the following continuation criterion.

\begin{lemma}[Continuation criterion]\label{lem:criterion}
Let $(x(s),y(s),\alpha(s))$ be the maximal solution of the system \eqref{eq:alpha}--\eqref{eq:alpha-general} with $(x(0),y(0))$ in the interior of $Q$, defined on the maximal interval $(s_-,s_+)$, and put $\gamma=(x,y)$. If $s_+<+\infty$, then $\gamma(s)$ converges, as $s\to s_+$, to a point of the boundary of $Q$, i.e., of one of the two boundary rays (the vertex included). The analogous statement holds for $s_-$. In particular, if there exists $m>0$ such that the distance from $\gamma(s)$ to the boundary of $Q$ is at least $m$ for every $s\in(s_-,s_+)$, then $(s_-,s_+)=\R$.
\end{lemma}

\begin{proof}
Assume $s_+<+\infty$. Since $|\dot\gamma|=1$, the curve $\gamma$ is $1$-Lipschitz, hence $\gamma(s)$ converges to some point $\xi_+\in Q$ as $s\to s_+$. Assume, by contradiction, that $\xi_+$ lies in the interior of $Q$, and let $\overline B$ be a closed disc centred at $\xi_+$ and contained in the interior of $Q$. The right-hand side of \eqref{eq:alpha-general}, solved for $\dot\alpha$, is continuous on $\overline B\times\R$ and $2\pi$-periodic in $\alpha$, hence bounded there by a constant $C$. For $s$ close to $s_+$ we have $\gamma(s)\in\overline B$, hence $|\dot\alpha|\le C$, so that $\alpha(s)$ converges to some $\alpha_+$ as $s\to s_+$. The point $(\xi_+,\alpha_+)$ belongs to the domain of the smooth system \eqref{eq:alpha}--\eqref{eq:alpha-general}, and the solution through it extends the given one beyond $s_+$, which contradicts maximality. The last assertion follows, since under its assumption the limit point $\xi_+$ would lie in the interior of $Q$.
\end{proof}

\begin{proposition}\label{prop:elementary}
Let $p,q\ge1$, $0<\vartheta_*<\pi/2$ and $r_0>0$, and let $(x(s),y(s),\alpha(s))$ be the maximal solution of \eqref{eq:alpha} and \eqref{eq:d2eq-new} with
\[
x(0)=x_0=r_0\cos\vartheta_*,\qquad y(0)=y_0=r_0\sin\vartheta_*,\qquad \alpha(0)=\vartheta_*+\frac\pi2,
\]
defined on the maximal interval $(s_-,s_+)$. Then
\begin{align*}
x(s)&\ge x_0\cos^{3/p}\vartheta_*, & y(s)&\ge y_0, && 0\le s<s_+,\\
x(s)&\ge x_0, & y(s)&\ge y_0\sin^{3/q}\vartheta_*, && s_-<s\le0.
\end{align*}
Consequently $(s_-,s_+)=\R$, and the corresponding $SO(p+1)\times SO(q+1)$-invariant hypersurface is complete.
\end{proposition}

\begin{proof}
By \eqref{eq:alpha} and \eqref{eq:d2eq-new},
\[
\dot I=\frac p3\,y^{q/3}\,\frac{\sin^2\alpha}{x}\ge0,
\qquad
\dot J=\frac q3\,x^{p/3}\,\frac{\cos^2\alpha}{y}\ge0,
\]
and $I(0)=-y_0^{q/3}\sin\vartheta_*<0$, $J(0)=x_0^{p/3}\cos\vartheta_*>0$.

For $s\ge0$ we have $J(s)\ge J(0)>0$, hence $\sin\alpha(s)>0$; thus $\dot y>0$ and $y(s)\ge y_0$. Moreover $x^{p/3}\ge x^{p/3}\sin\alpha=J(s)\ge J(0)$, that is, $x(s)\ge x_0\cos^{3/p}\vartheta_*$.

For $s\le0$ we have $I(s)\le I(0)<0$, hence $\cos\alpha(s)<0$; thus $\dot x<0$ on $(s_-,0]$ and $x(s)\ge x_0$. Moreover $y^{q/3}\ge-y^{q/3}\cos\alpha=-I(s)\ge-I(0)$, that is, $y(s)\ge y_0\sin^{3/q}\vartheta_*$.

Since $Q$ is the closed quarter plane, the distance from $\gamma(s)=(x(s),y(s))$ to the boundary of $Q$ is $\min\{x(s),y(s)\}$, which is bounded below by a positive constant on $(s_-,s_+)$. Hence $(s_-,s_+)=\R$ by Lemma~\ref{lem:criterion}, and the hypersurface is complete by Lemma~\ref{lem:completeness}.
\end{proof}

Proposition~\ref{prop:elementary} gives existence on $\R$ and completeness directly in the variables of \cite{MOR-BC}; the asymptotic behaviour of the profiles, i.e., \eqref{eq:profile-asymptotics-new} and Theorem~\ref{thm:D}, still requires the analysis of Sections~\ref{sec:global-phase} and~\ref{sec:focus-node}.

\section{Embedded hyperspheres in $\Sph^n$ and rotationally invariant hypersurfaces}\label{sec:spheres}

In this section $N^n=\Sph^n$ and $G$ is one of the groups of Table~\ref{tab:actions}, acting on $\Sph^n\subset\R^n\oplus\R$ as in Section~\ref{subsec:sphere-setting}. By Theorem~\ref{thm:A}, the proper biconservative profiles are the geodesics of $V^{2/3}g_{\Sph^2}$ on the interior of the lune $Q_{\Sph}$. We first write this equation in polar coordinates, then we discuss two situations in which compact examples can be obtained: the transfer from Hsiang's minimal hyperspheres for $d=2$, and the rotational case $d=1$. {In this section, for functions along a profile curve parametrized by arc length $s$, a prime denotes the derivative with respect to $s$, while for a function of a single coordinate, such as $\ell(\vartheta)$, it denotes the derivative with respect to that coordinate.}

\subsection{The first-order system}\label{subsec:sphere-system}

We use the coordinates \eqref{eq:sphere-polar} on $Q_{\Sph}$, in which
\[
 g_{\Sph^2}=dr^2+\sin^2r\,d\vartheta^2,
 \qquad
 V=\sin^{n-2}r\;v(\vartheta),
\]
with $v$ as in \eqref{eq:ell-new}. Put $e_r=\partial_r$, $e_\vartheta=(\sin r)^{-1}\partial_\vartheta$, and let $\delta$ be the oriented angle from $e_r$ to the unit tangent of a unit-speed profile,
\[
 T=\cos\delta\,e_r+\sin\delta\,e_\vartheta .
\]
The frame $(\gamma,e_r,e_\vartheta)$ is positively oriented in $\R^3$, so that $\nu=\gamma\times T=-\sin\delta\,e_r+\cos\delta\,e_\vartheta$. The connection of $\Sph^2$ in this frame is given by $\nabla_{e_r}e_r=\nabla_{e_r}e_\vartheta=0$, $\nabla_{e_\vartheta}e_r=\cot r\,e_\vartheta$, $\nabla_{e_\vartheta}e_\vartheta=-\cot r\,e_r$, and a direct computation gives
\[
 \kappa=\delta'+\cot r\,\sin\delta,
 \qquad
 \nu(\log V)=-(n-2)\cot r\,\sin\delta+\frac{\ell'(\vartheta)}{\sin r}\cos\delta .
\]
Hence:

\begin{proposition}\label{prop:sphere-system}
The proper biconservative profile equation $3\kappa=\nu(\log V)$ in $Q_{\Sph}$ is equivalent to
\begin{equation}\label{eq:sphere-system}
 r'=\cos\delta,\qquad
 \vartheta'=\frac{\sin\delta}{\sin r},\qquad
 \delta'=\frac13\Bigl(\frac{\ell'(\vartheta)}{\sin r}\cos\delta-(n+1)\cot r\,\sin\delta\Bigr).
\end{equation}
\end{proposition}

System \eqref{eq:sphere-system} is obtained from the Euclidean system \eqref{eq:polar-s-new} by replacing $1/r$ with $1/\sin r$ in the angular coefficient and with $\cot r$ in the radial one. The meridian $\vartheta=\vartheta_0$, $\delta\equiv0$, is a solution: its lift is the intersection of $\Sph^n$ with the cone $C_G\times\R$, {which for $d\ge2$ is} a minimal hypersurface with two singular points at the poles. {For $d=1$, instead, $\vartheta_0=\pi/2$, the minimal ``cone'' $C_G$ is the hyperplane $\{0\}\times\R^{n-1}$, and the lift is the totally geodesic great hypersphere $\Sph^n\cap(\{0\}\times\R^{n-1}\times\R)$: the poles are fixed points of the action, but regular points of this hypersurface.}

\begin{remark}\label{rem:no-homothety}
Unlike the Euclidean equation, \eqref{eq:sphere-system} is not invariant under a one-parameter group acting on $r$, so it does not reduce to a planar autonomous system and the methods of Sections~\ref{sec:compactified}--\ref{sec:focus-node} do not apply directly. Two mechanisms nevertheless give global information: the coincidence of the reduced metrics with those of minimal hypersurfaces when the multiplicities are divisible by three (Section~\ref{subsec:transfer}), and, for $d=1$, an additional rotational symmetry of $V^{2/3}g_{\Sph^2}$ (Section~\ref{sec:rotational}).
\end{remark}

\begin{remark}\label{rem:height}
Although \eqref{eq:sphere-system} does not reduce to a planar system, Lemma~\ref{lem:weighted} still gives a monotone quantity, the spherical counterpart of Corollary~\ref{cor:monotone}. Indeed, the height $z=\cos r$ satisfies $\Hess_{g_{\Sph^2}}z=-z\,g_{\Sph^2}$ and, since $\grad z=-\sin r\,\partial_r$ and $V=\sin^{n-2}r\,v(\vartheta)$, also $g_{\Sph^2}(\grad z,\grad\log V)=-(n-2)\cos r=-(n-2)z$. Hence Lemma~\ref{lem:weighted} applies with $\lambda=-z$ and $\mu=-(n-2)z$, and along every proper biconservative profile in the interior of $Q_{\Sph}$ we obtain
\[
\bigl(V^{1/3}z'\bigr)'=-\frac{n+1}{3}V^{1/3}z,
\qquad z'=-\sin r\cos\delta .
\]
Hence, on every interval on which $z>0$, the function $V^{1/3}z'$ is strictly decreasing, so that every critical point of $z$ is a strict local maximum and $z$ has at most one critical point there; similarly, on every interval on which $z<0$, $z$ has at most one critical point, a strict local minimum. In particular, no closed profile contained in the interior of $Q_{\Sph}$ lies in one of the open hemispheres $\{z>0\}$ and $\{z<0\}$.
\end{remark}

\subsection{A transfer principle from the spherical Bernstein problem}\label{subsec:transfer}

In this subsection $d=2$, $G=SO(p+1)\times SO(q+1)$ and $n=p+q+2$. The invariant hypersurface with profile $\gamma(s)=(x(s),y(s),z(s))\in Q_{\Sph}$ is
\begin{equation}\label{eq:immersion-sphere}
 \Sph^p\times\Sph^q\times I\longrightarrow\Sph^{p+q+2}\subset\R^{p+1}\times\R^{q+1}\times\R,
 \qquad
 (w,u,s)\longmapsto\bigl(x(s)w,\ y(s)u,\ z(s)\bigr),
\end{equation}
and $V=x^py^q$. We call a profile curve $\gamma:[a,b]\to Q_{\Sph}$ \emph{admissible} if:
\begin{enumerate}[label=\textup{(\roman*)}]
 \item $\gamma((a,b))$ lies in the interior of $Q_{\Sph}$ and $\gamma$ is simple;
 \item $\gamma(a)$ belongs to the singular side $\{x=0\}$ and $\gamma(b)$ to the singular side $\{y=0\}${, and neither of them is one of the poles $(0,0,\pm1)$};
 \item $\gamma$ meets both singular sides orthogonally;
 \item the functions $x,y,z$ satisfy the usual even--odd smoothness conditions{, to all orders,} at the endpoints which guarantee that the invariant lift extends smoothly across the two singular orbits.
\end{enumerate}
For example, at an endpoint where $x=0$, in terms of a signed normal coordinate $t$ centred at that endpoint, {the functions $x$, $y$, $z$ must extend to smooth functions of $t$ in a neighbourhood of $0$, with $x$ odd and $y$, $z$ even; in particular} $x(t)=t+O(t^3)$, $y(t)=y_0+O(t^2)$, $z(t)=z_0+O(t^2)$, {although these finite-order expansions alone do not guarantee the smoothness of the lift}; and analogously, with the roles of $x$ and $y$ interchanged, at an endpoint where $y=0$. These endpoint conditions are conditions on the profile curve alone: they do not depend on the dimensions $p,q$ of the orbit factors.

\begin{proposition}\label{prop:admissible-lift}
Let $\gamma$ be an admissible profile curve and let $p,q\geq1$. Then the invariant hypersurface determined by $\gamma$ extends smoothly across the two singular orbits. If $\gamma$ is embedded in the orbit space, the extended hypersurface is embedded and diffeomorphic to $\Sph^p*\Sph^q\cong\Sph^{p+q+1}$.
\end{proposition}

\begin{proof}
On $(a,b)$ the hypersurface is $\Sph^p\times\Sph^q\times(a,b)$. At $a$ the factor $\Sph^p$ collapses, at $b$ the factor $\Sph^q$ collapses, and the even--odd conditions are the standard smoothness conditions for the corresponding collapse; since they involve only $x,y,z$, the same profile gives a smooth extension for every $p,q$. The resulting compact manifold is obtained from $\Sph^p\times\Sph^q\times[a,b]$ by collapsing $\Sph^p$ at one end and $\Sph^q$ at the other, i.e., it is the join $\Sph^p*\Sph^q\cong\Sph^{p+q+1}$. If the profile has no self-intersections and meets the singular strata only at its endpoints, distinct points cannot be identified by the lift, so the extended hypersurface is embedded.
\end{proof}

\begin{theorem}[Tripling transfer]\label{thm:tripling-transfer}
Let $p_0,q_0\geq1$, and suppose that an admissible profile curve $\gamma$ determines a non-equatorial $SO(p_0+1)\times SO(q_0+1)$-invariant minimal embedding $\Sph^{p_0+q_0+1}\to\Sph^{p_0+q_0+2}$. Then the same profile curve determines an $SO(3p_0+1)\times SO(3q_0+1)$-invariant proper biconservative embedding
\[
 \Sph^{3p_0+3q_0+1}\longrightarrow\Sph^{3p_0+3q_0+2}.
\]
\end{theorem}

\begin{proof}
By Corollary~\ref{cor:tripling}, $\gamma$ satisfies the proper biconservative profile equation for $(p,q)=(3p_0,3q_0)$, and by Proposition~\ref{prop:admissible-lift} its lift is an embedded hypersphere of dimension $3p_0+3q_0+1$. By Theorem~\ref{thm:A}, $f=-2\kappa$, so it remains to show that $\kappa$ is not constant. If it were, $\gamma$ would lie on a circle $\langle a,\gamma\rangle=c$ of $\Sph^2$, $a\in\R^3$. At the endpoint on $\{x=0\}$ the profile is orthogonal to the singular side, so its tangent is parallel to the $x$-direction, and differentiating gives $a_1=0$; similarly $a_2=0$ at the other endpoint. Hence $\gamma$ is a latitude $z\equiv c_0$, whose lift is a geodesic hypersphere of $\Sph^{p_0+q_0+2}$. Since this lift is minimal, it is the equator, contradicting the assumption.
\end{proof}

For the outer action $O(m)\times O(m)$ on $\R^m\oplus\R^m\oplus\R$, Hsiang's construction \cite{HsiangI} gives infinitely many non-equatorial embedded minimal hyperspheres $\Sph^{2m-1}\to\Sph^{2m}$ for $m=2,3$; see also \cite{CarlottoSchulz}. In our notation these correspond to $(p_0,q_0)=(1,1)$ and $(2,2)$, and Theorem~\ref{thm:tripling-transfer} gives the following result.

\begin{corollary}\label{cor:S7S8}
There exist infinitely many proper biconservative embeddings $\Sph^7\to\Sph^8$, invariant under $SO(4)\times SO(4)$, and infinitely many proper biconservative embeddings $\Sph^{13}\to\Sph^{14}$, invariant under $SO(7)\times SO(7)$.
\end{corollary}

{In \cite{HsiangII} Hsiang constructed, for each of the outer actions
\[
SO(2)\times SO(3)\ \text{on}\ \Sph^5,\quad SO(2)\times SO(4)\ \text{on}\ \Sph^6,\quad SO(2)\times SO(5)\ \text{and}\ SO(3)\times SO(4)\ \text{on}\ \Sph^7,
\]
acting on $\R^{p_0+1}\oplus\R^{q_0+1}\oplus\R$, infinitely many mutually non-congruent non-equatorial embedded minimal hyperspheres of generalized rotational type. In our notation they correspond to $(p_0,q_0)=(1,2)$, $(1,3)$, $(1,4)$ and $(2,3)$, and Theorem~\ref{thm:tripling-transfer} gives the following further examples.

\begin{corollary}\label{cor:HsiangII}
There exist infinitely many proper biconservative embedded hyperspheres
\begin{itemize}
\item $\Sph^{10}\to\Sph^{11}$, invariant under $SO(4)\times SO(7)$;
\item $\Sph^{13}\to\Sph^{14}$, invariant under $SO(4)\times SO(10)$;
\item $\Sph^{16}\to\Sph^{17}$, invariant under $SO(4)\times SO(13)$;
\item $\Sph^{16}\to\Sph^{17}$, invariant under $SO(7)\times SO(10)$.
\end{itemize}
\end{corollary}

\begin{remark}[The action of $U(1)\times\mathrm{Spin}(10)$]\label{rem:U1Spin10}
Among the cohomogeneity two actions with $d\ge3$, the tripling of the multiplicities relates only $S(U(2)\times U(3))$, with multiplicities $(3,2)$, and $U(1)\times\mathrm{Spin}(10)$, with multiplicities $(9,6)$ (Corollary~\ref{cor:tripling}). In \cite{HsiangII} Hsiang also constructed infinitely many non-equatorial embedded minimal hyperspheres $\Sph^{11}\to\Sph^{12}$ invariant under $S(U(2)\times U(3))$, acting on $\Sph^{12}\subset(\mathbb{C}^2\otimes_{\mathbb{C}}\mathbb{C}^3)\oplus\R$. By Corollary~\ref{cor:tripling}, their profile curves are also profiles of proper biconservative $U(1)\times\mathrm{Spin}(10)$-invariant hypersurfaces of $\Sph^{32}$. In order to conclude that these are smooth embedded hyperspheres $\Sph^{31}\to\Sph^{32}$, one needs an analogue of Proposition~\ref{prop:admissible-lift} for $d=4$, that is, a description of the smoothness conditions and of the topology of the lift at the singular orbits of $U(1)\times\mathrm{Spin}(10)$, which differ from those of $S(U(2)\times U(3))$. We do not pursue this here.
\end{remark}}

\begin{proof}[{\bf Proof of Theorem~\ref{thm:E}}]
The first statement is Theorem~\ref{thm:tripling-transfer}, and the second one follows from Corollaries~\ref{cor:S7S8} and~\ref{cor:HsiangII}.
\end{proof}

\subsection{The rotational case and closed examples}\label{sec:rotational}

We now show that the closed non-CMC biconservative hypersurfaces of spheres obtained in \cite{MP-CAG} (surfaces in $\Sph^3$) and in \cite{MOP-JMAA} (hypersurfaces in $\Sph^n$, $n\geq3$) are the case $d=1$ of the present framework. For $d=1$ the reduced metric is rotationally symmetric; the Clairaut integral of its geodesics then gives the variational {characterization} of the profile curves, the closure condition and the non-embeddedness in a direct way. This is the spherical counterpart of the first integral {\eqref{eq:d1-first-integral}} found in the Euclidean case {in the proof of Theorem~\ref{thm:B}}.

For $d=1$ the group $G=SO(n-1)$ acts on $\Sph^n\subset\R\times\R^{n-1}\times\R$ on the middle factor; up to orbit equivalence, this is the action of $O(n-1)$ which defines the rotational hypersurfaces of $\Sph^n$ studied in \cite{MP-CAG,MOP-JMAA}. The orbit space $Q_{\Sph}$ is the closed hemisphere $\{(x,y,z)\in\Sph^2:\ y\ge0\}$, whose boundary $y=0$ corresponds to the singular orbits, namely the points of the totally geodesic $\Sph^1=\Sph^n\cap(\R\times\{0\}\times\R)$, and $V=y^{n-2}$. Since the orbit directions carry the single principal curvature $k_0=-\nu(\log y)$, the proper biconservative profile equation \eqref{eq:kappaV} reads
\begin{equation}\label{eq:rot-reduced}
 3\kappa-(n-2)\,\nu(\log y)=0.
\end{equation}

By Theorem~\ref{thm:A}, the profile curves are the geodesics of
\begin{equation}\label{eq:rot-metric}
 \widetilde g=y^{2c}g_{\Sph^2},
 \qquad
 c=\frac{n-2}{3}.
\end{equation}

\subsubsection{Polar coordinates around the pole and the Clairaut integral}

The metric \eqref{eq:rot-metric} is invariant under the rotations of $\Sph^2$ which fix the pole $N=(0,1,0)$; they are induced by the isometries of $\Sph^n$ which rotate the plane of the first and last coordinates and commute with $G$. It is therefore convenient to use geodesic polar coordinates centred at $N$,
\begin{equation}\label{eq:rot-coordinates}
 y=\cos\rho,\qquad
 x=\sin\rho\cos\psi,\qquad
 z=\sin\rho\sin\psi,
 \qquad 0<\rho<\frac\pi2,
\end{equation}
in which
\[
 g_{\Sph^2}=d\rho^2+\sin^2\rho\,d\psi^2,
 \qquad
 \widetilde g=\cos^{2c}\rho\,\bigl(d\rho^2+\sin^2\rho\,d\psi^2\bigr).
\]
The singular boundary $y=0$ is the circle $\rho=\pi/2$, and the pole $\rho=0$ is an interior point of $Q_{\Sph}$, over which lies the totally geodesic orbit $\Sph^{n-2}=\Sph^n\cap(\{0\}\times\R^{n-1}\times\{0\})$. Put $e_\rho=\partial_\rho$, $e_\psi=(\sin\rho)^{-1}\partial_\psi$, and write the unit tangent of a $g_{\Sph^2}$-unit-speed profile as
\[
 T=\cos\beta\,e_\rho+\sin\beta\,e_\psi .
\]
Since \eqref{eq:rot-reduced} does not depend on the orientation, we orient the hemisphere by the frame $(e_\rho,e_\psi)$; then, as in Section~\ref{subsec:sphere-system},
\[
 \nu=-\sin\beta\,e_\rho+\cos\beta\,e_\psi,
 \qquad
 \kappa=\beta'+\cot\rho\,\sin\beta .
\]
Since $\nu(\log y)=\nu(\log\cos\rho)=\tan\rho\,\sin\beta$, equation \eqref{eq:rot-reduced} is equivalent to the system
\begin{equation}\label{eq:rot-system}
 \rho'=\cos\beta,\qquad
 \psi'=\frac{\sin\beta}{\sin\rho},\qquad
 \beta'=\Bigl(c\tan\rho-\cot\rho\Bigr)\sin\beta .
\end{equation}
Moreover, along a solution,
\begin{equation}\label{eq:rot-kg}
 \kappa=c\tan\rho\,\sin\beta .
\end{equation}

Let
\begin{equation}\label{eq:rot-F}
 F(\rho)=\cos^{c}\rho\,\sin\rho,\qquad 0\le\rho\le\frac\pi2 .
\end{equation}
Then $F(0)=F(\pi/2)=0$, and $F$ has a unique critical point $\rho_0\in(0,\pi/2)$, a nondegenerate maximum, determined by, considering the value of $c$ in \eqref{eq:rot-metric},
\begin{equation}\label{eq:rot-rho0}
 \tan^2\rho_0=\frac1c=\frac{3}{n-2}.
\end{equation}

\begin{proposition}[Clairaut integral]\label{prop:clairaut}
Along every solution of \eqref{eq:rot-system} the function
\[
 E=F(\rho)\sin\beta=\cos^c\rho\,\sin\rho\,\sin\beta
\]
is constant. If $E\neq0$, then $\kappa$ never vanishes and
\begin{equation}\label{eq:rot-pelastic}
 y\,|\kappa|^{\frac{3}{n+1}}=(c|E|)^{\frac{3}{n+1}} .
\end{equation}
\end{proposition}

\begin{proof}
Using \eqref{eq:rot-system},
\[
 E'=F(\rho)\cos\beta\,\sin\beta\Bigl[\bigl(\cot\rho-c\tan\rho\bigr)+\bigl(c\tan\rho-\cot\rho\bigr)\Bigr]=0 .
\]
This is the Clairaut integral of the $\widetilde g$-geodesics: $E=\widetilde g(\widetilde T,\partial_\psi)$, where $\widetilde T=\cos^{-c}\rho\,T$ is the $\widetilde g$-unit tangent and $\partial_\psi$ is the Killing field of the rotations around $N$. By \eqref{eq:rot-kg},
\[
 \kappa=c\,\frac{\sin\rho\,\sin\beta}{\cos\rho}=\frac{cE}{\cos^{c+1}\rho}=\frac{cE}{y^{c+1}} .
\]
Hence $\kappa\neq0$ if $E\neq0$, and \eqref{eq:rot-pelastic} follows since $c+1=(n+1)/3$.
\end{proof}

\begin{remark}[$\sigma$-elastic profiles]\label{cor:pelastic}
Let $\sigma=(n-2)/(n+1)$. Since $3/(n+1)=1-\sigma$, relation \eqref{eq:rot-pelastic} shows that the profile curve of a proper biconservative $SO(n-1)$-invariant hypersurface of $\Sph^n$ with $E\neq0$ satisfies
\[
 y=a\,|\kappa|^{\sigma-1},
 \qquad a=(c|E|)^{3/(n+1)}>0 .
\]
{As a consequence, the profile is a critical curve of the curvature energy 
\[
\int_\gamma|\kappa|^{\sigma}\,ds,
\]
 i.e., a $\sigma$-elastic curve. To see this, recall that, for a smooth function $P$, the Euler--Lagrange equation of the functional 
 \[
 \int_\gamma P(\kappa)\,ds
 \]
  for curves in the unit sphere $\Sph^2$, parametrized by arc length, is
\begin{equation}\label{eq:EL-P}
\bigl(P'(\kappa)\bigr)''+\kappa^2P'(\kappa)-\kappa P(\kappa)+P'(\kappa)=0 .
\end{equation}
By Proposition~\ref{prop:clairaut}, $\kappa$ does not vanish and has the sign of $E$; with the orientation for which $E>0$ we have $\kappa>0$, so that $|\kappa|^\sigma=\kappa^\sigma$ and $y=a\kappa^{\sigma-1}$. Since $\gamma$ lies in the unit sphere $\Sph^2\subset\R^3$ and $y$ is the restriction of a linear function, $y''=\langle\gamma'',e_2\rangle=-y+\kappa\,\nu(y)$. By \eqref{eq:rot-reduced}, $\nu(y)=y\,\nu(\log y)=\kappa y/c$, so that
\[
y''=-y+\frac{\kappa^2}{c}\,y .
\]
Substituting $y=a\kappa^{\sigma-1}$, multiplying by $\sigma/a$ and using $\sigma/c=3/(n+1)=1-\sigma$, we obtain
\[
\sigma\bigl(\kappa^{\sigma-1}\bigr)''+(\sigma-1)\kappa^{\sigma+1}+\sigma\kappa^{\sigma-1}=0,
\]
which is \eqref{eq:EL-P} for $P(\kappa)=\kappa^\sigma$, since then $P'(\kappa)=\sigma\kappa^{\sigma-1}$ and $\kappa^2P'(\kappa)-\kappa P(\kappa)=(\sigma-1)\kappa^{\sigma+1}$.}
\end{remark}

In the present approach the Killing vector field used in \cite{MOP-JMAA} is simply the infinitesimal rotation $\partial_\psi$ around the pole, which is an isometry of the reduced metric $\widetilde g$. The solutions with $E=0$ are the great circles through $N$, whose lifts are totally geodesic hyperspheres.

\subsubsection{Periodicity and the closure condition}

\begin{proposition}\label{prop:rot-periodic}
{Let $\gamma$ be the profile curve of an $SO(n-1)$-invariant biconservative hypersurface of $\Sph^n$ which satisfies \eqref{eq:rot-reduced}, written in the polar coordinates \eqref{eq:rot-coordinates} around the pole $N$, and assume that its Clairaut constant $E$ is not zero; changing the orientation if necessary, we may take $E>0$. Then $0<E\le F(\rho_0)$, and the following hold.}
\begin{enumerate}[label=\textup{(\roman*)}]
\item If $E=F(\rho_0)$, the profile is the parallel $\rho=\rho_0$; its lift is a CMC hypersurface isometric to $\Sph^1(\sin\rho_0)\times\Sph^{n-2}(\cos\rho_0)$.
\item If $0<E<F(\rho_0)$, let $0<\rho_1<\rho_0<\rho_2<\pi/2$ be the two solutions of $F(\rho)=E$. Then the profile is defined for all $s\in\R$, $\rho(s)$ is a nonconstant periodic function oscillating between $\rho_1$ and $\rho_2$, and $\psi$ is strictly increasing. In particular the profile never meets the pole or the singular boundary, and its curvature $\kappa=cE/\cos^{c+1}\rho$ is periodic and nonconstant. During one period of $\rho$, the angle $\psi$ increases by
\begin{equation}\label{eq:rot-Phi}
 \Phi(E)=2\int_{\rho_1}^{\rho_2}\frac{E\,d\rho}{\sin\rho\,\sqrt{F(\rho)^2-E^2}} .
\end{equation}
\end{enumerate}
\end{proposition}

\begin{proof}
{Since $0<\sin\beta\le1$ and $F\le F(\rho_0)$, we have $0<E=F(\rho)\sin\beta\le F(\rho_0)$. }

(i) If $E=F(\rho_0)$ then $F(\rho)\sin\beta=\max F$ forces $\rho\equiv\rho_0$ and $\beta\equiv\pi/2$, which is a solution of \eqref{eq:rot-system} by \eqref{eq:rot-rho0}. Its curvature $\kappa=c\tan\rho_0$ is constant, so the lift is CMC.

{(ii) Let $0<E<F(\rho_0)$. We proceed in four steps.

\emph{The phase plane.} The pair $(\rho,\beta)$ solves the planar system
\begin{equation}\label{eq:rot-planar}
\rho'=\cos\beta,
\qquad
\beta'=\bigl(c\tan\rho-\cot\rho\bigr)\sin\beta,
\end{equation}
which does not involve $\psi$; once $(\rho,\beta)$ is known, $\psi$ is obtained by integrating $\psi'=\sin\beta/\sin\rho$. Since $F>0$ on $(0,\pi/2)$ and $E=F(\rho)\sin\beta>0$, we have $\sin\beta>0$, so that, $\beta$ being defined modulo $2\pi$, the solution lies in the strip $\mathcal S=(0,\pi/2)\times(0,\pi)$. In $\mathcal S$ an equilibrium of \eqref{eq:rot-planar} has $\cos\beta=0$, i.e., $\beta=\pi/2$, and then $c\tan\rho=\cot\rho$, i.e., $\rho=\rho_0$ by \eqref{eq:rot-rho0}. Thus $(\rho_0,\pi/2)$ is the only equilibrium in $\mathcal S$.

\emph{The level set.} By Proposition~\ref{prop:clairaut} the solution lies on the level set
\[
L_E=\bigl\{(\rho,\beta)\in\mathcal S:\ F(\rho)\sin\beta=E\bigr\}.
\]
Since $F$ is strictly increasing on $(0,\rho_0]$, strictly decreasing on $[\rho_0,\pi/2)$, and vanishes at $0$ and $\pi/2$, the inequality $F(\rho)\ge E$ holds exactly for $\rho\in[\rho_1,\rho_2]$, and $F'(\rho_1)>0>F'(\rho_2)$. For $\rho\in(\rho_1,\rho_2)$ the equation $\sin\beta=E/F(\rho)<1$ has exactly two solutions in $(0,\pi)$, namely $\beta_E(\rho)=\arcsin\bigl(E/F(\rho)\bigr)\in(0,\pi/2)$ and $\pi-\beta_E(\rho)$, while for $\rho=\rho_1,\rho_2$ it has the single solution $\beta=\pi/2$. Hence $L_E$ is the union of the graphs of $\beta_E$ and $\pi-\beta_E$ over $[\rho_1,\rho_2]$, joined at the points $(\rho_1,\pi/2)$ and $(\rho_2,\pi/2)$: it is a simple closed curve, and it is compact, being contained in $[\rho_1,\rho_2]\times[\beta_m,\pi-\beta_m]\subset\mathcal S$ with $\beta_m=\arcsin\bigl(E/F(\rho_0)\bigr)>0$. Moreover $L_E$ does not contain the equilibrium, because $F(\rho_0)\sin(\pi/2)=F(\rho_0)>E$.

\emph{Periodicity.} Since $L_E$ is a compact subset of $\mathcal S$, the solution of \eqref{eq:rot-planar} cannot leave every compact subset of $\mathcal S$, hence it is defined for all $s\in\R$. The vector field \eqref{eq:rot-planar} does not vanish on $L_E$, so its length has a positive minimum on the compact set $L_E$; therefore the solution travels along the closed curve $L_E$ with speed bounded from below, goes around it in a finite time $S_E$, and is periodic with period $S_E$, its trajectory being the whole of $L_E$. Consequently $\rho(s)$ is periodic, takes all the values in $[\rho_1,\rho_2]$, and is nonconstant because $\rho_1<\rho_2$. More precisely, $\rho'=\cos\beta>0$ on the branch $\beta=\beta_E(\rho)<\pi/2$ and $\rho'<0$ on the branch $\beta=\pi-\beta_E(\rho)>\pi/2$: during one period, $\rho$ increases from $\rho_1$ to $\rho_2$ and then decreases back to $\rho_1$. Since $\rho_1\le\rho\le\rho_2$ with $0<\rho_1$ and $\rho_2<\pi/2$, the profile never meets the pole $\rho=0$ nor the singular boundary $\rho=\pi/2$; moreover $\psi'=\sin\beta/\sin\rho>0$, so that $\psi$ is strictly increasing. By Proposition~\ref{prop:clairaut}, $\kappa=cE/\cos^{c+1}\rho$, which is periodic and, since $\cos$ is injective on $(0,\pi/2)$ and $\rho$ is not constant, nonconstant.

\emph{The rotation angle.} On each of the two halves of a period, $\rho$ is strictly monotone and can be used as a parameter. There $\sin\beta=E/F$ and $\cos\beta=\pm\sqrt{1-E^2/F^2}=\pm\sqrt{F^2-E^2}/F$, with the sign $+$ on the half where $\rho$ increases and $-$ on the other, so that
\[
\frac{d\psi}{d\rho}=\frac{\psi'}{\rho'}=\frac{\sin\beta}{\sin\rho\,\cos\beta}=\pm\frac{E}{\sin\rho\,\sqrt{F(\rho)^2-E^2}} .
\]
On the first half $\rho$ goes from $\rho_1$ to $\rho_2$, and on the second from $\rho_2$ back to $\rho_1$; in both cases $\psi$ increases by $\int_{\rho_1}^{\rho_2}E\,d\rho/\bigl(\sin\rho\sqrt{F^2-E^2}\bigr)$. This integral is finite: since $F'(\rho_1)\neq0$ and $F'(\rho_2)\neq0$, we have $F^2-E^2\ge C\,(\rho-\rho_1)(\rho_2-\rho)$ on $[\rho_1,\rho_2]$ for some $C>0$, so that the integrand has only integrable singularities of the type $|\rho-\rho_i|^{-1/2}$ at the endpoints. Adding the two contributions gives \eqref{eq:rot-Phi}.}
\end{proof}

\begin{figure}[!htbp]
\centering
\begin{tikzpicture}[line cap=round,
  guide/.style={dashed, line width=.5pt, gray},
  thetaline/.style={densely dotted, line width=.5pt, gray},
  vec/.style={line width=.8pt, -{Latex[length=2mm]}},
  lab/.style={font=\scriptsize, fill=white, inner sep=1pt}]
\draw[line width=1.1pt] (0,0) circle (3.2);
\draw[thetaline] (0,0) circle (0.4576);
\draw[guide] (0,0) circle (1.2538);
\draw[thetaline] (0,0) circle (2.149);
\node[lab] at (200:2.149) {$\rho_2$};
\node[lab] at (200:1.2538) {$\rho_0$};
\node[lab, anchor=east] at (180:0.4776) {$\rho_1$};
\node[font=\scriptsize, anchor=north west] at (-40:3.25) {$\rho=\frac\pi2$};
\draw[line width=.6pt, gray!70] plot coordinates {(-0.396,-0.229) (-0.382,-0.253) (-0.368,-0.277) (-0.353,-0.301) (-0.338,-0.325) (-0.324,-0.349) (-0.309,-0.372) (-0.293,-0.395) (-0.278,-0.419) (-0.262,-0.442) (-0.246,-0.465) (-0.230,-0.488) (-0.214,-0.510) (-0.198,-0.533) (-0.181,-0.555) (-0.164,-0.578) (-0.147,-0.600) (-0.130,-0.622) (-0.113,-0.644) (-0.095,-0.665) (-0.077,-0.687) (-0.059,-0.708) (-0.041,-0.730) (-0.023,-0.751) (-0.004,-0.772) (0.015,-0.792) (0.034,-0.813) (0.053,-0.833) (0.072,-0.853) (0.092,-0.873) (0.112,-0.893) (0.132,-0.913) (0.152,-0.932) (0.172,-0.951) (0.193,-0.970) (0.214,-0.989) (0.235,-1.008) (0.256,-1.026) (0.278,-1.044) (0.300,-1.062) (0.322,-1.079) (0.344,-1.096) (0.366,-1.113) (0.389,-1.130) (0.412,-1.147) (0.435,-1.163) (0.458,-1.178) (0.482,-1.194) (0.506,-1.209) (0.530,-1.224) (0.554,-1.238) (0.579,-1.252) (0.604,-1.266) (0.629,-1.279) (0.654,-1.292) (0.680,-1.305) (0.706,-1.317) (0.732,-1.328) (0.758,-1.339) (0.785,-1.350)};
\draw[line width=1.1pt, postaction={decorate}, decoration={markings, mark=at position 0.5 with {\arrow{Latex[length=2mm]}}}] plot coordinates {(0.396,-0.229) (0.411,-0.203) (0.426,-0.176) (0.441,-0.149) (0.456,-0.122) (0.470,-0.095) (0.484,-0.068) (0.498,-0.040) (0.511,-0.013) (0.524,0.015) (0.537,0.043) (0.550,0.070) (0.563,0.098) (0.575,0.127) (0.587,0.155) (0.599,0.183) (0.610,0.212) (0.621,0.240) (0.632,0.269) (0.643,0.298) (0.653,0.327) (0.663,0.356) (0.673,0.385) (0.682,0.414) (0.691,0.443) (0.700,0.473) (0.709,0.502) (0.717,0.532) (0.725,0.562) (0.732,0.592) (0.739,0.622) (0.746,0.652) (0.752,0.682) (0.758,0.712) (0.764,0.743) (0.769,0.773) (0.774,0.803) (0.779,0.834) (0.783,0.865) (0.786,0.895) (0.789,0.926) (0.792,0.957) (0.794,0.988) (0.796,1.019) (0.797,1.050) (0.798,1.081) (0.798,1.112) (0.798,1.144) (0.797,1.175) (0.796,1.206) (0.794,1.237) (0.791,1.268) (0.788,1.300) (0.784,1.331) (0.779,1.362) (0.773,1.393) (0.767,1.424) (0.760,1.455) (0.753,1.486) (0.744,1.517) (0.735,1.547) (0.725,1.578) (0.714,1.608) (0.702,1.638) (0.689,1.668) (0.675,1.697) (0.660,1.726) (0.644,1.755) (0.627,1.783) (0.608,1.811) (0.589,1.838) (0.569,1.865) (0.547,1.891) (0.524,1.916) (0.500,1.940) (0.475,1.964) (0.449,1.986) (0.421,2.008) (0.393,2.028) (0.363,2.047) (0.332,2.064) (0.300,2.081) (0.267,2.095) (0.234,2.108) (0.199,2.120) (0.164,2.129) (0.128,2.137) (0.092,2.143) (0.055,2.147) (0.018,2.149) (-0.018,2.149) (-0.055,2.147) (-0.092,2.143) (-0.128,2.137) (-0.164,2.129) (-0.199,2.120) (-0.234,2.108) (-0.267,2.095) (-0.300,2.081) (-0.332,2.064) (-0.363,2.047) (-0.393,2.028) (-0.421,2.008) (-0.449,1.986) (-0.475,1.964) (-0.500,1.940) (-0.524,1.916) (-0.547,1.891) (-0.569,1.865) (-0.589,1.838) (-0.608,1.811) (-0.627,1.783) (-0.644,1.755) (-0.660,1.726) (-0.675,1.697) (-0.689,1.668) (-0.702,1.638) (-0.714,1.608) (-0.725,1.578) (-0.735,1.547) (-0.744,1.517) (-0.753,1.486) (-0.760,1.455) (-0.767,1.424) (-0.773,1.393) (-0.779,1.362) (-0.784,1.331) (-0.788,1.300) (-0.791,1.268) (-0.794,1.237) (-0.796,1.206) (-0.797,1.175) (-0.798,1.144) (-0.798,1.112) (-0.798,1.081) (-0.797,1.050) (-0.796,1.019) (-0.794,0.988) (-0.792,0.957) (-0.789,0.926) (-0.786,0.895) (-0.783,0.865) (-0.779,0.834) (-0.774,0.803) (-0.769,0.773) (-0.764,0.743) (-0.758,0.712) (-0.752,0.682) (-0.746,0.652) (-0.739,0.622) (-0.732,0.592) (-0.725,0.562) (-0.717,0.532) (-0.709,0.502) (-0.700,0.473) (-0.691,0.443) (-0.682,0.414) (-0.673,0.385) (-0.663,0.356) (-0.653,0.327) (-0.643,0.298) (-0.632,0.269) (-0.621,0.240) (-0.610,0.212) (-0.599,0.183) (-0.587,0.155) (-0.575,0.127) (-0.563,0.098) (-0.550,0.070) (-0.537,0.043) (-0.524,0.015) (-0.511,-0.013) (-0.498,-0.040) (-0.484,-0.068) (-0.470,-0.095) (-0.456,-0.122) (-0.441,-0.149) (-0.426,-0.176) (-0.411,-0.203) (-0.396,-0.229)};
\draw[line width=.5pt] (0,0) -- (-30.08:0.4576);
\draw[line width=.5pt] (0,0) -- (210.08:0.4576);
\draw[line width=.5pt, -{Latex[length=1.4mm]}] (-30.08:0.28) arc[start angle=-30.08, end angle=210.08, radius=0.28];
\node[font=\scriptsize] at (90:0.62) {$\Phi$};
\fill (-30.08:0.4576) circle (1.3pt);
\fill (210.08:0.4576) circle (1.3pt);
\fill (0,0) circle (1.4pt);
\node[font=\scriptsize, anchor=north] at (0,-0.06) {$N$};
\coordinate (X) at (0.7933,1.2409);
\draw[guide] (X) -- ++(57.41:1.0) node[font=\scriptsize, anchor=south west, inner sep=1pt] {$e_\rho$};
\draw[vec] (X) -- ++(94.51:1.0) node[font=\small, anchor=south] {$T$};
\draw[line width=.5pt] ($(X)+(57.41:0.55)$) arc[start angle=57.41, end angle=94.51, radius=0.55];
\node[font=\scriptsize] at ($(X)+(76:0.75)$) {$\beta$};
\fill (X) circle (1.3pt);
\node[font=\small, anchor=east] at (-0.62,1.75) {$\gamma$};
\end{tikzpicture}
\caption{The rotational case $d=1$ in $\Sph^n$ (here $n=8$ and $E=0.55\,F(\rho_0)$), drawn in geodesic polar coordinates $(\rho,\psi)$ around the pole $N$; the thick circle $\rho=\pi/2$ is the singular boundary $y=0$. By the Clairaut integral $E=F(\rho)\sin\beta$, where $\beta$ is the angle between the unit vector $e_\rho=\partial_\rho$ along the meridian through the point (dashed) and the unit tangent $T$, the profile oscillates between the parallels $\rho_1$ and $\rho_2$ (dotted) defined by $F(\rho)=E$, around the CMC parallel $\rho_0$ (dashed). The thick arc is one period of $\rho$, between two consecutive points of minimal distance from $N$ (dots); during this period the profile turns around $N$ by the angle $\Phi=\Phi(E)$ of \eqref{eq:rot-Phi}, and the next period, of which only the first part is drawn (grey), is obtained by rotating it by $\Phi$ about $N$. Hence the profile closes up exactly when $r\Phi=2\pi l$.}
\label{fig:rot-Phi}
\end{figure}
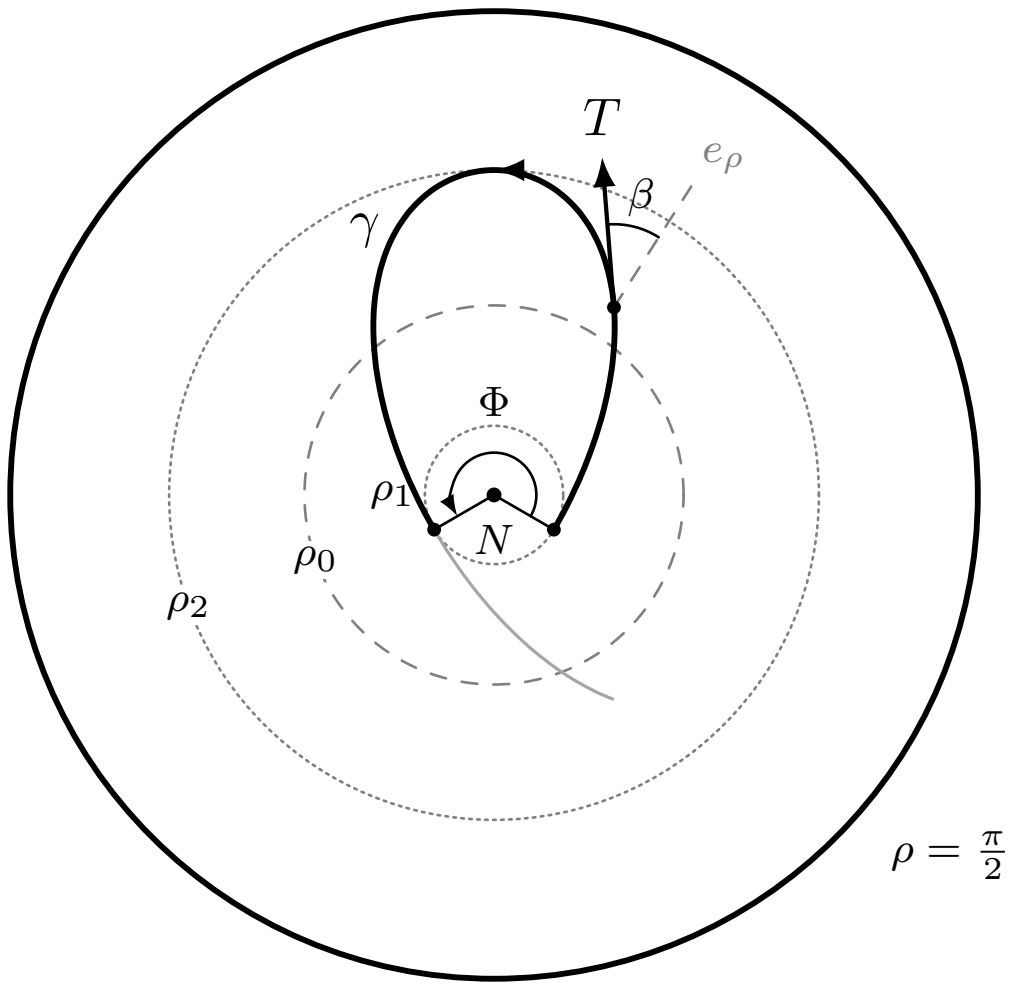

Since $\rho$ is periodic and $\psi$ increases by $\Phi(E)$ along each period, the profile is closed if and only if $r\Phi(E)=2\pi l$ for some natural numbers $l,r$, which we may assume coprime. Then $l$ is the number of turns around the pole and $r$ is the number of periods of the curvature (lobes) contained in the closed curve. When this happens, the lift is a compact non-CMC biconservative immersion of $\Sph^1\times\Sph^{n-2}$ into $\Sph^n$. The geometric meaning of $\Phi$ is illustrated in Figure~\ref{fig:rot-Phi}.

\begin{lemma}\label{lem:Phi-limits}
The function $\Phi$ is continuous on $(0,F(\rho_0))$, and
\[
 \lim_{E\to F(\rho_0)}\Phi(E)=\sqrt2\,\pi,
 \qquad
 \lim_{E\to0}\Phi(E)=\pi,
\]
for every $n\geq3$.
\end{lemma}

\begin{proof}
Continuity follows by writing \eqref{eq:rot-Phi} on a fixed interval through the substitution $\rho=\rho_1+(\rho_2-\rho_1)(1-\cos\tau)/2$, $\tau\in[0,\pi]$, which removes the square-root singularities at the endpoints.

\emph{The limit $E\to F(\rho_0)$.} Put $G=F^2$ and $e=E^2$. Since $G$ has a nondegenerate maximum at $\rho_0$, the Morse lemma gives a local coordinate $u$ with $G=G(\rho_0)-u^2$, $u(\rho_0)=0$ and $d\rho/du(0)=\sqrt{2/(-G''(\rho_0))}$. Then
\[
 \Phi(E)=2E\int_{-\sqrt{G(\rho_0)-e}}^{\sqrt{G(\rho_0)-e}}\frac{1}{\sin\rho}\frac{d\rho}{du}\frac{du}{\sqrt{G(\rho_0)-e-u^2}}
 \longrightarrow
 \frac{2\pi F(\rho_0)}{\sin\rho_0}\sqrt{\frac{2}{-G''(\rho_0)}} .
\]
Since $F'(\rho_0)=0$, $G''(\rho_0)=2F(\rho_0)F''(\rho_0)$, and differentiating $\log F=c\log\cos\rho+\log\sin\rho$ twice gives, by \eqref{eq:rot-rho0},
\[
 -\frac{F''(\rho_0)}{F(\rho_0)}=\frac{c}{\cos^2\rho_0}+\frac{1}{\sin^2\rho_0}=c(1+\tan^2\rho_0)+1+\cot^2\rho_0=2(1+c),
 \qquad
 \sin^2\rho_0=\frac{1}{1+c}.
\]
Hence the limit equals $2\pi/\bigl(\sin\rho_0\sqrt{2(1+c)}\bigr)=\sqrt2\,\pi$.

\emph{The limit $E\to0$.} Now $\rho_1\to0$ and $\rho_2\to\pi/2$. Fix ${\rho_*}\in(0,\rho_0)$ so small that $F$ is increasing on $[0,{\rho_*}]$ and decreasing on $[\pi/2-{\rho_*},\pi/2]$, put ${u_*}=F({\rho_*})$, ${\bar u_*}=F(\pi/2-{\rho_*})$, and let $E<\min({u_*},{\bar u_*})$. Split the integral in \eqref{eq:rot-Phi} into the parts over $[\rho_1,{\rho_*}]$, $[{\rho_*},\pi/2-{\rho_*}]$ and $[\pi/2-{\rho_*},\rho_2]$. On the middle interval $F$ is bounded below by a positive constant, so that part is $O(E)$.

On each of the two outer intervals we use $u=F(\rho)$ as a new variable; this is legitimate since $F$ is monotone there, and it turns the singular factor into exactly $\sqrt{u^2-E^2}$, with lower endpoint $u=E$. Near the pole, $F(\rho)=\rho\,(1+O(\rho^2))$, $F'(\rho)=1+O(\rho^2)$ and $\sin\rho=\rho\,(1+O(\rho^2))$, so that
\[
 \frac{d\rho}{\sin\rho}=\frac{du}{F'(\rho)\sin\rho}=\bigl(1+\varepsilon(u)\bigr)\frac{du}{u},
 \qquad |\varepsilon(u)|\le C u^2 ,
\]
with $C$ independent of $E$. Hence the first part equals
\[
 2\int_E^{{u_*}}\bigl(1+\varepsilon(u)\bigr)\frac{E\,du}{u\sqrt{u^2-E^2}}
 =2\arccos\frac{E}{{u_*}}+R,
 \qquad |R|\le C{u_*^2}\cdot2\int_E^{\infty}\frac{E\,du}{u\sqrt{u^2-E^2}}=\pi C {u_*^2} .
\]
Near the boundary, with $\tau=\pi/2-\rho$, we have $F=\tau^c(1+O(\tau^2))$, $|F'|=c\,\tau^{c-1}(1+O(\tau^2))$ and $\sin\rho=1+O(\tau^2)$; hence $d\rho/\sin\rho=(1+O(u^{2/c}))\,du/(c\,u^{1-1/c})$, and the last part is bounded by a constant times
\[
 E\int_E^{{\bar u_*}}\frac{u^{1/c-1}\,du}{\sqrt{u^2-E^2}}
 =E^{1/c}\int_1^{{\bar u_*}/E}\frac{v^{1/c-1}\,dv}{\sqrt{v^2-1}} ,
\]
which tends to $0$ as $E\to0$: it is $O(E^{1/c})$ if $c>1$, $O(E\log(1/E))$ if $c=1$, and $O(E)$ if $c<1$. Therefore $\limsup_{E\to0}|\Phi(E)-\pi|\le\pi C {u_*^2}$, and letting ${\rho_*}\to0$ gives $\lim_{E\to0}\Phi(E)=\pi$.
\end{proof}

We can now recover the main results of \cite{MP-CAG,MOP-JMAA}.

\begin{theorem}\label{thm:rot-closed}
Let $n\geq3$, and let $l,r$ be coprime natural numbers with
\[
 r<2l<\sqrt2\,r .
\]
Then there is a closed profile curve with $l$ turns around the pole and $r$ lobes, and its lift is a compact non-CMC biconservative $SO(n-1)$-invariant immersion $\Sph^1\times\Sph^{n-2}\to\Sph^n$. More generally, a closed profile with $\pi<\Phi(E)<2\pi$ is never simple; in particular, none of the hypersurfaces above is embedded.
\end{theorem}

\begin{figure}[!htbp]
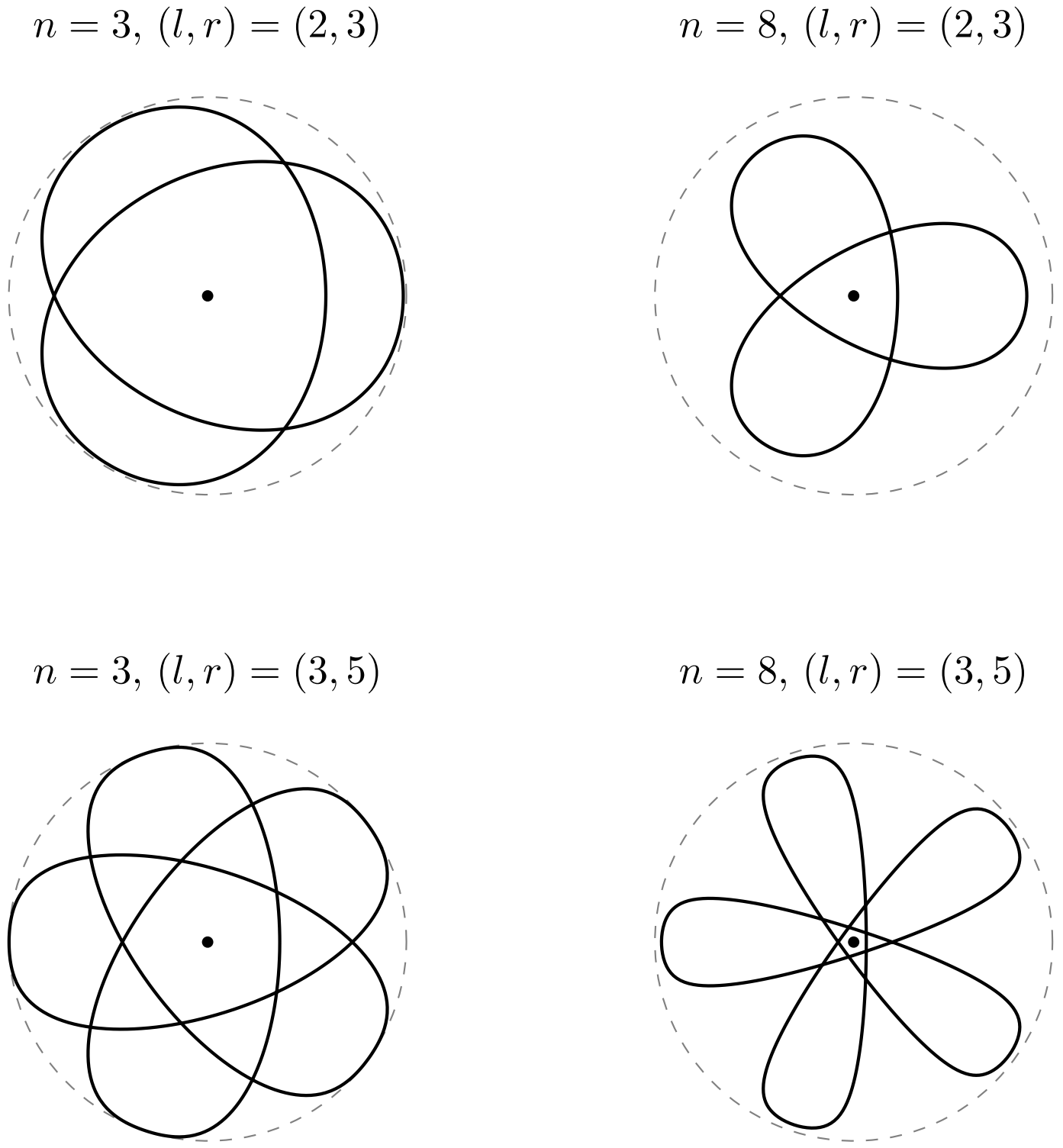

\centering

\caption{Closed profiles of Theorem~\ref{thm:rot-closed}, seen from above the pole $N$ (orthogonal projection of the hemisphere $Q_{\Sph}$ onto the $(x,z)$-plane). The dashed circle is the singular boundary $y=0$ and the dot is the pole; the profiles stay in the annulus $\rho_1\le\rho\le\rho_2$ (Proposition~\ref{prop:rot-periodic}) and never reach the boundary, although for $n=3$ they come close to it. Top row: $(l,r)=(2,3)$; bottom row: $(l,r)=(3,5)$. Left: $n=3$, the surfaces of \cite{MP-CAG}; right: $n=8$. The curves were obtained by solving $\Phi(E)=2\pi l/r$ numerically and integrating \eqref{eq:rot-system}.}
\label{fig:rot-closed}
\end{figure}

\begin{proof}
The condition on $l,r$ means $\pi<2\pi l/r<\sqrt2\pi$, so by Lemma~\ref{lem:Phi-limits} and the intermediate value theorem there is $E\in(0,F(\rho_0))$ with $\Phi(E)=2\pi l/r$, and the profile closes up. The lift is non-CMC because $\kappa$ is not constant (Proposition~\ref{prop:rot-periodic}), and $f=-2\kappa$ by Theorem~\ref{thm:A}. The closed profiles obtained for $(l,r)=(2,3)$ and $(3,5)$, with $n=3$ and $n=8$, are shown in Figure~\ref{fig:rot-closed}, where the self-intersections discussed below are also visible.

Let now $\gamma$ be a closed profile with $E\neq0$ and $\pi<\Phi(E)<2\pi$. By Proposition~\ref{prop:rot-periodic} it lies in the regular part of the orbit space, so its lift is embedded if and only if $\gamma$ is simple. Since $\psi$ is strictly increasing, we can write $\rho=R(\psi)$, where $R$ is periodic with minimal period $\Phi=\Phi(E)$, and the closed curve corresponds to $\psi\in[0,2\pi l)$. The continuous function $R(\psi+2\pi)-R(\psi)$ has zero mean over a period. It does not vanish identically, because $2\pi/\Phi\in(1,2)$ is not an integer; hence it vanishes at some $\psi_*$, and $\gamma(\psi_*)=\gamma(\psi_*+2\pi)$. If $l\geq2$, these are two different points of the closed curve, so $\gamma$ is not simple. If $l=1$, then $\Phi=2\pi/r$ with $r\in\mathbb N$; since $\Phi<2\pi$ we get $r\ge2$ and $\Phi\le\pi$, a contradiction.
\end{proof}

\subsubsection{Monotonicity of the period function}\label{subsec:Phi-monotone}

The monotonicity of $\Phi$ follows from a theorem of Otsuki. In \cite{Otsuki93} (see also \cite{Otsuki76,Otsuki79}), Otsuki studied, for every \emph{real} number $\mathsf n>1$, the geodesics of the metric
\[
 O^2_{\mathsf n}=(1-x^2-z^2)^{\mathsf n-2}\bigl\{(1-z^2)\,dx^2+2xz\,dx\,dz+(1-x^2)\,dz^2\bigr\}
\]
on the unit disk $x^2+z^2<1$. For integer $\mathsf n$ these geodesics are the profile curves of the rotational minimal hypersurfaces of $\Sph^{\mathsf n+1}${, since $O^2_{\mathsf n}=y^{2(\mathsf n-1)}g_{\Sph^2}$ is the Hsiang--Lawson metric $V^2g_{\Sph^2}$ of the action of $SO(\mathsf n)$ on $\Sph^{\mathsf n+1}$}, but Otsuki observed that the problem has a geometric meaning for every real $\mathsf n$. We identify the disk with our hemisphere $\{y>0\}$ by the orthogonal projection $(x,y,z)\mapsto(x,z)$, so that $y^2=1-x^2-z^2$ and the pole $N$ corresponds to the origin. In these coordinates
\[
 g_{\Sph^2}=dx^2+dz^2+\frac{(x\,dx+z\,dz)^2}{1-x^2-z^2}=\frac{(1-z^2)\,dx^2+2xz\,dx\,dz+(1-x^2)\,dz^2}{1-x^2-z^2},
\]
so that $O^2_{\mathsf n}=y^{2(\mathsf n-1)}g_{\Sph^2}$. Comparing with \eqref{eq:rot-metric}, \emph{the profile curves of the proper biconservative $SO(n-1)$-invariant hypersurfaces of $\Sph^n$ are exactly the geodesics of $O^2_{\mathsf n}$ with nonconstant curvature, where}
\[
 \mathsf n=c+1=\frac{n+1}{3}>1 .
\]

Otsuki showed that{, along} the non-trivial geodesics of $O^2_{\mathsf n}$ which do not pass through the origin, {the distance from the origin is a periodic function of the polar angle (compare Proposition~\ref{prop:rot-periodic}), and that its period, i.e., the angle swept around the origin while this distance goes from one minimum to the next,} is given by \cite[(2.9)]{Otsuki93}
\begin{equation}\label{eq:Otsuki-T}
 T({\mathsf c})=\sqrt{\mathsf n\,{\mathsf c}}\int_{x_0}^{x_1}\frac{dx}{x\sqrt{(\mathsf n-x)\bigl(x(\mathsf n-x)^{\mathsf n-1}-{\mathsf c}\bigr)}},
 \qquad 0<{\mathsf c}<B=(\mathsf n-1)^{\mathsf n-1},
\end{equation}
where $0<x_0<1<x_1<\mathsf n$ are the two solutions of $x(\mathsf n-x)^{\mathsf n-1}={\mathsf c}$. {Otsuki denotes this constant by $c$; as for $\mathsf n$, we use a sans-serif letter for it, to avoid confusion with the exponent $c$ of \eqref{eq:rot-metric}. The geodesic itself is closed if and only if $T(\mathsf c)/2\pi$ is rational, so that, in general, it is not periodic.} Here $x=\mathsf n\,r^2$, where $r$ is the Euclidean distance from the origin in the disk \cite[(2.8)]{Otsuki93}; in our coordinates $r=\sin\rho$. Accordingly, put
\[
 x=\mathsf n\sin^2\!\rho,\qquad \mathsf n=c+1 .
\]
Then $\mathsf n-x=\mathsf n\cos^2\!\rho$ and $x(\mathsf n-x)^{\mathsf n-1}=\mathsf n^{\mathsf n}\sin^2\!\rho\,\cos^{2c}\!\rho=\mathsf n^{\mathsf n}F(\rho)^2$, so that $x_0$ and $x_1$ correspond to the turning parallels $\rho_1$ and $\rho_2$, the critical point $x=1$ to the CMC parallel $\rho_0$ (since $\sin^2\rho_0=1/(c+1)$ by \eqref{eq:rot-rho0}), and the maximum $B$ of $x(\mathsf n-x)^{\mathsf n-1}$ to $\mathsf n^{\mathsf n}F(\rho_0)^2$. Moreover $dx=2\mathsf n\sin\rho\cos\rho\,d\rho$, and a direct substitution in \eqref{eq:Otsuki-T} with
\begin{equation}\label{eq:Otsuki-c}
 {\mathsf c}=\mathsf n^{\mathsf n}E^2
\end{equation}
gives
\[
 \sqrt{\mathsf n\,{\mathsf c}}\,\frac{dx}{x\sqrt{(\mathsf n-x)\bigl(x(\mathsf n-x)^{\mathsf n-1}-{\mathsf c}\bigr)}}
 =\frac{2E\,d\rho}{\sin\rho\sqrt{F(\rho)^2-E^2}} .
\]
Comparing with \eqref{eq:rot-Phi}, we obtain the exact identity
\begin{equation}\label{eq:Phi-T}
 \Phi(E)=T\bigl(\mathsf n^{\mathsf n}E^2\bigr),\qquad 0<E<F(\rho_0),\qquad \mathsf n=\frac{n+1}{3},
\end{equation}
in which $E\mapsto\mathsf n^{\mathsf n}E^2$ is an increasing bijection of $(0,F(\rho_0))$ onto $(0,B)$. Otsuki's Theorem states that $T=T({\mathsf c})$ is monotone increasing on $0<{\mathsf c}<(\mathsf n-1)^{\mathsf n-1}$ for every real $\mathsf n$ with $1<\mathsf n<\infty$; the proof is given in \cite[\S5]{Otsuki93}{, where the strict inequality $T'(\mathsf c)>0$ is established \cite[(5.3$'$)]{Otsuki93}} (see also \cite{Otsuki76,Otsuki79}). Since $\mathsf n=(n+1)/3>1$ for every $n\ge3$, we obtain the following result.

\begin{theorem}[Otsuki {\cite[Theorem~B]{Otsuki93}}]\label{thm:Phi-monotone}
Let $n\ge3$. The function $\Phi$ is strictly increasing on $(0,F(\rho_0))$; hence it is an increasing homeomorphism of $(0,F(\rho_0))$ onto $(\pi,\sqrt2\,\pi)$.
\end{theorem}

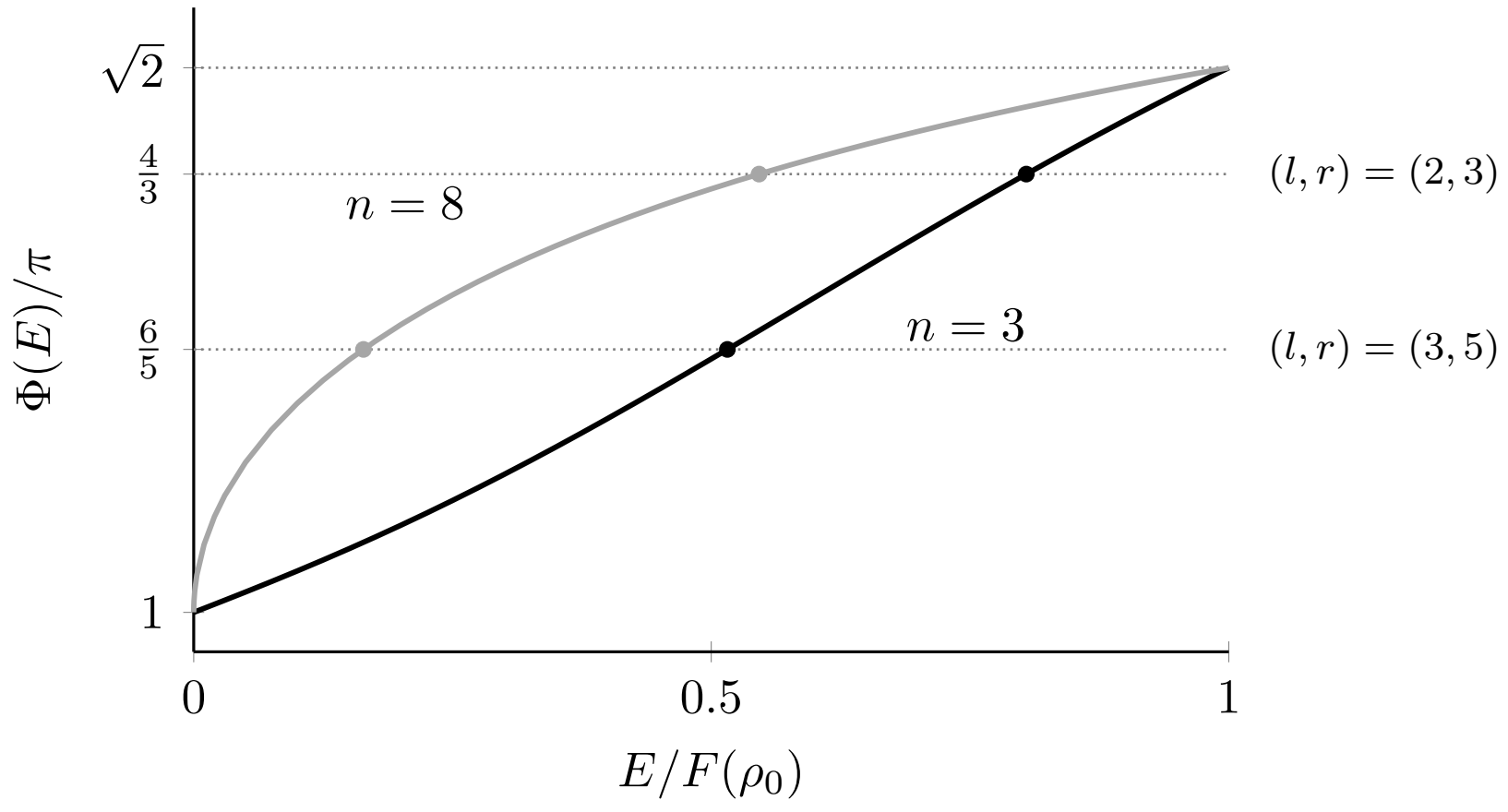
\begin{figure}[!htbp]
\centering
\begin{tikzpicture}
\begin{axis}[width=9cm, height=6.2cm, xmin=0, xmax=1, ymin=0.97, ymax=1.46,
  axis lines=left, axis line style={line width=.6pt, -}, clip=false,
  xlabel={$E/F(\rho_0)$}, ylabel={$\Phi(E)/\pi$}, label style={font=\small},
  xtick={0,0.5,1}, ytick={1,1.2,1.33333,1.41421}, yticklabels={$1$,$\frac65$,$\frac43$,$\sqrt2$},
  tick label style={font=\small}]
  \addplot[densely dotted, gray, line width=.5pt] coordinates {(0,1.41421) (1,1.41421)};
  \addplot[densely dotted, gray, line width=.5pt] coordinates {(0,1.33333) (1,1.33333)};
  \addplot[densely dotted, gray, line width=.5pt] coordinates {(0,1.2) (1,1.2)};
  \addplot[line width=1.1pt, black] coordinates {(0.0000,1.00000) (0.0001,1.00003) (0.0003,1.00009) (0.0010,1.00031) (0.0030,1.00092) (0.0100,1.00307) (0.0200,1.00615) (0.0300,1.00924) (0.0500,1.01547) (0.0750,1.02339) (0.1000,1.03150) (0.1250,1.03982) (0.1500,1.04837) (0.1750,1.05717) (0.2000,1.06621) (0.2250,1.07552) (0.2500,1.08508) (0.2750,1.09490) (0.3000,1.10497) (0.3250,1.11527) (0.3500,1.12580) (0.3750,1.13654) (0.4000,1.14747) (0.4250,1.15858) (0.4500,1.16984) (0.4750,1.18124) (0.5000,1.19275) (0.5250,1.20435) (0.5500,1.21602) (0.5750,1.22772) (0.6000,1.23945) (0.6250,1.25117) (0.6500,1.26286) (0.6750,1.27451) (0.7000,1.28609) (0.7250,1.29759) (0.7500,1.30898) (0.7750,1.32024) (0.8000,1.33137) (0.8250,1.34236) (0.8500,1.35318) (0.8750,1.36382) (0.9000,1.37429) (0.9250,1.38457) (0.9500,1.39465) (0.9700,1.40257) (0.9900,1.41037) (0.9970,1.41306) (1.0000,1.41421)};
  \addplot[line width=1.1pt, gray!70] coordinates {(0.0000,1.00000) (0.0001,1.00518) (0.0003,1.00897) (0.0010,1.01637) (0.0030,1.02834) (0.0100,1.05162) (0.0200,1.07278) (0.0300,1.08887) (0.0500,1.11405) (0.0750,1.13867) (0.1000,1.15898) (0.1250,1.17650) (0.1500,1.19202) (0.1750,1.20601) (0.2000,1.21878) (0.2250,1.23055) (0.2500,1.24147) (0.2750,1.25167) (0.3000,1.26124) (0.3250,1.27026) (0.3500,1.27879) (0.3750,1.28688) (0.4000,1.29458) (0.4250,1.30192) (0.4500,1.30894) (0.4750,1.31565) (0.5000,1.32209) (0.5250,1.32828) (0.5500,1.33423) (0.5750,1.33996) (0.6000,1.34549) (0.6250,1.35083) (0.6500,1.35599) (0.6750,1.36098) (0.7000,1.36582) (0.7250,1.37051) (0.7500,1.37505) (0.7750,1.37947) (0.8000,1.38375) (0.8250,1.38792) (0.8500,1.39198) (0.8750,1.39593) (0.9000,1.39977) (0.9250,1.40352) (0.9500,1.40717) (0.9700,1.41003) (0.9900,1.41283) (0.9970,1.41380) (1.0000,1.41421)};
  \addplot[only marks, mark=*, mark size=1.5pt, black] coordinates {(0.8044,1.33333) (0.5156,1.20000)};
  \addplot[only marks, mark=*, mark size=1.5pt, gray!70] coordinates {(0.5462,1.33333) (0.1640,1.20000)};
  \node[font=\small, anchor=north west] at (axis cs:0.67,1.245) {$n=3$};
  \node[font=\small, anchor=south east] at (axis cs:0.28,1.285) {$n=8$};
  \node[font=\scriptsize, anchor=west] at (axis cs:1.02,1.33333) {$(l,r)=(2,3)$};
  \node[font=\scriptsize, anchor=west] at (axis cs:1.02,1.2) {$(l,r)=(3,5)$};
\end{axis}
\end{tikzpicture}
\caption{The period function $\Phi$ of \eqref{eq:rot-Phi}, computed numerically, for $n=3$ (solid) and $n=8$ (grey). In accordance with Lemma~\ref{lem:Phi-limits} and Theorem~\ref{thm:Phi-monotone}, $\Phi$ increases from $\pi$ to $\sqrt2\,\pi$. For $n=8$ the graph leaves the value $\pi$ with a vertical tangent, because of the contribution of order $E^{1/c}$, $c=2$, coming from the neighbourhood of the singular boundary in the proof of Lemma~\ref{lem:Phi-limits}; for $n=3$, where $c=1/3$, this contribution is $O(E)$. The dotted levels $\Phi=2\pi l/r$ correspond to $(l,r)=(2,3)$ and $(3,5)$: by monotonicity each of them is met exactly once, and the dots give the values of $E$ of the closed profiles shown in Figure~\ref{fig:rot-closed}.}
\label{fig:Phi}
\end{figure}

The graph of $\Phi$ for $n=3$ and $n=8$ is shown in Figure~\ref{fig:Phi}, together with the values of $E$ which give the closed profiles of Figure~\ref{fig:rot-closed}.

\begin{corollary}\label{cor:rot-classification}
Let $n\ge3$. A compact proper biconservative $SO(n-1)$-invariant hypersurface of $\Sph^n$ is, up to an isometry of $\Sph^n$ and up to composition with a finite covering of the factor $\Sph^1$ (that is, up to traversing its closed profile several times), one of the hypersurfaces of Theorem~\ref{thm:rot-closed}. For each pair of coprime natural numbers $l,r$ with $r<2l<\sqrt2\,r$ there is exactly one of them whose profile is traversed once, up to isometries, and none of them is embedded.
\end{corollary}

\begin{proof}
By {Theorem~\ref{thm:A}} the profile of a proper biconservative $SO(n-1)$-invariant hypersurface satisfies \eqref{eq:rot-reduced} along the part of it which lies in the interior of the orbit space, and its curvature is not constant. By Proposition~\ref{prop:rot-periodic} it therefore has $0<|E|<F(\rho_0)$ and stays in the regular part of the orbit space, so the hypersurface is compact if and only if its profile is closed, i.e., $\Phi(|E|)=2\pi l/r$; we may take $l,r$ coprime, the profile being then traversed a finite number of times. By Theorem~\ref{thm:Phi-monotone}, this happens if and only if $\pi<2\pi l/r<\sqrt2\,\pi$, and then $|E|$ is uniquely determined by $l/r$. Two profiles with the same $|E|$ differ by a rotation about the pole and possibly a change of orientation, so the corresponding hypersurfaces are congruent. Finally, $\pi<\Phi<2\pi$, so none of them is embedded by Theorem~\ref{thm:rot-closed}.
\end{proof}

\begin{proof}[{\bf Proof of Theorem~\ref{thm:F}}]
By {Theorem~\ref{thm:A}}, the profile of a proper biconservative $SO(n-1)$-invariant hypersurface satisfies \eqref{eq:rot-reduced} and has nonconstant $\kappa$, hence $E\neq0$ and $E\neq F(\rho_0)$ (Proposition~\ref{prop:rot-periodic}). Its profile is $\sigma$-elastic by Remark~\ref{cor:pelastic}, and the description of the profile and of its curvature is Proposition~\ref{prop:rot-periodic}. The limits of $\Phi$ are Lemma~\ref{lem:Phi-limits} and its monotonicity is Otsuki's Theorem~\ref{thm:Phi-monotone}; the existence{ and non-embeddedness of the compact examples are Theorem~\ref{thm:rot-closed}, and their uniqueness and the classification up to finite coverings are Corollary~\ref{cor:rot-classification}.}
\end{proof}

\begin{remark}[Constant scalar curvature]\label{rem:scalar}
The problem of determining the biconservative hypersurfaces with constant scalar curvature was raised by Fetcu and Oniciuc \cite{FetcuOniciuc}. Fu, Hong, Yang and Zhan \cite{FHYZ} proved that such hypersurfaces have constant mean curvature in four-dimensional space forms, and that in five-dimensional space forms they are either CMC or an open part of a certain rotational hypersurface. Our families illustrate this in a precise way. By the Gauss equation, the scalar curvature of a hypersurface of $\R^n$ ($\varepsilon=0$) or $\Sph^n$ ($\varepsilon=1$) is 
\[
\mathrm{scal}=\varepsilon(n-1)(n-2)+f^2-|A|^2,
\]
and along a proper biconservative profile $f=-2\kappa$ and $\sum_im_ik_i=-3\kappa$, so that
\[
f^2-|A|^2=3\kappa^2-\sum_{i=0}^{d-1}m_ik_i^2 .
\]
\begin{enumerate}[label=\textup{(\alph*)},leftmargin=2em]
\item For $d=1$ the principal curvature of the orbits is $k_0=-3\kappa/(n-2)$, so that
\[
f^2-|A|^2=\frac{3(n-5)}{n-2}\,\kappa^2 .
\]
Since $\kappa$ is not constant, the scalar curvature is constant if and only if $n=5$. In this case the principal curvatures are $\kappa,-\kappa,-\kappa,-\kappa$, and $\mathrm{scal}\equiv0$ in $\R^5$, $\mathrm{scal}\equiv12$ in $\Sph^5$. In $\R^5$ the first integral {\eqref{eq:d1-first-integral}} is $y\cos\alpha$, so that the profiles are catenaries, in accordance with Corollary~\ref{cor:tripling}: the hypersurfaces $\{(x,Y)\in\R\times\R^4:\ |Y|=a\cosh((x-b)/a)\}$ are scalar-flat proper biconservative hypersurfaces. In $\Sph^5$ we have $\mathsf n=2$, the profiles are those of the Otsuki tori{, i.e., of the compact rotational minimal surfaces of $\Sph^3$ which are neither totally geodesic nor Clifford tori; indeed, for $\mathsf n=2$ the reduced metric coincides with the Hsiang--Lawson metric of the rotational minimal surfaces of $\Sph^3$ (Corollary~\ref{cor:tripling})}, and Theorem~\ref{thm:F} provides infinitely many compact proper biconservative immersions $\Sph^1\times\Sph^3\to\Sph^5$ with constant scalar curvature $12$, none of them embedded.
\item For $d\ge2$, the hypersurfaces of Theorem~\ref{thm:B} and of Proposition~\ref{prop:boundary-profile} do not have constant scalar curvature. Indeed, in the polar variables of Section~\ref{sec:compactified}, writing $M=n-2$, $c_i=\cot(i\pi/d-\vartheta)$, $C_1=\sum_im_ic_i=-\ell'(\vartheta)$ and $C_2=\sum_im_ic_i^2$, one finds
\[
r^2\bigl(f^2-|A|^2\bigr)=\Bigl(\frac{C_1^2}3-C_2\Bigr)\cos^2\delta+2C_1\Bigl(\frac M3-1\Bigr)\sin\delta\cos\delta+M\Bigl(\frac M3-1\Bigr)\sin^2\delta ,
\]
which tends to $-C_2(\vartheta_0)$ as the profile approaches $P_0$ or $P_1$. Since $dc_i/d\vartheta=1+c_i^2$, we have $\ell''=-(n-2)-C_2$, so that, by Lemma~\ref{lem:ell2-universal-new}, $C_2(\vartheta_0)=(d-1)(n-2)$; thus
\[
r^2\bigl(f^2-|A|^2\bigr)\longrightarrow-(d-1)(n-2)<0 .
\] Hence $f^2-|A|^2$ is negative for large $r$, but tends to zero at infinity, so it is not constant. For $d=1$ one has $C_2(\vartheta_0)=0$, in agreement with~(a).
\end{enumerate}
\end{remark}

{\em Use of AI}.\quad {During the preparation of this work the authors used Claude (Anthropic), through the Claude Code environment, to check symbolic and numerical computations, to carry out the numerical integrations used for the figures, to verify bibliographic references, and to revise parts of the text. After using this tool, the authors reviewed and edited the content as needed and take full responsibility for the content of the article.}

\end{document}